\documentclass[12pt]{amsart} 
\usepackage{amssymb}
\usepackage{fullpage}
\usepackage[colorlinks=true]{hyperref}
\usepackage{graphicx}

\newcommand{\bbC}{\mathbb{C}}
\newcommand{\bbK}{\mathbb{K}}

\newcommand{\bbQ}{\mathbb{Q}}

\newcommand{\bbZ}{\mathbb{Z}}

\newcommand{\cC}{\mathcal{C}}
\newcommand{\cD}{\mathcal{D}}
\newcommand{\cN}{\mathcal{N}}
\newcommand{\cO}{\mathcal{O}}

\DeclareMathOperator{\core}{sf}
\DeclareMathOperator{\Imag}{Im}

\DeclareMathOperator{\rad}{rad}
\DeclareMathOperator{\Real}{Re}

\DeclareMathOperator{\Tr}{Tr}

\numberwithin{equation}{section}

\newtheorem{theorem}{Theorem}[section]
\newtheorem{conjecture}[theorem]{Conjecture}

\newtheorem{lemma}[theorem]{Lemma}
\newtheorem{proposition}[theorem]{Proposition}

\theoremstyle{definition}

\newtheorem*{remark-nonum}{Remark}

\begin{document}

\title[Bounds on the number of squares \ldots]{Bounds on the number of squares in recurrence sequences}

\author{Paul M Voutier}
\address{London, UK}
\email{Paul.Voutier@gmail.com}

\date{}

\begin{abstract}
We investigate the number of squares in a very broad family of binary recurrence
sequences with $u_{0}=1$. We show that there are at most two distinct squares in
such sequences (the best possible result), except under very special conditions
where we prove there are at most three such squares.
\end{abstract}


\maketitle

\section{Introduction}

\subsection{Background}

The study of the arithmetic properties of recurrence sequences has a long history.
Important problems include the open question of whether there are infinitely many
primes in such sequences, lower bounds for the largest prime divisor of the $n$-th
element \cite{St}, the zero-multiplicity of such sequences (how often zero occurs
as an element) \cite{Sch} and also the occurrence of powers in them \cite{BMS1}.

While some progress on arithmetic questions has been achieved for sequences of
Lucas numbers (numbers of the form $u_{n}=\left( \alpha^{n}-\beta^{n} \right) / (\alpha-\beta)$,
where $\alpha$ and $\beta$ are algebraic numbers such that $\alpha+\beta$ and
$\alpha \beta$ are non-zero coprime rational numbers and $\alpha/\beta$ is not
a root of unity) and their associated sequences ($v_{n}=\alpha^{n}+\beta^{n}$),
much less is known for other binary recurrence sequences.

As with sequences of Lucas numbers, other binary recurrence sequences are
connected with the solutions of Diophantine equations such as
$aX^{2}-bY^{4}=c$. Such equations are important in their own right, and also
as quartic models of elliptic curves.

This paper is concerned with such problems, in particular how many distinct
squares can occur in non-degenerate binary recurrence sequences.

\subsection{Notation}
\label{subsect:notation}

To define the sequences we are interested in and to express our results, we
start with some notation.

Let $a$, $b$ and $d$ be positive integers such
that $d$ is not a square. Suppose $\alpha=a+b^{2}\sqrt{d}$ has $N_{\alpha}
=N_{\bbQ \left( \sqrt{d} \right)/\bbQ}(\alpha)=a^{2}-b^{4}d$ and let
$\varepsilon=\left( t+u\sqrt{d} \right)/2$ be a unit in
$\cO_{\bbQ \left( \sqrt{d} \right)}$ with $t$ and $u$ positive integers.

We define the two sequences $\left( x_{k} \right)_{k=-\infty}^{\infty}$
and $\left( y_{k} \right)_{k=-\infty}^{\infty}$ by
\begin{equation}
\label{eq:yk-defn}
x_{k}+y_{k}\sqrt{d}
=\alpha \varepsilon^{2k}.
\end{equation}

Observe that $x_{0}=a$,
$y_{0}=b^{2}$,
\begin{equation}
\label{eq:yPM1}
y_{1}= \frac{b^{2}\left( t^{2}+du^{2} \right)+2atu}{4}, \hspace*{3,0mm}
y_{-1}=\frac{b^{2}\left( t^{2}+du^{2} \right)-2atu}{4}
\end{equation}
and that both sequences satisfy
the recurrence relation
\begin{equation}
\label{eq:yk-recurrence}
u_{k+1}= \frac{t^{2}+du^{2}}{2} u_{k}-u_{k-1},
\end{equation}
for all $k \in \bbZ$. Note that $\left( t^{2}+du^{2} \right)/2
=\Tr_{\bbQ \left( \sqrt{d} \right)/\bbQ} \left( \varepsilon^{2} \right)$.

To relate such sequences back to the quartic equations mentioned in the previous
section, observe that from \eqref{eq:yk-defn},
\[
x_{k}^{2}-dy_{k}^{2}=N_{\alpha}.
\]

We restrict to the coefficient of $\sqrt{d}$ in $\alpha$ being a square since
here we are interested in squares in the sequence of $y_{k}$'s. For such problems,
we will choose $\alpha$ such that $b^{2}$ is the smallest square among the $y_{k}$'s.

For any non-zero integer, $n$, let $\core(n)$ be the unique squarefree integer
such that $n/\core(n)$ is a square. We will put $\core(1)=1$.

\subsection{Conjectures}
\label{subsect:conj}

We start with some conjectures regarding squares in the sequence of $y_{k}$'s.
The dependence on the arithmetic of $N_{\alpha}$ is noteworthy.

\begin{conjecture}
\label{conj:1-seq}
There are at most four distinct integer squares among the $y_{k}$'s.

If $\core \left( \left| N_{\alpha} \right| \right) =2^{\ell}p^{m}$ where
$\ell, m \in \{ 0,1 \}$ with $\ell+m \geq 1$ and $p$ is an odd prime, then there
are at most three distinct integer squares among the $y_{k}$'s.

Furthermore, if $\left| N_{\alpha} \right|$ is a perfect square, then there are
at most two distinct integer squares among the $y_{k}$'s.
\end{conjecture}

When $N_{\alpha}<0$ and $b=1$, a stronger result appears to hold.

\begin{conjecture}
\label{conj:3-seq}
Suppose that $b=1$ and $N_{\alpha}<0$.
There are at most three distinct integer squares among the $y_{k}$'s.

If $-N_{\alpha}$ is a perfect square, then there are at most two distinct
integer squares among the $y_{k}$'s.
\end{conjecture}

In fact, a more general result than Conjecture~\ref{conj:1-seq} also appears to be true.

\begin{conjecture}
\label{conj:2-seq}
Let $\left( y_{k} \right)_{k=-\infty}^{\infty}$ be the sequence formed by replacing
$\varepsilon^{2k}$ in \eqref{eq:yk-defn} with $\varepsilon^{k}$.

There are at most four distinct integer squares among the $y_{k}$'s.

If $\left| N_{\alpha} \right|$ is a prime power or a perfect square,
then there are at most three distinct integer squares among the $y_{k}$'s.
\end{conjecture}

Examples showing that if these conjectures are true, then they are best possible
are provided in Subsections~\ref{subsect:conj1-egs}, \ref{subsect:conj3-egs} and
\ref{subsect:conj2-egs}.

\vspace*{1.0mm}

We used PARI/GP \cite{Pari} to search for squares among $y_{-80},\ldots, y_{80}$
from Conjecture~\ref{conj:2-seq}, where $\varepsilon=\left( t+u\sqrt{d} \right)/2>1$
is the fundamental unit in $\cO_{\bbQ \left( \sqrt{d} \right)}$ with $t, u \in \bbZ$ and
$u \leq 10^{30}$. Note that this also includes searching among $y_{-40},\ldots, y_{40}$
from Conjectures~\ref{conj:1-seq} and \ref{conj:3-seq}.
The range of $b$ and $d$ in this search was $1 \leq b \leq 2000$ and $ 1\leq d \leq 1000$
with $d$
squarefree. For $a$, we searched over two ranges:
(1) $\max \left( 1, \lfloor \sqrt{db^{4}} \rfloor -1000 \right)
\leq a \leq \lceil \sqrt{db^{4}} \rceil + 1000$ (so $\left| N_{\alpha} \right|$ is small)
and (2) $1 \leq a \leq 2000$.
For both ranges of $a$, we considered only
$\gcd \left( a,b^{2} \right)$ squarefree. The search across these ranges took just
under $319$ hours on a Windows 11 laptop with an Intel i7-13700H 2.40 GHz processor
and 32~GB of RAM.
No counterexamples to any of these conjectures were found in this search.

\begin{remark-nonum}
The distinctness condition in these conjectures, and in our results below, is
important, as such sequences can have repeated elements. E.g., $(a,b,d,t,u)=(42,4,7,16,6)$
where $y_{-k}=y_{k-1}$ for all $k \geq 1$.
But this can only happen when $\alpha$ divided by its algebraic conjugate is a
unit in the ring of integers.
\end{remark-nonum}

\subsection{Results}
\label{subsect:results}

In this paper, we prove Conjecture~\ref{conj:3-seq} when $-N_{\alpha}$ is a square,
except for a very limited set of sequences, where we prove a slightly weaker result.

\begin{theorem}
\label{thm:1.3-seq-new}
Let $b=1$, $a$ and $d$ be positive integers, where $d$ is not a square, $N_{\alpha}<0$
and $-N_{\alpha}$ is a square.

\noindent
{\rm (a)} If $u=1$, $t^{2}-du^{2}=-4$, $N_{\alpha} \equiv 12 \pmod{16}$, $\gcd \left( a^{2},d \right)=1,4$
and one of $y_{\pm 1}$ is a perfect square, then there are at most three
distinct squares among the $y_{k}$'s.

\noindent
{\rm (b)} If $u=2$, $t^{2}-du^{2}=-4$, $N_{\alpha}$ is odd, $\gcd \left( a^{2},d \right)=1$
and one of $y_{\pm 1}$ is a perfect square, then there are at most three
distinct squares among the $y_{k}$'s.

\noindent
{\rm (c)} Otherwise, there are at most two distinct squares among the $y_{k}$'s.
\end{theorem}

In Subsection~\ref{subsect:thm14a-egs}, we present an infinite family of examples
of sequences satisfying the conditions in part~(a) with $y_{1}$ roughly one-tenth
the size of the bounds in Proposition~\ref{prop:4.1}. The same is possible for
$y_{-1}$ in part~(a) and $y_{\pm 1}$ in part~(b).

\subsection{Our method of proof}

Our proof has its basis in the work of Siegel \cite{Siegel2} (also see
\cite{Evert1}). Foremost with our use of hypergeometric functions. But also with
how we use them here. The definition and use of our quantity $r_{0}$
in Section~\ref{sect:prop-11} has similarities to the definition and use of
$\ell_{1}$ and $\ell_{2}$ from Lemma~7 onwards in \cite{Evert1}. The idea is to
treat ranges of hypothetical squares beyond those that can be treated by our
gap principle alone.

We state and prove our results in Section~\ref{sect:prelim} more
generally than required here. This is because they will be useful in further
work on binary recurrence sequences. An example is \cite{V6} where we prove
Conjecture~\ref{conj:3-seq} when $N_{\alpha}$ is a prime power.

\subsection{Structure of this article}

Section~\ref{sect:dio} contains results on
diophantine approximation, with a focus on the use of hypergeometric functions
to do so.
In Section~\ref{sect:prelim}, we collect various facts that we will require
about elements of the sequence $\left( y_{k} \right)_{k=-\infty}^{\infty}$
defined by \eqref{eq:yk-defn}. These
two sections are independent of each other. Their results are brought together
in Section~\ref{sect:prop-11}, where we state and prove Proposition~\ref{prop:4.1}.
Our theorem follows from this proposition. Its proof is given in Section~\ref{sect:proofs}.
Finally, in Section~\ref{sect:egs},
we give examples showing that our results and conjectures are best possible.

\subsection{Acknowledgements}

The author deeply appreciates all the time and effort spent by Mihai Cipu,
very carefully reading an earlier version of this work. His extensive comments,
questions and our follow-up discussions led to significant improvements in the
presentation of this work, as well as numerous corrections. The author is also
grateful to the referee for their notes and suggestions, which improved this paper
too.

\section{Diophantine Approximation via Hypergeometric Functions}
\label{sect:dio}

Recall that by an {\it effective irrational measure} for an irrational number,
$\alpha$, we mean an inequality of the form
\[
\left| \alpha - \frac{p}{q} \right|>\frac{c}{|q|^{\mu}},
\]
for all $p/q \in \bbQ$ with $\gcd(p,q)=1$ and $|q|>Q$, where $c$, $Q$ and $\mu$
are all effectively computable.

By Liouville's famous result \cite{Liou}, where he constructed the first
examples of numbers proven to be transcendental, we have such effective irrational
measures for algebraic numbers of degree $n$, with $\mu=n$. But for practically
all applications we require $\mu<n$.

We can sometimes use the hypergeometric method to obtain effective irrationality measures
that improve on Liouville's result for the algebraic numbers that arise
here. This was first done by Baker \cite{Baker1, Baker2} building on earlier
applications of hypergeometric functions to diophantine approximation such as
those of Siegel \cite{Siegel1}. However, that does not suffice
for us to prove our results here. The problem here arises not because of the exponent,
$\mu$, in the effective irrationality measure, but because the constant, $c$,
is too small. Upon investigating this further, we found that we can complete
the proof of our results if we use not the effective irrationality measures from
the hypergeometric method, but rather consider more carefully the actual results
that we obtain from the use of hypergeometric functions.

The means of doing so is the following lemma.

\begin{lemma}
\label{lem:2.1}
Let $\theta \in \bbC$ and let $\bbK$ be an imaginary quadratic field. Suppose that
there exist $k_{0},\ell_{0} > 0$ and $E,Q > 1$ such that for all non-negative
integers $r$, there are algebraic integers $p_{r}$ and $q_{r}$ in $\bbK$ with
$\left| q_{r} \right| < k_{0}Q^{r}$ and
$\left| q_{r} \theta - p_{r} \right| \leq \ell_{0}E^{-r}$ satisfying 
$p_{r}q_{r+1} \neq p_{r+1}q_{r}$.

For any algebraic integers $p$ and $q$ in $\bbK$, let $r_{0}$ be the smallest
positive integer such that $\left( Q-1/E  \right)\ell_{0}|q|/\left( Q-1 \right)<cE^{r_{0}}$,
where $0<c<1$.

\noindent
{\rm (a)} We have 
\[
\left| q\theta - p \right| 
> \frac{1-c/E}{k_{0}Q^{r_{0}+1}}.
\]

\noindent
{\rm (b)} When $p/q \neq p_{r_{0}}/q_{r_{0}}$, we have
\[
\left| q\theta - p \right| 
> \frac{1-c}{k_{0}Q^{r_{0}}}.
\]
\end{lemma}

\begin{remark-nonum}
This is a modification of Lemma~6.1 of \cite{V2} (and other similar results).

First, we replace $1/\left( 2k_{0} \right)$ in both parts and define $r_{0}$
somewhat differently too.

Second, the terms in the lower bounds in both parts are no longer converted into
ones involving $|q|^{-(\kappa+1)}$. For our use here, where an effective irrationality
measure is not required, this results in further improvements.
\end{remark-nonum}

\begin{proof}
Let $p$, $q$ be algebraic integers in $\bbK$.
Let $R=q\theta-p$ and $R_{r}=q_{r}\theta-p_{r}$. We can write
\[
q_{r}R=q_{r}q\theta-q_{r}p=q \left( q_{r}\theta-p_{r} \right)+p_{r}q-q_{r}p
=qR_{r} + p_{r}q-q_{r}p.
\]

We consider two cases according to whether $p_{r_{0}}q-q_{r_{0}}p=0$ or not.

If $p_{r}q-q_{r}p \neq 0$, then $\left| p_{r}q-q_{r}p \right| \geq 1$, since it
is an algebraic integer in an imaginary quadratic field. From our upper bounds
for $\left| q_{r} \right|$ and $\left| q_{r}\theta-p_{r} \right|$ in the
statement of the lemma, we have
\[
k_{0}Q^{r}|R| > 1 - \ell_{0}E^{-r}|q|.
\]

With $r=r_{0}$, from the definition of $r_{0}$, we have
\[
k_{0}Q^{r}|R| > 1 - c\frac{Q-1}{Q-1/E} > 1 - c.
\]

Part~(b) now follows from the upper bound for $\left| q_{r_{0}} \right|$ in the
statement of the lemma.

If $p_{r_{0}}q-q_{r_{0}}p=0$, then we use $r=r_{0}+1$. From the definition of
$r_{0}$, we have
\[
k_{0}Q^{r_{0}+1}|R| > 1 - c\frac{Q-1}{E(Q-1/E)}=\frac{QE-1-cQ+c}{QE-1}
>1-c/E.
\]

Part~(a) now follows.
\end{proof}

\subsection{Construction of Approximations}
\label{subsect:const}

Let $t'$, $u_{1}$ and $u_{2}$ be rational integers with $t'<0$
\footnote{To avoid confusion with $t$ in the expression for $\varepsilon$, we
use $t'$ here to denote what was $t$ in our previous works like \cite{V3}. Similarly,
we will use $d'$ where $d$ was used in previous works.} such that
$u=\left( u_{1}+u_{2}\sqrt{t'} \right)/2$ be an algebraic integer in
$\bbK=\bbQ \left( \sqrt{t'} \right)$ with $\sigma(u)=\left( u_{1}-u_{2}\sqrt{t'} \right)/2$
as its algebraic (and complex) conjugate. Put $\omega = u/\sigma(u)$ and write
$\omega=e^{i\varphi}$, where $-\pi<\varphi \leq \pi$. For any real number
$\nu$, we shall put $\omega^{\nu}= e^{i\nu\varphi}$ -- unless otherwise stated,
we will use this convention throughout this paper.

Suppose that $\alpha$, $\beta$ and $\gamma$ are complex numbers and $\gamma$ is
not a non-positive integer. We denote by ${}_{2}F_{1}(\alpha, \beta, \gamma, z)$
the classical (or Gauss) hypergeometric function of the complex variable $z$.

For integers $m$ and $n$ with $0 < m < n$, $(m,n) = 1$ and $r$ a non-negative
integer, put $\nu=m/n$ and
\[
X_{m,n,r}(z)={}_{2}F_{1}(-r-\nu, -r, 1-\nu, z), \quad
Y_{m,n,r}(z)=z^{r}X_{m,n,r} \left(z^{-1} \right)
\]
and
\[
R_{m,n,r}(z)
= (z-1)^{2r+1} \frac{\nu \cdots (r+\nu)}{(r+1) \cdots (2r+1)} 
		   {} _{2}F_{1} \left( r+1-\nu, r+1; 2r+2; 1-z \right).
\]

We collect here some facts about these functions that we will require.

\begin{lemma}
\label{lem:hypg}
\noindent
{\rm (a)} Suppose that $|\omega-1|<1$. We have
\[
\omega^{\nu}Y_{m,n,r}(\omega)-X_{m,n,r}(\omega)=R_{m,n,r}(\omega).
\]

\noindent
{\rm (b)} We have
\[
X_{m,n,r}(\omega)Y_{m,n,r+1}(\omega) \neq X_{m,n,r+1}(\omega)Y_{m,n,r}(\omega).
\]

\noindent
{\rm (c)} If $|\omega|=1$ and $|\omega-1|<1$, then
\[
\left| R_{m,n,r}(\omega) \right|
\leq \frac{\Gamma(r+1+\nu)}{r!\Gamma(\nu)} |\varphi| \left| 1-\sqrt{\omega} \right|^{2r}.
\]

\noindent
{\rm (d)} If $|\omega|=1$ and $|\omega-1|<1$, then
\[
\left| X_{m,n,r}(\omega) \right| = \left| Y_{m,n,r}(\omega) \right|
< 1.072\frac{r!\Gamma(1-\nu)}{\Gamma(r+1-\nu)} \left| 1 + \sqrt{\omega} \right|^{2r}.
\]

\noindent
{\rm (e)} For $|\omega|=1$ and $\Real(\omega) \geq 0$, we have
\[
\left| {} _{2}F_{1} \left( r+1-\nu, r+1; 2r+2; 1-\omega \right) \right| \geq 1,
\]
with the minimum value occurring at $\omega=1$.
\end{lemma}

\begin{proof}
Part~(a) is established in the proof of Lemma~2.3 of \cite{Chen1}.

Part~(b) is Lemma~4 of \cite{Baker2}.

Part~(c) is Lemma~2.5 of \cite{Chen1}.

Part~(d) is Lemma~4 of \cite{V4}.

Part~(e) is Lemma~5 of \cite{V4}.
\end{proof}

We let $D_{n,r}$ denote the smallest positive integer such that $D_{n,r} X_{m,n,r}(x) \in \bbZ[x]$
for all $m$ as above. For $d' \in \bbZ$, we define $N_{d',n,r}$ to be the largest
integer such that $\left( D_{n,r}/ N_{d',n,r} \right)X_{m,n,r}\left( 1-\sqrt{d'}\,x \right)
\in \bbZ \left[ \sqrt{d'} \right] [x]$, again for all $m$ as above. We will use
$v_{p}(x)$ to denote the largest power of a prime $p$ which divides
the rational number $x$. We put
\begin{equation}
\label{eq:ndn-defn}
\cN_{d',n} =\prod_{p|n} p^{\min(v_{p}(d')/2, v_{p}(n)+1/(p-1))}.
\end{equation}

In what follows, we shall restrict our attention to $m=1$ and $n=4$, so $\nu=1/4$.

\begin{lemma}
\label{lem:denom-est}
\noindent
{\rm (a)} We have
\begin{equation}
\label{eq:cndn-defn}
\frac{\Gamma(3/4) \, r!}{\Gamma(r+3/4)} \frac{D_{4,r}}{N_{d',4,r}}
<\cC_{4,1} \left( \frac{\cD_{4}}{\cN_{d',4}} \right)^{r}
\text{ and } \hspace*{1.0mm}
\frac{\Gamma(r+5/4)}{\Gamma(1/4)r!} \frac{D_{4,r}}{N_{d',4,r}}
< \cC_{4,2} \left( \frac{\cD_{4}}{\cN_{d',4}} \right)^{r}
\end{equation}
for all non-negative integers $r$, where $\cC_{4,1}=0.83$, $\cC_{4,2}=0.2$
and $\cD_{4}=e^{1.68}$.

\noindent
{\rm (b)} For any positive integer, $r$, we have
\begin{equation}
\label{eq:gamma-bnds}
\frac{5}{24 \cdot 4^{r} r^{1/4}} \leq \frac{(1/4) \cdots (r+1/4)}{(r+1) \cdots (2r+1)}
\quad \text{and} \quad
\frac{r!\Gamma(3/4)}{\Gamma(r+3/4)} \leq 4r^{1/4}/3.
\end{equation}
\end{lemma}

\begin{proof}
(a) From Lemma~7.4(c) of \cite{V2}, we have
\[
\max \left( 1, \frac{\Gamma(3/4) \, r!}{\Gamma(r+3/4)},
4\frac{\Gamma(r+5/4)}{\Gamma(1/4)r!} \right)
\frac{D_{4,r}}{N_{d',4,r}} < 100 \left( \frac{e^{1.64}}{\cN_{d',4}} \right)^{r}
\]

However, the value $100$ in this inequality results in us requiring a lot of
computation to complete the proof of our results here (in particular,
Proposition~\ref{prop:4.1}). Therefore, we seek a
smaller value at the expense of replacing $1.64$ by a larger value, whose value
has less of an impact on our proof. For $r \geq 156$, we have $100\exp(1.64r)<0.2\exp(1.68r)$,
so we compute directly the left-hand sides of \eqref{eq:cndn-defn} for $r \leq 155$.
We find that the maximum values of the left-hand sides of \eqref{eq:cndn-defn}
divided by $\exp(1.68r)$ both occur for $r=3$.
Part~(a) follows.

(b) We prove both of these by induction. A quick calculation
shows that we have equality in both cases for $r=1$. Then we take the
expression for $r+1$ and divide it by the expression for $r$. For
the first inequality, this ratio is $(1/4)\left( r+5/4 \right)/\left( r+3/2 \right)$,
so the inequality will hold if $\left( r+5/4 \right)/\left( r+3/2 \right)>\left( r/\left( r+1 \right) \right)^{1/4}$.
We take the fourth-power of both sides and subtract them, this gives us
\[
\frac{224r^{3}+944r^{2}+1329r+625}{16\left( 2r+3 \right)^{4} \left( r+1 \right)},
\]
which is clearly positive for $r \geq 1$. Since both $\left( r+5/4 \right)/\left( r+3/2 \right)$
and $\left( r/\left( r+1 \right) \right)^{1/4}$ are positive real numbers
for $r \geq 1$, the required inequality hlds.

We proceed in the very same way to prove the second inequality.
\end{proof}

As in \cite[Theorem~1]{V3}, put
\begin{align}
\label{eq:g-defn}
g_{1}  & = \gcd \left( u_{1}, u_{2} \right), \\
g_{2}  & = \gcd \left( u_{1}/g_{1}, t' \right), \nonumber \\
g_{3}  & = \left\{
			 \begin{array}{ll}
	             1 & \text{if $t' \equiv 1 \pmod{4}$ and $\left( u_{1}-u_{2} \right)/g_{1} \equiv 0 \pmod{2}$}, \\
	             2 & \text{if $t' \equiv 3 \pmod{4}$ and $\left( u_{1}-u_{2} \right)/g_{1} \equiv 0 \pmod{2}$},\\
	             4 & \text{otherwise,}
             \end{array}
             \right. \nonumber \\
g      & = g_{1}\sqrt{g_{2}/g_{3}}. \nonumber
\end{align}

Then we can put
\begin{align}
\label{eq:7}
p_{r} &= \frac{D_{4,r}}{N_{d',4,r}} \left( \frac{u_{1}-u_{2}\sqrt{t'}}{2g} \right)^{r} X_{1,4,r}(\omega),\\
q_{r} &= \frac{D_{4,r}}{N_{d',4,r}} \left( \frac{u_{1}-u_{2}\sqrt{t'}}{2g} \right)^{r} Y_{1,4,r}(\omega)
\hspace*{3.0mm} \text{ and} \nonumber \\
R_{r} &= \frac{D_{4,r}}{N_{d',4,r}} \left( \frac{u_{1}-u_{2}\sqrt{t'}}{2g} \right)^{r} R_{1,4,r}(\omega),
\nonumber
\end{align}
where
\begin{equation}
\label{eq:d-defn}
d'=\left( u-\sigma(u) \right)^{2}/g^{2} = u_{2}^{2}t'/g^{2}.
\end{equation}

From Lemma~\ref{lem:hypg}(a), we have
\[
q_{r}\omega^{1/4} - p_{r}=R_{r}.
\]

By the definitions of $D_{4,r}$ and $N_{d',4,r}$, we see that $X_{1,4,r}(1-z)$
is a polynomial of degree $r$ with coefficients in $\bbZ \left[ \sqrt{d'} \right]$ and
\begin{align*}
p_{r} &= \frac{D_{4,r}}{N_{d',4,r}} \left( \frac{u_{1}-u_{2}\sqrt{t'}}{2g} \right)^{r} X_{1,4,r}(\omega) \\
&= \frac{D_{4,r}}{N_{d',4,r}} \left( \frac{u_{1}-u_{2}\sqrt{t'}}{2g} \right)^{r} X_{1,4,r} \left( 1-u_{2}\sqrt{t'} \frac{-2}{u_{1}-u_{2}\sqrt{t'}} \right).
\end{align*}

Since $d'=u_{2}^{2}t'/g^{2}$ and $\left( u_{1}-u_{2}\sqrt{t'} \right)/(2g)$ is an
algebraic integer, it follows that $p_{r}$ is an algebraic integer in $\bbQ \left( \sqrt{t'} \right)$.
The analogous expression for $q_{r}$ shows that it is also an algebraic integer.

So by applying these quantities, with Lemma~\ref{lem:hypg}(d) and Lemma~\ref{lem:denom-est}(a),
in Lemma~\ref{lem:2.1}, we can take
\begin{equation}
\label{eq:q-defn}
Q = \frac{\cD_{4} \left| \left| u_{1} \right| + \sqrt{u_{1}^{2}-t'u_{2}^{2}} \right|}{|g|\cN_{d',4}}
\end{equation}
and
\begin{equation}
\label{eq:k-UB}
k_{0}<1.072\cC_{4,1}<0.89.
\end{equation}

Using Lemma~\ref{lem:hypg}(c) instead of Lemma~\ref{lem:hypg}(d), we also have
\begin{equation}
\label{eq:e-defn}
E = \frac{|g|\cN_{d',4} \left|  \left| u_{1} \right|  + \sqrt{u_{1}^{2}-t'u_{2}^{2}} \right|}{\cD_{4}u_{2}^{2}|t'|}
\end{equation}
and
\begin{equation}
\label{eq:ell0-defn}
\ell_{0}=\cC_{4,2}|\varphi|=0.2|\varphi|.
\end{equation}

\section{Lemmas about $\left( x_{k} \right)_{k=-\infty}^{\infty}$ and $\left( y_{k} \right)_{k=-\infty}^{\infty}$}
\label{sect:prelim}

\subsection{Representation Proposition}

The proposition in this section plays a crucial role in this work. It is the
expression for $f^{2} \left( x + N_{\varepsilon}\sqrt{N_{\alpha}} \right)$ in
this proposition that permits us to use the hypergeometric method described in the
previous section. For this reason, we call it our representation proposition.

It is also because we have $\varepsilon^{2}$, rather than $\varepsilon$,
in \eqref{eq:quad-rep-assumption} that we use $\varepsilon^{2k}$ in
\eqref{eq:yk-defn}.

For any non-zero integer, $n$, we let $\rad(n)$ be the product of all distinct
prime divisors of $n$. We will put $\rad(\pm 1)=1$.

\begin{proposition}
\label{prop:quad-rep}
Let $a \neq 0$, $b>0$ and $d$ be rational integers such that $d$ is not a square.
Put $\alpha=a+b^{2} \sqrt{d}$ and denote $N_{\bbQ(\sqrt{d})/\bbQ}(\alpha)$
by $N_{\alpha}$. Suppose that $N_{\alpha}$ is not a square, $x \neq 0$ and $y>0$
are rational integers with
\begin{equation}
\label{eq:quad-rep-assumption}
x+y^{2} \sqrt{d} = \alpha \varepsilon^{2},
\end{equation}
where $\varepsilon=\left( t+u\sqrt{d} \right)/2 \in \cO_{\bbQ \left( \sqrt{d} \right)}$
with $t$ and $u$ non-zero rational integers, norm $N_{\varepsilon}$.

\noindent
{\rm (a)} We can write
\begin{align}
f^{2} \left( x + N_{\varepsilon}\sqrt{N_{\alpha}} \right)
&= \left( a+\sqrt{N_{\alpha}} \right) \left( r + s\sqrt{\core \left( N_{\alpha} \right)} \right)^{4}
\quad \text{and} \nonumber \\
\label{eq:fy-rel}
fy &= b \left( r^{2}-\core \left( N_{\alpha} \right)s^{2} \right),
\end{align}
for some integers $f$, $r$ and $s$ satisfying $f \neq 0$,
\[
f | \left( 4b^{2} \rad \left( f' \gcd \left( uN_{\alpha}/\core \left( N_{\alpha} \right), N_{\varepsilon} \right) \right) \right),
\]
where $f'| \core \left( N_{\alpha} \right)$ and
$0<f'< \max \left( 2, \sqrt{\left| \core \left( N_{\alpha} \right) \right|} \right)$.

\noindent
{\rm (b)}
If $-N_{\alpha}$ is a square, then $f | \left( b^{2} \rad \left( \gcd \left( uN_{\alpha}, N_{\varepsilon} \right) \right) \right)$,
where the first relationship in part~{\rm (a)} is replaced by
\[
\pm f^{2} \left( x + N_{\varepsilon}\sqrt{N_{\alpha}} \right)
= \left( a+\sqrt{N_{\alpha}} \right) \left( r + s\sqrt{\core \left( N_{\alpha} \right)} \right)^{4}.
\]

\noindent
{\rm (c)}
If $\core \left( \left| N_{\alpha} \right| \right) = 2^{\ell}p^{m}$ where $\ell, m \in \{ 0, 1 \}$
with $\ell+m \geq 1$
and $p$ is an odd prime, then we have $f \mid \left( 4b^{2} \rad \left( \gcd \left( uN_{\alpha}/\core \left( N_{\alpha} \right), N_{\varepsilon} \right) \right) \right)$
when $N_{\alpha} \equiv 1 \pmod{4}$ and $4|d$, and $f \mid \left( 2b^{2} \rad \left( \gcd \left( uN_{\alpha}/\core \left( N_{\alpha} \right), N_{\varepsilon} \right) \right) \right)$
otherwise.
\end{proposition}

\begin{remark-nonum}
For work with binary recurrence sequences, $\varepsilon$ will be a power of an
algebraic integer $\varepsilon'$. So $\rad \left( N_{\varepsilon} \right)
=\rad \left( N_{\varepsilon'} \right)$. This makes our representation proposition
useful for more general sequences than those considered in this paper.
\end{remark-nonum}

To prove this proposition, we will use the following three lemmas.

\begin{lemma}
\label{lem:quad-rep-1}
Suppose that $a, b, d, t, u, x, y, \alpha$ and $\varepsilon$ are as in
Proposition~\ref{prop:quad-rep}.
Put $r_{1}=tb^{2}+au \pm 2by$ and $s_{1}'=-u\sqrt{N_{\alpha}/\core \left( N_{\alpha} \right)}$.

Then $\gcd \left( 4b^{2}r_{1}/\core \left( r_{1} \right), r_{1}^{2} \right)$
is a divisor of $s_{1}'^{2}$.
\end{lemma}

\begin{proof}
Put $v_{p}(b)=k_{1}$ and $v_{p} \left( r_{1} \right)=k_{2}$ with $k_{1}, k_{2} \geq 0$.

Suppose first that $p$ is an odd prime.

With $g_{1}^{2}=\gcd \left( 4b^{2}r_{1}/\core \left( r_{1} \right), r_{1}^{2} \right)$,
we have
\begin{equation}
\label{eq:val-32-pOdd}
v_{p} \left( g_{1}^{2} \right)
= \left\{
\begin{array}{ll}
2k_{1}+k_{2},
& \text{if $k_{2}$ is even and $2k_{1} \leq k_{2}$,} \\
2k_{2},
& \text{if $k_{2}$ is even and $2k_{1} \geq k_{2}$,} \\
2k_{1}+k_{2}-1,
& \text{if $k_{2}$ is odd and $2k_{1}-1 \leq k_{2}$,} \\
2k_{2},
& \text{if $k_{2}$ is odd and $2k_{1}+1 \geq k_{2}$.}
\end{array}
\right.
\end{equation}

From the definitions of $r_{1}$ and $k_{2}$, it follows that
$v_{p}(au) \geq \min \left( k_{2}, v_{p} \left( tb^{2} \pm 2by \right) \right)$.

Expanding the expression for $y^{2}$ in \eqref{eq:quad-rep-assumption},
we obtain
\begin{equation}
\label{eq:rel-y}
4y^{2}=2atu+b^{2}t^{2}+b^{2}du^{2}=2atu+4b^{2}N_{\varepsilon}+2b^{2}du^{2}.
\end{equation}

Applying this equation, for odd primes, $p$, we have
\begin{align*}
2v_{p}(y) & \geq \min \left( v_{p}(atu), v_{p} \left( b^{2}t^{2} \right), v_{p} \left( b^{2}du^{2} \right) \right) \\
& \geq \min \left( v_{p} \left( t^{2}b^{2} \pm 2bty \right), k_{2}+v_{p}(t), v_{p} \left( b^{2}t^{2} \right), v_{p} \left( b^{2}du^{2} \right) \right) \\
& \geq \min \left( v_{p} \left( bty \right), k_{2}+v_{p}(t), v_{p} \left( b^{2}t^{2} \right), v_{p} \left( b^{2}du^{2} \right) \right) \\
& \geq \min \left( k_{1}+v_{p}(t)+v_{p}(y), k_{2}+v_{p}(t), 2k_{1}+2v_{p}(t), 2k_{1}+v_{p}(d)+2v_{p}(u) \right).
\end{align*}

This implies that
$v_{p}(y) \geq \min \left( k_{1}, k_{2}/2 \right)$.

We consider separately the two possibilities for the gcd defining $g_{1}^{2}$
for each odd prime, $p$.

(1) Suppose that $v_{p}\left( r_{1}^{2} \right)
\leq v_{p} \left( 4b^{2}r_{1}/\core \left( r_{1} \right) \right)
=v_{p} \left( b^{2}r_{1}/\core \left( r_{1} \right) \right)$.

If $k_{2}$ is even, then we have $v_{p}\left( r_{1}^{2} \right)
\leq v_{p} \left( b^{2}r_{1} \right)$. I.e., $k_{2} \leq 2k_{1}$.
Similarly, if $k_{2}$ is odd, then we have $v_{p}\left( r_{1}^{2} \right)
\leq v_{p} \left( b^{2}r_{1} \right)-1$. So $k_{2} \leq 2k_{1}-1$.
In both cases, we have $k_{1} \geq k_{2}/2$ and thus $v_{p}(y) \geq k_{2}/2$.

(1-a) If $v_{p}(a)<2v_{p}(b)+v_{p}(d)/2$, then $v_{p} \left( N_{\alpha} \right)
=v_{p} \left( N_{\alpha}/\core \left( N_{\alpha} \right) \right)=2v_{p}(a)$.
So $v_{p} \left( s_{1}'^{2} \right)=2v_{p}(au)$.

From $v_{p}(y) \geq k_{2}/2$ and $k_{1} \geq k_{2}/2$, it follows that
$v_{p} \left( tb^{2} \pm 2by \right) \geq k_{2}$. So from the definition of
$r_{1}$, we have $v_{p}(au) \geq k_{2}$. Hence
$v_{p} \left( s_{1}'^{2} \right) \geq 2k_{2}$. From \eqref{eq:val-32-pOdd},
we have $v_{p} \left( g_{1}^{2} \right)=2k_{2} \leq v_{p} \left( s_{1}'^{2} \right)$,
as desired.

(1-b) If $v_{p}(a) \geq 2v_{p}(b)+v_{p}(d)/2$, then $v_{p} \left( N_{\alpha} \right) \geq v_{p} \left( N_{\alpha}/\core \left( N_{\alpha} \right) \right)
\geq 4v_{p}(b)$. So $v_{p}\left( s_{1}'^{2} \right) \geq 4v_{p}(b)+2v_{p}(u)$.

From \eqref{eq:val-32-pOdd}, we have $v_{p} \left( g_{1}^{2} \right)=2k_{2}$.
Also $v_{p}\left( s_{1}'^{2} \right) \geq 4v_{p}(b)+2v_{p}(u) \geq 2k_{2}+2v_{p}(u)$.
So $v_{p} \left( g_{1}^{2} \right) \leq v_{p}\left( s_{1}'^{2} \right)$.

\vspace*{1.0mm}

(2) Now suppose that $v_{p} \left( 4b^{2}r_{1}/\core \left( r_{1} \right) \right)
< v_{p}\left( r_{1}^{2} \right)$.

Put $r_{1}'=tb^{2}+au+2by$ and $r_{1}''=tb^{2}+au-2by$. We shall prove the lemma
in this case (i.e., when $v_{p} \left( 4b^{2}r_{1}/\core \left( r_{1} \right) \right)
< v_{p}\left( r_{1}^{2} \right)$) for $r_{1}=r_{1}'$. The proof is identical
for $r_{1}=r_{1}''$.

If $k_{2}$ is even, then we have $2k_{1}<k_{2}$.
If $k_{2}$ is odd, then we have $2k_{1}-1<k_{2}$. So in both cases,
$2k_{1} \leq k_{2}$ (i.e., $v_{p} \left( 4b^{2} \right) \leq v_{p} \left( r_{1}' \right)$)
and $v_{p}(y) \geq k_{1}$.

We now show that in this case we have
$v_{p} \left( 4b^{2} \right)<v_{p}\left( r_{1}'' \right)$.

Since $v_{p}(y) \geq v_{p}(b)$, from the definition of $r_{1}'$ and
$v_{p}\left( r_{1}' \right)>v_{p} \left( 4b^{2} \right)$, we also have
$v_{p}(au)>v_{p} \left( b^{2} \right)$. So the $p$-adic valuation of
each term in the definition of $r_{1}'$ is greater than $v_{p} \left( b^{2} \right)$.
Thus $v_{p} \left( r_{1}'' \right)>v_{p} \left( b^{2} \right)$ too.

Combining this with $r_{1}'r_{1}''=u^{2}N_{\alpha}$, we have
$v_{p} \left( b^{2}r_{1}' \right)<v_{p} \left( r_{1}'r_{1}'' \right)=v_{p} \left( u^{2}N_{\alpha} \right)$,
as desired.

\vspace*{1.0mm}

The proof when $p=2$ is similar. But we will use the following information later
in the proof of Proposition~\ref{prop:quad-rep}. We have
\begin{equation}
\label{eq:val-32-p2}
v_{2} \left( g_{1}^{2} \right)
= \left\{
\begin{array}{ll}
2k_{1}+k_{2}+2,
& \text{if $k_{2}$ is even and $2k_{1}+2 \leq k_{2}$,} \\
2k_{2},
& \text{if $k_{2}$ is even and $2k_{1}+2 \geq k_{2}$,} \\
2k_{1}+k_{2}+1,
& \text{if $k_{2}$ is odd and $2k_{1}+1 \leq k_{2}$,} \\
2k_{2},
& \text{if $k_{2}$ is odd and $2k_{1}+1 \geq k_{2}$.}
\end{array}
\right.
\end{equation}

As above, we have
$v_{2}(au) \geq \min \left( k_{2}, v_{2} \left( tb^{2} \pm 2by \right) \right)$.

Applying this and the right-hand expression in \eqref{eq:rel-y} for $4y^{2}$,
we obtain
\begin{align*}
2v_{2}(y)+2 & \geq \min \left( v_{2}(2atu), v_{2} \left( b^{2} \left( 4N_{\varepsilon}+2du^{2} \right) \right) \right) \\
& \geq \min \left( v_{2} \left( 2t^{2}b^{2} \pm 4bty \right), k_{2}+v_{p}(t)+1, v_{2} \left( 2b^{2} \right) \right) \\
& \geq \min \left( v_{2} \left( 4bty \right), k_{2}+v_{2}(t)+1, 2k_{1}+1 \right) \\
& \geq \min \left( k_{1}+v_{2}(t)+v_{2}(y)+2, k_{2}+v_{2}(t)+1, 2k_{1}+1 \right).
\end{align*}

Since $v_{2}(y)$ is an integer, this implies that
\begin{equation}
\label{eq:v2y-LB}
v_{2}(y) \geq \min \left( k_{1}, \left( k_{2}-1 \right)/2 \right)
\end{equation}
and
$v_{2}(y) \geq \min \left( k_{1}, k_{2}/2 \right)$, if $k_{2}$ is even.
\end{proof}

\begin{lemma}
\label{lem:quad-rep-2}
In the notation of Proposition~\ref{prop:quad-rep}, put
$r_{1}'=tb^{2}+au+2by$, $g_{1}'^{2}=\gcd \left( 4b^{2}r_{1}'/\core \left( r_{1}' \right), r_{1}'^{2} \right)$,
$r_{1}''=tb^{2}+au-2by$ and $g_{1}''^{2}=\gcd \left( 4b^{2}r_{1}''/\core \left( r_{1}'' \right), r_{1}''^{2} \right)$.

\noindent
{\rm (a)} For all primes, $p$,
\[
\max \left( v_{p} \left( r_{1}'/g_{1}'^{2} \right),
v_{p} \left( r_{1}''/g_{1}''^{2} \right) \right)
\leq v_{p} \left( \rad \left( \gcd \left( uN_{\alpha}, N_{\varepsilon}\core \left( N_{\alpha} \right) \right) \right) \right).
\]

\noindent
{\rm (b)} For all primes, $p$, if $v_{p} \left( \core \left( N_{\alpha} \right) \right)>0$,
then
\begin{equation}
\label{eq:rep-c}
\min \left( v_{p} \left( r_{1}'/g_{1}'^{2} \right),
v_{p} \left( r_{1}''/g_{1}''^{2} \right) \right)
\leq 0.
\end{equation}

\noindent
{\rm (c)} If $r_{1}'$ is even, then $v_{2} \left( r_{1}'/g_{1}'^{2} \right)<0$.
The same result holds for $r_{1}''$.
\end{lemma}

\begin{proof}
(a) We will prove that for all primes, $p$, we have
$v_{p} \left( r_{1}'/g_{1}'^{2} \right) \leq
v_{p} \left( \rad \left( \gcd \left( uN_{\alpha}, N_{\varepsilon}\core \left( N_{\alpha} \right) \right) \right) \right)$.
The proof is the same for $r_{1}''/g_{1}''^{2}$.

From \eqref{eq:g1Sqr-simp}, we can write
\begin{equation}
\label{eq:quad-rep-2a}
\frac{r_{1}'}{g_{1}'^{2}}= \frac{\core \left( r_{1}' \right)}{\gcd \left( 4b^{2}, r_{1}'\core \left( r_{1}' \right) \right)}.
\end{equation}

Since $v_{p} \left( \core \left( r_{1}' \right) \right) \leq 1$, for all primes,
$p$, it follows from this expression that $v_{p} \left( r_{1}'/g_{1}'^{2} \right) \leq 1$,
for all primes, $p$. So it remains to show that if $v_{p} \left( r_{1}'/g_{1}'^{2} \right)=1$,
then we also have $v_{p} \left( \rad \left( \gcd \left( uN_{\alpha}, N_{\varepsilon}\core \left( N_{\alpha} \right) \right) \right) \right)=1$.

From the expression above for $r_{1}'/g_{1}'^{2}$, we see that we are considering
primes, $p$, such that $p| \core \left( r_{1}' \right)$ and $p \nmid \left( 4b^{2} \right)$
(i.e., $p \nmid (2b)$).

From the observation that $r_{1}'r_{1}''
=N_{\alpha}u^{2}$, we have $r_{1}'| \left( N_{\alpha}u^{2} \right)$.

We start by showing that if $v_{p} \left( r_{1}'/g_{1}^{2} \right)=1$, then
$p| \left( N_{\varepsilon}\core \left( N_{\alpha} \right) \right)$.

If $p| \core \left( r_{1}' \right)$, but $p \nmid \gcd \left( r_{1}', r_{1}'' \right)$,
then $v_{p} \left( u^{2}N_{\alpha} \right)$ must be odd. Hence
$p | \core \left( N_{\alpha} \right)$.

So we need only consider what happens when
$p|\core \left( r_{1}' \right)$, $p \nmid (2b)$ and
$p | \gcd \left( r_{1}', r_{1}'' \right)$.

Suppose that $v_{p} \left( r_{1}' \right)=k_{1}$ and
$v_{p} \left( r_{1}'' \right)=k_{2}$ for positive integers $k_{1}$ and $k_{2}$.
Then $v_{p} \left( u^{2}N_{\alpha} \right)=k_{1}+k_{2}$. If $k_{1}+k_{2}$ is
odd, then $p | \core \left( N_{\alpha} \right)$.

We also know that $k_{1}$ is odd (otherwise, $p \nmid \core \left( r_{1}' \right)$).
So we need only consider when $k_{2}$ is also odd.

Since $p|\left( r_{1}'-r_{1}'' \right)$, $\left( r_{1}'-r_{1}'' \right)=4by$ and $p \nmid (2b)$, it follows that $p|y$.
We also have $p$ is a divisor of $r_{1}'+r_{1}''=2au+2b^{2}t$.

We can write $4y^{2}=2atu+\left( t^{2}+du^{2} \right)b^{2}$ and find that
$p$ is a divisor of $t \left( r_{1}'+r_{1}'' \right)-4y^{2}
=b^{2} \left( t^{2}-du^{2} \right)$.
Since $p \nmid b$, we have $p|N_{\varepsilon}$,
completing the proof that $p| \left( N_{\varepsilon}\core \left( N_{\alpha} \right) \right)$.

Now we show that $p| \left( uN_{\alpha} \right)$ holds if
$v_{p} \left( r_{1}'/g_{1}^{2} \right)=1$. But this is immediate from
$r_{1}'r_{1}''=u^{2}N_{\alpha}$.

This completes the proof of part~(a).

\vspace*{1.0mm}

(b) Suppose $v_{p} \left( \core \left( N_{\alpha} \right) \right)>0$. That is,
$v_{p} \left( \core \left( N_{\alpha} \right) \right)=1$. We will show that at
least one of $v_{p} \left( r_{1}'/g_{1}'^{2} \right)$ or $v_{p} \left( r_{1}''/g_{1}''^{2} \right)$
must be non-positive. Otherwise, by \eqref{eq:quad-rep-2a},
$v_{p} \left( \core \left( r_{1}' \right) \right)=
v_{p} \left( \core \left( r_{1}'' \right) \right)=1$. So $v_{p} \left( r_{1}'r_{1}'' \right)$
is even.
But from $r_{1}'r_{1}''=u^{2}N_{\alpha}$ and
$v_{p} \left( \core \left( N_{\alpha} \right) \right)=1$, we see that
$v_{p} \left( r_{1}'r_{1}'' \right)$ is odd.
Hence at least one of $v_{p} \left( r_{1}'/g_{1}'^{2} \right)$ or $v_{p} \left( r_{1}''/g_{1}''^{2} \right)$
is non-positive, as required.

\vspace*{1.0mm}

(c) This follows from \eqref{eq:quad-rep-2a}. We have
$v_{2} \left( \core \left( r_{1}' \right) \right)=0$ or $1$.
In both cases, since $r_{1}'$ is even, we find that $v_{2}$ of the denominator
of the right-hand side of \eqref{eq:quad-rep-2a} is strictly greater than
$v_{2}$ of its numerator.
\end{proof}

We shall need the following lemma for the proofs of parts~(b) and (c)
of Proposition~\ref{prop:quad-rep}.

\begin{lemma}
\label{lem:quad-rep-3}
With the same notation as in Lemma~\ref{lem:quad-rep-2},
$\gcd \left( r_{1}', r_{1}'' \right)$ is even, unless
$N_{\alpha} \equiv 1 \pmod{4}$ and $4|d$.
\end{lemma}

\begin{proof}
From the expressions for $r_{1}'$ and $r_{1}''$, we need only show that $tb^{2}+au$
is even unless $N_{\alpha} \equiv 1 \pmod{4}$ and $4 \mid d$. We do so by
showing that if either
(i) $tb^{2}$ even and $au$ odd or (ii) $tb^{2}$ odd and $au$ even, then
$N_{\alpha} \equiv 1 \pmod{4}$ and $4|d$.

(i) We first suppose that $t$ is even in case~(i). Since $u$ is odd and $\left( t^{2}-du^{2} \right)/4
\in \bbZ$, we have $4|d$. But then $N_{\alpha}=a^{2}-db^{4} \equiv 1 \pmod{4}$,
since $a$ must be odd.

Now suppose that $t$ is odd. Since $\left( t^{2}-du^{2} \right)/4 \in \bbZ$
and $u$ is odd, it follows that $d \equiv 1 \pmod{4}$. Here $b$ must be even,
since $tb^{2}$ is
even. Again $N_{\alpha}=a^{2}-db^{4} \equiv 1 \pmod{4}$,
since $a$ must be odd. Applying these to the middle expression in \eqref{eq:rel-y},
we find that $4y^{2} \equiv 2 \pmod{4}$, so this case is not possible.

(ii) If $u$ is even, then since $t$ is odd, $N_{\varepsilon}$ cannot be an integer.
So we must have $a$ even and $b$, $t$ and $u$ all odd. Since $t$ and $u$ are both
odd, $d \equiv 1 \pmod{4}$. So we have $N_{\alpha} \equiv 3 \pmod{4}$. Hence
$\left( b^{2}t+au \right)^{2}-N_{\alpha}u^{2} \equiv 2 \pmod{4}$.
But this quantity is $b^{2}$ times the middle expression in \eqref{eq:rel-y}, so
it is $4b^{2}y^{2}$. Thus (ii) is not possible.
\end{proof}

\begin{proof}[Proof of Proposition~\ref{prop:quad-rep}]
(a) We start by writing $\left( x + N_{\varepsilon}\sqrt{N_{\alpha}} \right)/
\left( a+\sqrt{N_{\alpha}} \right)$ as a square.

By rationalising the denominator of $\left( x + N_{\varepsilon}\sqrt{N_{\alpha}} \right)/
\left( a+\sqrt{N_{\alpha}} \right)$ and then applying the expressions for
$N_{\varepsilon}$ and $N_{\alpha}$, we find that
\[
\frac{x + N_{\varepsilon}\sqrt{N_{\alpha}}}{a+\sqrt{N_{\alpha}}}
=\frac{4ax-a^{2}t^{2}+b^{4}dt^{2}+a^{2}du^{2}-b^{4}d^{2}u^{2}
+\left( at^{2}-adu^{2}-4x\right)\sqrt{N_{\alpha}}}{4b^{4}d}.
\]

Substituting
\begin{equation}
\label{eq:rel-x}
4x=a \left( t^{2}+ du^{2} \right) + 2b^{2}dtu,
\end{equation}
from the expression for $x+y^{2}\sqrt{d}$ in \eqref{eq:quad-rep-assumption},
we obtain
\[
\frac{x + N_{\varepsilon}\sqrt{N_{\alpha}}}{a+\sqrt{N_{\alpha}}}
=\frac{2a^{2}u^{2}+2ab^{2}tu+b^{4}t^{2}-db^{4}u^{2}-2u\left( au+b^{2}t \right)\sqrt{N_{\alpha}}}{4b^{4}}.
\]

We can write $2a^{2}u^{2}+2ab^{2}tu+b^{4}t^{2}-db^{4}u^{2}=\left( au+b^{2}t \right)^{2}+u^{2}N_{\alpha}$,
so
\[
\frac{x + N_{\varepsilon}\sqrt{N_{\alpha}}}{a+\sqrt{N_{\alpha}}}
=\left( \frac{\left( au+b^{2}t \right)-u\sqrt{N_{\alpha}}}{2b^{2}} \right)^{2}.
\]

With $r_{1}=tb^{2}+au\pm 2by$ and $s_{1}=-u$, a routine calculation, along with
\eqref{eq:rel-y}, shows that
\[
\left( r_{1}+s_{1}\sqrt{N_{\alpha}} \right)^{2}
=2r_{1} \left( \left( au+b^{2}t \right)-u\sqrt{N_{\alpha}} \right).
\]

Hence
\begin{align}
\label{eq:rep-lemma-expr1}
\left( 4b^{2}r_{1} \right)^{2} \left( x+N_{\varepsilon}\sqrt{N_{\alpha}} \right)
&= \left( a+\sqrt{N_{\alpha}} \right) \left( r_{1}+s_{1}\sqrt{N_{\alpha}} \right)^{4} \nonumber \\
&= \left( a+\sqrt{N_{\alpha}} \right) \left( r_{1}+s_{1}'\sqrt{\core \left( N_{\alpha} \right)} \right)^{4},
\end{align}
where $s_{1}'=s_{1}\sqrt{N_{\alpha}/\core \left( N_{\alpha} \right)}$.

From this relationship, we now show that we can find such a relationship for
$x+N_{\varepsilon}\sqrt{N_{\alpha}}$ with the conditions on $f$ in part~(a)
satisfied.

We want to take the largest common factor, $g_{1}$, of $r_{1}$ and $s_{1}'$ such
that $4b^{2}r_{1}/g_{1}^{2}$ is also an integer. That is,
\[
g_{1}^{2}= \gcd \left( \gcd \left( 4b^{2} r_{1}/\core \left( r_{1} \right), r_{1}^{2} \right), s_{1}'^{2} \right).
\]

First, note that we can write
\begin{equation}
\label{eq:g1Sqr-simp}
\gcd \left( 4b^{2} r_{1}/\core \left( r_{1} \right), r_{1}^{2} \right)
= r_{1} \gcd \left( 4b^{2}, r_{1}\core \left( r_{1} \right) \right) /\core \left( r_{1} \right).
\end{equation}

From Lemma~\ref{lem:quad-rep-1}, we see that
\[
g_{1}^{2}=\gcd \left( 4b^{2} r_{1}/\core \left( r_{1} \right), r_{1}^{2} \right).
\]

From Lemma~\ref{lem:quad-rep-2}(a), $4b^{2}r_{1}/g_{1}^{2}$ is a divisor
of
\[
4b^{2} \rad \left( \gcd \left( uN_{\alpha}, N_{\varepsilon}\core \left( N_{\alpha} \right) \right) \right)
=4b^{2} \rad \left( \core \left( N_{\alpha} \right) \gcd \left( uN_{\alpha}/\core \left( N_{\alpha} \right), N_{\varepsilon} \right) \right).
\]

We put $r=r_{1}/g_{1}$, for either choice of $r_{1}$, $s=s_{1}'/g_{1}$ and
$f=4b^{2}r_{1}/g_{1}^{2}$.
By the definition of $g_{1}^{2}$, we see that $f,r,s \in \bbZ$.

Lastly, we choose which of the two possible values of $r_{1}$ to use. As in the
statement of Lemma~\ref{lem:quad-rep-2}, we let $r_{1}'=tb^{2}+au+2by$ and
$r_{1}''=tb^{2}+au-2by$ be the two possibilities. We let $g_{1}'$ and $g_{1}''$
be the associated values of $g_{1}$. We will use
Lemma~\ref{lem:quad-rep-2}(b). We divide the prime divisors of $\core \left( N_{\alpha} \right)$
into three disjoint sets:
\begin{align*}
I'   & = \left\{ p : p \text{ prime, with } p|\core \left( N_{\alpha} \right)
\text{ and }
v_{p} \left( r_{1}'/g_{1}'^{2} \right)<v_{p} \left( r_{1}''/g_{1}''^{2} \right) \right\}, \\
I''  & = \left\{ p : p \text{ prime, with } p|\core \left( N_{\alpha} \right)
\text{ and }
v_{p} \left( r_{1}''/g_{1}''^{2} \right)<v_{p} \left( r_{1}'/g_{1}'^{2} \right) \right\}, \\
I''' & = \left\{ p : p \text{ prime, with } p|\core \left( N_{\alpha} \right)
\text{ and }
v_{p} \left( r_{1}'/g_{1}'^{2} \right) = v_{p} \left( r_{1}''/g_{1}''^{2} \right) \right\}.
\end{align*}

With $d'=\prod_{p \in I'} p$, $d''=\prod_{p \in I''} p$ and
$d'''=\prod_{p \in I'''} p$, we have
$d' d'' d'''= \core \left( N_{\alpha} \right)$.

We put $r_{1}=r_{1}'$ if $d'>d''$ and $r_{1}=r_{1}''$ otherwise.
By Lemma~\ref{lem:quad-rep-2}(b),
\begin{align*}
v_{p} \left( r_{1}/g_{1}^{2} \right) \leq 0
& \leq v_{p} \left( \core \left( N_{\alpha} \right)/p \right) \\
& \leq v_{p} \left( \rad \left( \left( \core \left( N_{\alpha} \right)/p \right) \gcd \left( uN_{\alpha}/\core \left( N_{\alpha} \right), N_{\varepsilon} \right) \right) \right),
\end{align*}
for all $p \in I' \cup I'''$ in the first case and
for all $p \in I'' \cup I'''$ in the second case.
Hence with
$f'=\core \left( N_{\alpha} \right)/ \max \left( d'd''', d''d''' \right)$,
we see that $f'<\sqrt{\core \left( N_{\alpha} \right)}$ and that
$4b^{2}r_{1}/g_{1}^{2}$ is a divisor of
\[
4b^{2} \rad \left( f' \gcd \left( uN_{\alpha}/\core \left( N_{\alpha} \right), N_{\varepsilon} \right) \right).
\]

\vspace*{3.0mm}

(b)
If $4b^{2}r_{1}/g_{1}^{2}$ is odd, then there is nothing to prove, as
$4b^{2}r_{1}/g_{1}^{2}$ dividing $b^{2} \rad \left( \gcd \left( uN_{\alpha}, N_{\varepsilon} \right) \right)$
follows from
$4b^{2}r_{1}/g_{1}^{2}$ being a divisor of $4b^{2} \rad \left( \gcd \left( uN_{\alpha}, N_{\varepsilon} \right) \right)$.
So we put $f=4b^{2}r_{1}/g_{1}^{2}$ and need only consider when
$4b^{2}r_{1}/g_{1}^{2}$ is even.

We consider the four cases for the parity of $r$ and $s$.

(b-i) Suppose that $r$ and $s$ are both odd.

Here we will replace $r$ and $s$.

We can write
$(r+si)/(1+i)=(r+s)/2+(s-r)i/2$.
Since $r$ and $s$ are both odd and $(1+i)^{4}=-4$, we find that $(s \pm r)/2$ are
both integers and
\[
-\left( \frac{2b^{2}r}{g_{1}^{2}} \right)^{2} \left( x+N_{\varepsilon} \sqrt{N_{\alpha}} \right)
= \left( a + \sqrt{N_{\alpha}} \right) \left( (r+s)/2 + (s-r)i/2 \right)^{4}.
\]

Since $4b^{2}r/g_{1}^{2}$ is even, we know that $2b^{2}r/g_{1}^{2} \in \bbZ$.

Since $-N_{\alpha}$ is a square, we have $N_{\alpha} \not\equiv 1 \pmod{4}$, so
by Lemmas~\ref{lem:quad-rep-3} and \ref{lem:quad-rep-2}(c), we see that
$\left( 4b^{2}r/g_{1}^{2} \right) | \left( 2b^{2} \rad \left( \gcd \left( uN_{\alpha}, N_{\varepsilon} \right) \right) \right)$.
Hence $2b^{2}r/g_{1}^{2}$ is a divisor of
$b^{2} \rad \left( \gcd \left( uN_{\alpha}, N_{\varepsilon} \right) \right)$
and so we take $f=2b^{2}r/g_{1}^{2}$ in this case.

\vspace*{3.0mm}

(b-ii) Suppose that $r=r_{1}/g_{1}$ is odd and $s=s_{1}'/g_{1}$ is even.

In this case, as well as cases~(b-iii) and (b-iv), we will assume that $r_{1}=r_{1}'$.
The proof is identical using $r_{1}=r_{1}''$ instead.

Since $r$ is odd, we have $v_{2} \left( g_{1} \right)=v_{2} \left( r_{1} \right)$.
So from the definition of $g_{1}^{2}$ in the proof of part~(a), it follows that
$v_{2} \left( 4b^{2}/ \core \left( r_{1} \right) \right) \geq v_{2} \left( r_{1} \right)$.

For this case and case~(b-iii), if $v_{2} \left( r_{1}' \right) \neq v_{2} \left( r_{1}'' \right)$,
then we must have $v_{2} \left( au+b^{2}t \right)=v_{2}(2by)$, since $r_{1}''=r_{1}'+4by$.

(b-ii-1) Suppose first that $v_{2} \left( r_{1} \right)$ is even.

In this case, $v_{2} \left( 4b^{2} \right) \geq v_{2} \left( r_{1} \right)=v_{2} \left( g_{1} \right)$.
If we have equality,
then $f=4b^{2}r_{1}/g_{1}^{2}$ is odd and so $f| \left( b^{2} \rad \left( \gcd \left( uN_{\alpha}, N_{\varepsilon} \right) \right) \right)$,
as required. So we consider $v_{2} \left( 4b^{2} \right) > v_{2} \left( r_{1} \right)$.
Since both of these quantities are even, we must have $v_{2} \left( 4b^{2} \right) \geq v_{2} \left( r_{1} \right)+2$.
That is, $v_{2} \left( b^{2} \right) \geq v_{2} \left( r_{1} \right)$.

Since $s$ is even, we have
\[
v_{2} \left( s_{1}'^{2} \right)=v_{2} \left( u^{2}N_{\alpha} \right)
=v_{2} \left( r_{1}' \right) + v_{2} \left( r_{1}'' \right)>2v_{2} \left( g_{1} \right).
\]

Since $r$ is odd, we also have $2v_{2} \left( g_{1} \right)
=2v_{2} \left( r_{1}' \right)$. Hence $v_{2} \left( r_{1}'' \right)>v_{2} \left( r_{1}' \right)$.
As we stated near the start of this case (case~(b-ii)), since $v_{2} \left( r_{1}' \right) \neq v_{2} \left( r_{1}'' \right)$,
we have $v_{2} \left( au+b^{2}t \right)=v_{2}(2by)$.

We saw in the proof of Lemma~\ref{lem:quad-rep-1} that $v_{2}(y) \geq \min \left( v_{2}(b), v_{2} \left( r_{1}' \right)/2 \right)$.
Since $v_{2} \left( b^{2} \right) \geq v_{2} \left( r_{1} \right)$, it follows
that $v_{2}(y) \geq v_{2} \left( r_{1}' \right)/2$. Hence
\[
v_{2}(2by) \geq 1+v_{2} \left( r_{1} \right)/2+v_{2} \left( r_{1} \right)/2=v_{2} \left( r_{1} \right)+1.
\]

Thus $v_{2} \left( au+b^{2}t \right)=v_{2}(2by) \geq v_{2} \left( r_{1} \right)+1$,
which implies that $v_{2} \left( r_{1} \right)>v_{2} \left( r_{1} \right)$. Therefore,
the case when $v_{2} \left( r_{1} \right)$ is even and $v_{2} \left( b^{2} \right) \geq v_{2} \left( r_{1} \right)$
never occurs and we find that in case~(b-ii-1), 
$f| \left( b^{2} \rad \left( \gcd \left( uN_{\alpha}, N_{\varepsilon} \right) \right) \right)$,
as required.

(b-ii-2) Now suppose that $v_{2} \left( r_{1} \right)$ is odd.

Here $v_{2} \left( 2b^{2} \right) \geq v_{2} \left( r_{1} \right)=v_{2} \left( g_{1} \right)$.

As in case~(b-ii-1), since $s$ is even, we have $v_{2} \left( r_{1}' \right) \neq
v_{2} \left( r_{1}'' \right)$. Hence $v_{2} \left( au+b^{2}t \right)=v_{2} \left( 2by \right)$.

Recall from \eqref{eq:v2y-LB} that
$v_{2}(y) \geq \min \left( v_{2}(b), \left( v_{2} \left( r_{1} \right)-1 \right) /2 \right)$.
Hence, $v_{2} \left( au+b^{2}t \right)=v_{2} \left( 2by \right) \geq v_{2} \left( r_{1} \right)$.
Since $v_{2} \left( au+b^{2}t \right)=v_{2} \left( 2by \right)$, we have
$v_{2} \left( r_{1} \right) \geq v_{2} \left( 2by \right)+1>v_{2} \left( r_{1} \right)$,
which is not possible.

This completes the proof for case~(b-ii).

\vspace*{3.0mm}

(b-iii) Suppose that $r=r_{1}/g_{1}$ is even and $s=s_{1}'/g_{1}$ is odd.

Since $r$ is even, we must have $v_{2} \left( g_{1}^{2} \right)
= v_{2} \left( 4b^{2}r_{1}/\core \left( r_{1} \right) \right)< v_{2} \left( r_{1}^{2} \right)$.

Since $v_{2} \left( s_{1}'^{2} \right)=v_{2} \left( u^{2} N_{\alpha} \right)
=v_{2} \left( r_{1}' \right) + v_{2} \left( r_{1}'' \right)$, by the assumption
that $s_{1}'/g_{1}$ is odd, we have
$v_{2} \left( r_{1}' \right) + v_{2} \left( r_{1}'' \right)
=v_{2} \left( g_{1}^{2} \right)$. Since $r$ is even, we also have
$v_{2} \left( g_{1}^{2} \right) < 2v_{2} \left( r_{1} \right)$, it follows that
$v_{2} \left( r_{1}'' \right) \neq v_{2} \left( r_{1} \right)$. Hence, by the
comment at
the start of case~(b-ii), we have $v_{2} \left( au+b^{2}t \right) = v_{2} (2by)$.

(b-iii-1) If $v_{2} \left( r_{1} \right)$ is even, then
$v_{2} \left( g_{1}^{2} \right)= v_{2} \left( 4b^{2}r_{1} \right)$.
Hence $f=4b^{2}r_{1}/g_{1}^{2}$ is odd and
$f| \left( b^{2} \rad \left( \gcd \left( uN_{\alpha}, N_{\varepsilon} \right) \right) \right)$
as required.

\vspace*{1.0mm}

(b-iii-2) If $v_{2} \left( r_{1} \right)$ is odd, then
$v_{2} \left( 2b^{2}r_{1} \right)<v_{2} \left( r_{1}^{2} \right)$. So
$v_{2}(b)< \left( v_{2} \left( r_{1} \right)-1 \right)/2$.
Recall from \eqref{eq:v2y-LB} that
$v_{2}(y) \geq \min \left( v_{2}(b), \left( v_{2} \left( r_{1} \right)-1 \right) /2 \right)$.
So $v_{2}(y) \geq v_{2}(b)$.

Hence, $v_{2} \left( au+b^{2}t \right)=v_{2} \left( 2by \right) \geq 2v_{2}(b)+1$.
Since $v_{2} \left( au+b^{2}t \right)=v_{2} \left( 2by \right)$, we have
$v_{2} \left( r_{1} \right) \geq v_{2} \left( 2by \right)+1 \geq 2v_{2}(b)+2$.
Using the same argument, we have
$v_{2} \left( r_{1}'' \right) \geq 2v_{2}(b)+2$ too.

But $v_{2} \left( 2b^{2}r_{1} \right)=v_{2} \left( g_{1}^{2} \right)
=v_{2} \left( s_{1}'^{2} \right)=v_{2} \left( r_{1}' \right)+v_{2} \left( r_{1}'' \right)$.
Hence $v_{2} \left( r_{1}'' \right) = v_{2} \left( 2b^{2} \right)$.
But this contradicts $v_{2} \left( r_{1}'' \right) \geq 2v_{2}(b)+2$.
Hence $v_{2} \left( r_{1} \right)$ cannot be odd.

\vspace*{3.0mm}

(b-iv) Suppose that $r$ and $s$ are both even.

Since $r$ is even, we know that $v_{2} \left( g_{1}^{2} \right)
=v_{2} \left( 4b^{2}r_{1}/\core \left( r_{1} \right) \right)
<v_{2} \left( r_{1}^{2} \right)$.
So $v_{2}(f)=0$ or $1$, depending on whether
$v_{2} \left( r_{1} \right)$ is even or odd. So, as stated at the start of the proof of
part~(b), we need only consider the latter case.

Since $s$ is even, we have
\[
2 \leq v_{2} \left( s^{2} \right)=v_{2} \left( u^{2}N_{\alpha} \right)-v_{2} \left( 2b^{2}r_{1}' \right)
=v_{2} \left( r_{1}'' \right)-v_{2} \left( 2b^{2} \right).
\]
So $v_{2} \left( r_{1}'' \right) \geq 2v_{2}(b)+3$.

Similarly, $2 \leq v_{2} \left( r_{1}'^{2} \right)-v_{2} \left( 2b^{2}r_{1}' \right)
=v_{2} \left( r_{1}' \right) - v_{2} \left( 2b^{2} \right)$, so
$v_{2} \left( r_{1}' \right) \geq 2v_{2}(b)+3$.

A consequence of these two inequalities is that $64| \left( u^{2}N_{\alpha} \right)$.

From \eqref{eq:rel-y}, we have $4y^{2}=2atu+2b^{2}du^{2}+4b^{2}N_{\varepsilon}$.
So $atu+b^{2}du^{2}$ must be even. I.e., $atu$ and $b^{2}du^{2}$ have the same parity.

(b-iv-1) If $b^{2}du^{2}$ is odd, then $a$ and $t$ are also odd. Since $4N_{\varepsilon}
=t^{2}-du^{2}$ is
divisible by $4$ and $t$ and $u$ are both odd, it must be the case that $d \equiv 1 \pmod{4}$.

If $d \equiv 1 \pmod{8}$, then $4N_{\varepsilon}
=t^{2}-du^{2} \equiv 0 \pmod{8}$ and so $N_{\varepsilon}$ is even. Furthermore, since $abd$
is odd, we have $N_{\alpha}$ is even, so
$b^{2} \rad \left( \gcd \left( uN_{\alpha}, N_{\varepsilon} \right) \right)$ is
also even and the desired conclusion holds.

If $d \equiv 5 \pmod{8}$, then $N_{\alpha}=a^{2}-db^{4} \equiv 4 \pmod{8}$. Since
$u$ is odd, we get a contradiction with $64| \left( u^{2}N_{\alpha} \right)$.

(b-iv-2) We now consider $b^{2}du^{2}$ even. If $b$ is even, then
\[
v_{2} \left( b^{2} \rad \left( \gcd \left( uN_{\alpha}, N_{\varepsilon} \right) \right) \right) \geq 2
> v_{2}(f)=1.
\]

So the desired conclusion holds. Therefore, we may assume that $b$ is odd and $du^{2}$
is even in what follows. We showed above that $uN_{\alpha}$ is even. We will show
here that $N_{\varepsilon}$ is also even. These facts will suffice to show that
$v_{2} \left( b^{2} \rad \left( \gcd \left( uN_{\alpha}, N_{\varepsilon} \right) \right) \right) \geq 1
=v_{2}(f)$ and hence to prove part~(b) in this case.

To prove that $N_{\varepsilon}$ is even, we will assume that it is odd and obtain
a contradiction.

Since $4N_{\varepsilon}=t^{2}-du^{2} \in \bbZ$ and $du^{2}$ is even, we have $t$
is even and hence $4| \left( du^{2} \right)$.

If $u$ is odd, then $4|d$. Also $N_{\alpha}$ is even, so $a$ is even. As a result,
we find that $v_{2}(2atu), v_{2} \left( 2b^{2}du^{2} \right) \geq 3$. So under
our assumption that $N_{\varepsilon}$ is odd, from \eqref{eq:rel-y} we have
$v_{2} \left( 4y^{2} \right)=2$. I.e., $y$ is odd.

The same reasoning shows that $y$ is odd if $u$ is even under the assumption that
$N_{\varepsilon}$ is odd.

Suppose that $v_{2} \left( r_{1}' \right)=v_{2} \left( r_{1}'' \right)$. Then $v_{2} \left( 4by \right)=v_{2} \left( \pm \left( r_{1}'-r_{1}'' \right) \right)>v_{2} \left( r_{1}' \right)$. So $v_{2}(y)>v_{2} \left( r_{1}' \right)-v_{2}(b)-2>v_{2}(b)+1$.

Suppose that $v_{2} \left( r_{1}' \right)<v_{2} \left( r_{1}'' \right)$. Then $v_{2} \left( 4by \right)=v_{2} \left( r_{1}' \right)$. So $v_{2}(y) \geq v_{2}(b)+1$.

In both cases, we find that $y$ is even. This contradicts what we obtained
($y$ is odd) under the assumption that $N_{\varepsilon}$ is odd. Hence we know
that in the case when $b^{2}du^{2}$ is even,
$N_{\varepsilon}$ is even and the assertion in part~(b) of the lemma holds.

\vspace*{3.0mm}

(c) We consider two cases.

(c-i) Suppose that $\left| \core \left( N_{\alpha} \right) \right|=p$,
for a prime, $p$. Then we have $f'=1$ in part~(a) and so
$f| \left( 4b^{2} \rad \left( \gcd \left( uN_{\alpha}/\core \left( N_{\alpha} \right), N_{\varepsilon} \right)
\right) \right)$. Hence the result holds if $N_{\alpha} \equiv 1 \pmod{4}$ and
$4|d$.

From Lemma~\ref{lem:quad-rep-3}, $r$ is even unless
$N_{\alpha} \equiv 1 \pmod{4}$ and $4|d$. So, from Lemma~\ref{lem:quad-rep-2}(c),
we have $v_{2} \left( r/g_{1}^{2} \right) \leq -1$.
With $f=4b^{2}r/g_{1}^{2}$, we have
$f| \left( 2b^{2} \rad \left( \gcd \left( uN_{\alpha}/\core \left( N_{\alpha} \right), N_{\varepsilon} \right)
\right) \right)$.

(c-ii) Suppose that $\left| \core \left( N_{\alpha} \right) \right|=2p$, where
$p$ is an odd prime. From Lemma~\ref{lem:quad-rep-3}, we see that $r_{1}'$ and
$r_{1}''$ are both even. Hence, from Lemma~\ref{lem:quad-rep-2}(c), we have
\[
\max \left( v_{2} \left( r_{1}'/g_{1}'^{2} \right),
v_{2} \left( r_{1}''/g_{1}''^{2} \right) \right)<0
\leq v_{2} \left( \rad \left( \gcd \left( uN_{\alpha}/\core \left( N_{\alpha} \right), N_{\varepsilon} \right) \right) \right).
\]

Similar to the end of the proof of part~(a), we let $r_{1}=r_{1}'$ if
$v_{p} \left( r_{1}'/g_{1}'^{2} \right) \leq v_{p} \left( r_{1}''/g_{1}^{2} \right)$
and $r_{1}=r_{1}''$ otherwise. From Lemma~\ref{lem:quad-rep-2}(b),
we find that
\[
v_{p} \left( r_{1}/g_{1}^{2} \right)
\leq 0 \leq v_{p} \left( \rad \left( \gcd \left( uN_{\alpha}/\core \left( N_{\alpha} \right), N_{\varepsilon} \right) \right) \right).
\]
So with $f=4b^{2}r_{1}/g_{1}^{2}$, we have
$f| \left( 2b^{2} \rad \left( \gcd \left( uN_{\alpha}/\core \left( N_{\alpha} \right), N_{\varepsilon} \right)
\right) \right)$.
\end{proof}

\subsection{Lower bounds for $y_{k}$'s}

Next we bound the $y_{k}$'s from below.

\begin{lemma}
\label{lem:Y-LB}
Let the $y_{k}$'s be defined by $\eqref{eq:yk-defn}$ with the notation and assumptions
there. Suppose that $N_{\alpha}<0$.
Let $K$ be the largest negative integer such that $y_{K}>b^{2}$.
 
\noindent
{\rm (a)} Put $\overline{\alpha}=a-b^{2}\sqrt{d}$. We have
\begin{equation}
\label{eq:yLB1-gen}
y_{k}>
\left\{
\begin{array}{ll}
\dfrac{\alpha\varepsilon^{2k}}{2\sqrt{d}} & \text{for $k \geq 0$,} \\
\dfrac{-\overline{\alpha} \, \varepsilon^{2|k|}}{2\sqrt{d}} & \text{for $k<0$.}
\end{array}
\right.
\end{equation}

\noindent
{\rm (b)}
For all $k$, $2y_{k}$ is a positive integer.
The sequences $\left( y_{k} \right)_{k \geq 0}$ and 
$\left( y_{K+1}, y_{K}, y_{K-1}, y_{K-2}, \ldots \right)$ are increasing sequences of positive numbers.

\noindent
{\rm (c)} We have
\begin{equation}
\label{eq:yLB3-gen}
y_{k} \geq
\left\{
\begin{array}{ll}
\left( \left| N_{\alpha} \right|u^{2} / \left( 4b^{2} \right) \right) \left( 2du^{2}/5 \right)^{k-1} & \text{for $k>0$,} \\
\left( \left| N_{\alpha} \right|u^{2} / \left( 4b^{2} \right) \right) \left( 2du^{2}/5 \right)^{\max(0,K-k)} & \text{for $k<0$.}
\end{array}
\right.
\end{equation}

In fact, if $(d,t,u) \neq (5,1,1)$, then we can replace $2du^{2}/5$ by $5du^{2}/8$
and if $N_{\varepsilon}=1$, then we can replace $2du^{2}/5$ by $du^{2}$.
\end{lemma}

\begin{remark-nonum}
The condition $N_{\alpha}<0$ is needed, since if $N_{\alpha}>0$, then
$y_{k}<0$ can occur for $k<0$.

Also, since $du^{2} \geq 5$ holds under the conditions here, the lower bounds in
\eqref{eq:yLB3-gen} are increasing as $|k|$ increases.
\end{remark-nonum}

\begin{proof}
(a) From \eqref{eq:yk-defn}, we can write
\[
y_{k}= \frac{\alpha\varepsilon^{2k}-\overline{\alpha}\,\overline{\varepsilon}^{2k}}{2\sqrt{d}},
\]
where $\overline{\varepsilon}=\left( t-u\sqrt{d} \right)/2$.
Since $N_{\alpha}=\alpha\overline{\alpha}<0$ and $\alpha>0$, we have
$\overline{\alpha}<0$. So the inequality~\eqref{eq:yLB1-gen} follows.

\vspace*{1.0mm}

(b) From part~(a), it follows that all the $y_{k}$'s are positive.
Since $\alpha \varepsilon^{2k}$ is an algebraic integer, we have $2y_{k} \in \bbZ$.

Since $d$, $t^{2}$ and $u^{2}$ are all positive integers, a quick search over
small values of $d$, $t$ and $u$ with $t^{2}- du^{2}=\pm 4$ shows that
$t^{2} + du^{2} \geq 6$ (the minimum occurs for $(d,t,u)=(5, \pm 1, \pm 1)$).
From this, \eqref{eq:yk-recurrence} and $y_{k}>0$, we have
$y_{k+1} \geq 3y_{k} - y_{k-1}$ for $k \geq 1$.

From the expression for $y_{1}$ after \eqref{eq:yk-defn} and $t^{2} + du^{2} \geq 6$,
we have $y_{1}>(3/2)b^{2}>b^{2}=y_{0}$. So using induction and
$y_{k+1} \geq 3y_{k} - y_{k-1}$, we find that
$y_{k+1}>y_{k}$ for all $k \geq 0$.

By the definition of $K$, $y_{K}>y_{K+1}$. Also $y_{K}>b^{2}>0$, so we can
proceed as in the case of $k \geq 0$, using
$y_{k-1} \geq 3y_{k} - y_{k+1}$ and induction.

\vspace*{1.0mm}

(c) From \eqref{eq:yPM1}, we find that
\begin{equation}
\label{eq:rel-sqrs}
4b^{2}y_{1}=b^{4} \left( t^{2}+du^{2} \right) + 2ab^{2}tu
= \left( b^{2}t + au \right)^{2} -N_{\alpha}u^{2}.
\end{equation}

From the second equality in \eqref{eq:rel-sqrs} and $N_{\alpha}<0$,
equation~\eqref{eq:yLB3-gen} follows for $k=1$.

Writing $\varepsilon^{\ell}=t_{\ell}+u_{\ell}\sqrt{d}$ for $\ell \geq 1$ (note
that $u_{1}=u/2$ in the notation of this lemma), we can
show by an easy induction that $\left( u_{\ell} \right)_{\ell \geq 1}$ is an
increasing sequence. Hence, for $k<0$,
\begin{align*}
4b^{2}y_{k} &=b^{4} \left( t_{|2k|}^{2}+du_{|2k|}^{2} \right) - 2ab^{2}t_{|2k|} u_{|2k|}
= \left( b^{2}t_{|2k|} - au_{|2k|} \right)^{2} -N_{\alpha}u_{|2k|}^{2} \\
&\geq \left| N_{\alpha} \right| u^{2}.
\end{align*}

So equation~\eqref{eq:yLB3-gen} holds for $k$ satisfying $K \leq k<0$.

Recalling the recurrence in \eqref{eq:yk-recurrence} for $k \geq 0$ and using the
monotonicity established in part~(b), we obtain
\[
y_{k+1}=\left( du^{2} + 2N_{\varepsilon} \right)y_{k}-y_{k-1}
\geq \left( du^{2} + 2N_{\varepsilon}-1 \right)y_{k}.
\]

If $N_{\varepsilon}=1$, then the stronger form of the desired inequality stated
after equation~\eqref{eq:yLB3-gen} in the lemma holds,
so we need only consider $N_{\varepsilon}=-1$.
We always have $du^{2} \geq 5$ (with equality only if $(d,t,u)=(5,1,1)$.
Otherwise, $du^{2} \geq 8$, which we use to establish the other inequality
after equation~\eqref{eq:yLB3-gen} in the lemma). Hence the desired inequality holds.

Using the recurrence
$y_{k-1}=\left( du^{2} + 2N_{\varepsilon} \right)y_{k}-y_{k+1}$,
we obtain equation~\eqref{eq:yLB3-gen} for $k \leq K$ in the same way.
\end{proof}

We will also need to know when $K$ in Lemma~\ref{lem:Y-LB} is not equal to $-1$
for $b=1$.

\begin{lemma}
\label{lem:K-value}
Let the $y_{k}$'s be defined by $\eqref{eq:yk-defn}$ with the notation and assumptions
there. Suppose that $b=1$ and $N_{\alpha}<0$.
 
We have $y_{-1}>1$, except for\\
{\rm (i)} $a \geq 1$, $d=a^{2}+4$, $t=a$ and $u=1$, where $\alpha=2\varepsilon$ and $N_{\alpha}=-4$, \\
{\rm (ii)} $a \geq 1$, $d=a^{2}+1$, $t=2a$ and $u=2$, where $\alpha=\varepsilon$ and $N_{\alpha}=-1$.

In these cases, we have $y_{-2}=a^{2}+1$ and $y_{-2}=4a^{2}+1$, respectively.
So in each of these cases, $K=-2$.
\end{lemma}

\begin{proof}
From \eqref{eq:yLB3-gen}, we have $y_{1} \geq \left| N_{\alpha} \right| u^{2}/4$.
So if $y_{-1}=1$, we must have either (1) $N_{\alpha}=-4$ and $u=1$, or (2)
$N_{\alpha}=-1$ and $u=2$.

In case~(1), we have $d=t^{2} \pm 4$ and $a^{2}-d=-4$. Substituting the
expression for $d$ into the second equation, we get $a^{2}-t^{2}=-4 \pm 4=0,-8$.

$a^{2}-t^{2}=-8$ can only occur for $a=1$ and $t=3$ (since we are assuming that
$a$ and $t$ are positive integers). But in this case, $b^{2}t-au=2$, so
$4y_{-1}=\left( b^{2}t-au \right)^{2}-N_{\alpha}u^{2}=8$.

Otherwise, we have $t=a$. So $d=a^{2}+4$, $\alpha=2\varepsilon$ and $N_{\varepsilon}=-1$.
Here $b^{2}t-au=0$, so $y_{-1}=1$.

In case~(2), we have $d= \left( t^{2} \pm 4 \right)/4$ and $a^{2}-d=-1$.
We proceed similarly, obtaining $4a^{2}-t^{2}=-8$ (which is never possible)
when $N_{\varepsilon}=1$ and $4a^{2}-t^{2}=0$ when $N_{\varepsilon}=-1$.
Here we have $t=2a$ and $d=a^{2}+1$, so $\alpha=\varepsilon$.
Once again, $b^{2}t-au=0$, so $y_{-1}=1$.

A direct calculation, done using Maple, provides the values of $y_{-2}$.
%
%
\end{proof}

\subsection{Gap Principle}

In Lemma~\ref{lem:gap} below, we establish a gap principle separating distinct
squares in the sequence \eqref{eq:yk-defn}. We first need the following technical
lemma to help us prove our gap principle.

\begin{lemma}
\label{lem:omega-bnd}
Let $\omega=e^{i\theta}$ with $-\pi < \theta \leq \pi$.

\noindent
{\rm (a)} Put $\omega^{1/4}=e^{i\theta/4}$. If
\[
0 < \left| \omega^{1/4} - z \right| <c_{1},
\]
for some $z \in \bbC$ with $|z|=1$ and $0 < c_{1} < \sqrt{2}$, then
\[
\left| \omega - z^{4} \right|
> c_{2} \left| \omega^{1/4} - z \right|,
\]
where $c_{2}=\left( 2-c_{1}^{2} \right) \sqrt{4-c_{1}^{2}}$.

\noindent
{\rm (b)} If
\[
0 < \left| \omega - z^{4} \right| <c_{0},
\]
for some $z \in \bbC$ with $|z|=1$ and $0 < c_{0} \leq 2$, then
\[
0 < \left| \omega^{1/4} - z \right| <c_{3},
\]
for some choice of $\omega^{1/4}$, where $c_{3}$ is the smallest positive real
root of $x^{8}-8x^{6}+20x^{4}-16x^{2}+c_{0}^{2}$.
\end{lemma}

\begin{proof}
(a) We can write
\[
\left| \omega - z^{4} \right|
= \left| \omega^{1/4} - z \right|
  \times \prod_{k=1}^{3} \left| \omega^{1/4} - e^{2\pi ik/4}z \right|.
\]

Multiplying by $\omega^{-3/4}$ and expanding the resulting expression, the
product above equals
\[
\prod_{k=1}^{3} \left| e^{2\pi ik/4-i\theta/4}z-1 \right|
= \left| e^{3i\varphi}+e^{2i\varphi}+e^{i\varphi}+1 \right|,
\]
for $-\pi < \varphi \leq \pi$ such that $e^{i\varphi}=e^{-i\theta/4}z$. Squaring
this quantity and simplifying, we obtain
\begin{equation}
\label{eq:trig-a}
8\cos^{2}(\varphi) \left( \cos(\varphi)+1 \right).
\end{equation}

If $\left| \omega^{1/4} - z \right|=\left| 1 - \omega^{-1/4}z \right|=c_{1}$,
we have $2-2\cos(\varphi)=c_{1}^{2}$ and, substituting this expression for
$\cos(\varphi)$ into \eqref{eq:trig-a}, we find that
$\left| \omega - z^{4} \right| = c_{2} \left| \omega^{1/4} - z \right|$.
Since $c_{2}$ is a decreasing function
of $c_{1}$, part~(a) of the lemma holds
(i.e., if $\left| \omega^{1/4} - z \right|<c_{1}$).

(b) Using the same argument as in part~(a), but supposing that
$\left| \omega - z^{4} \right|=\left| 1 - \omega^{-1}z^{4} \right|=c_{0}$,
then we have $2-2\cos(4\varphi)=c_{0}^{2}$. Thus
\begin{align*}
\left| \omega^{1/4}-z \right|^{2}
& =\left( \frac{\left| \omega-z^{4} \right|}{\left| e^{3i\varphi}+e^{2i\varphi}+e^{i\varphi}+1 \right|} \right)^{2}
=\frac{c_{0}^{2}}{8\cos^{2}(\varphi)(\cos(\varphi)+1)} \\
& =\frac{2-2\cos (4 \varphi)}{8\cos^{2}(\varphi)(\cos(\varphi)+1)}.
\end{align*}

Since $\cos(4\varphi)=8\cos^{4}(\varphi)-8\cos^{2}(\varphi)+1$, we have
$16\cos^{4}(\varphi)-16\cos^{2}(\varphi)+c_{0}^{2}=0$ and using this relation,
we find that
$c_{0}^{2}/ \left( 8\cos^{2}(\varphi)(\cos(\varphi)+1) \right)$ is a root of
the polynomial $x^{4}-8x^{3}+20x^{2}-16x+c_{0}^{2}$. Part~(b) now follows and
follows with inequalities too, since the smallest positive real root of $x^{4}-8x^{3}+20x^{2}-16x+c_{0}^{2}$
grows with $c_{0}$.
\end{proof}

\begin{lemma}
\label{lem:gap}
Let the $y_{k}$'s be defined as in $\eqref{eq:yk-defn}$ with $N_{\alpha}<0$.

{\rm (a)} Suppose that $-N_{\alpha}$ is a square. If $y_{i}$ and $y_{j}$ are
distinct squares with $i,j \neq 0$ and
$y_{j}>y_{i} \geq \max \left( 4\sqrt{\left| N_{\alpha} \right|/d}, b^{2}\left| N_{\alpha} \right|/d \right)$,
then
\[
y_{j} > 57.32 \left( \frac{d}{b^{2}\left| N_{\alpha} \right|} \right)^{2} y_{i}^{3}.
\]

{\rm (b)} Suppose $-N_{\alpha}$ is not a square. If $y_{k_{1}}$, $y_{k_{2}}$ and
$y_{k_{3}}$ are three distinct squares with $k_{1}, k_{2}, k_{3} \neq 0$ and
$y_{k_{3}}>y_{k_{2}}>y_{k_{1}} \geq 4\sqrt{\left| N_{\alpha} \right|/d}$,
then there exist distinct $i,j \in \{ k_{1},k_{2}, k_{3} \}$ such that
\[
y_{j} > 15.36 \left( \frac{b^{2}d}{f_{i}f_{j} \left| N_{\alpha} \right|} \right)^{2} y_{i}^{3},
\]
where $f_{i}$ and $f_{j}$ are the values of $f$ in Proposition~$\ref{prop:quad-rep}$
associated with $y_{i}$ and $y_{j}$, respectively.

If $y_{j}>y_{i} \geq \max \left( 4\sqrt{\left| N_{\alpha} \right|/d}, 4.27b^{2}\left| N_{\alpha} \right|^{2}/d \right)$,
then we can replace $15.36$ with $182$.
\end{lemma}

\begin{proof}
We start by considering the two parts together and what is common to their proofs.

Since we have assumed that $i \neq 0$ and $j \neq 0$, we can apply Proposition~\ref{prop:quad-rep}
to show that there are integers $f_{j}$, $r_{j}$ and $s_{j}$,
and choices of signs such that
\begin{align}
\label{eq:gap-rel1}
\pm f_{j}^{2} \left( x_{j} + N_{\varepsilon^{j}}\sqrt{N_{\alpha}} \right)
&= \left( a+\sqrt{N_{\alpha}} \right) \left( r_{j} + s_{j}\sqrt{\core \left( N_{\alpha} \right)} \right)^{4}
\quad \text{and} \nonumber \\
f_{j}\sqrt{y_{j}} &= b \left( r_{j}^{2}-\core \left( N_{\alpha} \right)s_{j}^{2} \right),
\end{align}
where $f_{j}$ satisfies $f_{j} | \left( 4b^{2}\core \left( \left| N_{\alpha} \right| \right) \right)$
and $f_{j}^{2}< 16b^{4} \core \left( \left| N_{\alpha} \right| \right)$. Recall
that the other terms in the relationships for $f_{j}$ in Proposition~\ref{prop:quad-rep}
are not present here since $N_{\varepsilon^{j}}= \pm 1$.

For any two distinct squares among the $y_{k}$'s, say $y_{i}$ and $y_{j}$ with
$i \neq 0$ and $j \neq 0$, we will assume that the $\pm$ on the left-hand side
of \eqref{eq:gap-rel1} is always $+$. That is,
\begin{align}
\label{eq:3.2}
f_{i}^{2} \left( x_{i} + N_{\varepsilon^{i}}\sqrt{N_{\alpha}} \right)
&= \left( a+\sqrt{N_{\alpha}} \right) \left( r_{i} + s_{i}\sqrt{\core \left( N_{\alpha} \right)} \right)^{4}
\text{ and } \\
f_{j}^{2} \left( x_{j} + N_{\varepsilon^{j}}\sqrt{N_{\alpha}} \right)
&= \left( a+\sqrt{N_{\alpha}} \right) \left( r_{j} + s_{j}\sqrt{\core \left( N_{\alpha} \right)} \right)^{4},
\nonumber
\end{align}
as the argument for the other cases is exactly the same.
Subtracting the complex conjugate of one of these equations from the equation
itself, we obtain
\begin{align}
\label{eq:3.3}
&  \left( a+\sqrt{N_{\alpha}} \right) \left( r_{j} + s_{j}\sqrt{\core \left( N_{\alpha} \right)} \right)^{4}
      - \left( a-\sqrt{N_{\alpha}} \right) \left( r_{j} - s_{j}\sqrt{\core \left( N_{\alpha} \right)} \right)^{4} \\
= &
2i \Imag \left( f_{j}^{2} \left( x_{j} + N_{\varepsilon^{j}}\sqrt{N_{\alpha}} \right) \right)
= \pm 2f_{j}^{2}\sqrt{\left| N_{\alpha} \right|}, \nonumber
\end{align}
and the analogous equation with the index $j$ replaced by $i$ also holds.

Putting $\omega = \left( a-\sqrt{N_{\alpha}} \right) / \left( a+\sqrt{N_{\alpha}} \right)$,
by \eqref{eq:3.3} and the relationship for $f_{j}\sqrt{y_{j}}$ in Proposition~\ref{prop:quad-rep}(a),
we have
\begin{align}
\label{eq:3.4}
\left| \omega - \left( \frac{r_{j}+s_{j}\sqrt{\core \left( N_{\alpha} \right)}}
{r_{j}-s_{j}\sqrt{\core \left( N_{\alpha} \right)}} \right)^{4} \right|
& =\left| \frac{\pm 2f_{j}^{2}\sqrt{ \left| N_{\alpha} \right|}}
{\left( a+\sqrt{N_{\alpha}} \right)\left( r_{j}-s_{j}\sqrt{\core \left( N_{\alpha} \right)} \right)^{4}} \right| \\
&= \frac{2b^{2}\sqrt{\left| N_{\alpha} \right|}}{\sqrt{a^{2}+\left| N_{\alpha} \right|} \, y_{j}}
\leq \frac{1}{2}, \nonumber
\end{align}
the last inequality holds because $a^{2}+\left| N_{\alpha} \right|=db^{4}$ and
$y_{j} \geq 4\sqrt{\left| N_{\alpha} \right|/d}$.

Let $\zeta_{4}^{(j)}$ be the $4$-th root of unity such that
\begin{equation}
\label{eq:zetaj-defn}
\left| \omega^{1/4} - \zeta_{4}^{(j)} \frac{r_{j}+s_{j}\sqrt{\core \left( N_{\alpha} \right)}}{r_{j}-s_{j}\sqrt{\core \left( N_{\alpha} \right)}} \right|
= \min_{0 \leq k \leq 3} \left| \omega^{1/4} - e^{2k \pi i/4} \frac{r_{j}+s_{j}\sqrt{\core \left( N_{\alpha} \right)}}{r_{j}-s_{j}\sqrt{\core \left( N_{\alpha} \right)}} \right|.
\end{equation}

From \eqref{eq:3.4}, we can apply Lemma~\ref{lem:omega-bnd}(b) with $c_{0}=1/2+0.0001$
to obtain
\[
\left| \omega^{1/4} - \zeta_{4}^{(j)} \frac{r_{j}+s_{j}\sqrt{\core \left( N_{\alpha} \right)}}
{r_{j}-s_{j}\sqrt{\core \left( N_{\alpha} \right)}} \right|
<0.1263.
\]

Applying Lemma~\ref{lem:omega-bnd}(a) with $c_{1}=0.1263$, we obtain
\[
\left| \omega - \left( \frac{r_{j}+s_{j}\sqrt{\core \left( N_{\alpha} \right)}}
{r_{j}-s_{j}\sqrt{\core \left( N_{\alpha} \right)}} \right)^{4} \right|
>3.96\left| \omega^{1/4} - \zeta_{4}^{(j)} \frac{r_{j}+s_{j}\sqrt{\core \left( N_{\alpha} \right)}}
{r_{j}-s_{j}\sqrt{\core \left( N_{\alpha} \right)}} \right|
\]
and combining this with the equalities in \eqref{eq:3.4} yields
\begin{equation}
\label{eq:3.6}
\left| \omega^{1/4} - \zeta_{4}^{(j)} \frac{r_{j}+s_{j}\sqrt{\core \left( N_{\alpha} \right)}}
{r_{j}-s_{j}\sqrt{\core \left( N_{\alpha} \right)}} \right|
<0.5051 \frac{b^{2} \sqrt{\left| N_{\alpha} \right|}}{\sqrt{a^{2}+\left| N_{\alpha} \right|}} \, \frac{1}{y_{j}}.
\end{equation}

In the same way, we define $\zeta_{4}^{(i)}$ for any square $y_{i}$ with $i \neq 0$
in our sequence and \eqref{eq:3.6} also holds with $j$ replaced by $i$. Hence, for
any two distinct squares, $y_{i},y_{j} \geq 4\sqrt{\left| N_{\alpha} \right|/d}$, among the $y_{k}$'s, we have
\begin{align}
\label{eq:3.7}
      & \left| \zeta_{4}^{(i)} \frac{r_{i}+s_{i}\sqrt{\core \left( N_{\alpha} \right)}}
{r_{i}-s_{i}\sqrt{\core \left( N_{\alpha} \right)}}
         - \zeta_{4}^{(j)} \frac{r_{j}+s_{j}\sqrt{\core \left( N_{\alpha} \right)}}
{r_{j}-s_{j}\sqrt{\core \left( N_{\alpha} \right)}} \right| \\
 \leq & \left| \omega^{1/4} - \zeta_{4}^{(i)} \frac{r_{i}+s_{i}\sqrt{\core \left( N_{\alpha} \right)}}
{r_{i}-s_{i}\sqrt{\core \left( N_{\alpha} \right)}} \right|
         + \left| \omega^{1/4} - \zeta_{4}^{(j)} \frac{r_{j}+s_{j}\sqrt{\core \left( N_{\alpha} \right)}}
{r_{j}-s_{j}\sqrt{\core \left( N_{\alpha} \right)}} \right| \nonumber \\
   <  & 0.5051b^{2} \sqrt{\frac{\left| N_{\alpha} \right|}{a^{2}+\left| N_{\alpha} \right|}}
         \left( \frac{1}{y_{i}}+\frac{1}{y_{j}} \right). \nonumber
\end{align}

Next we obtain a lower bound for this same quantity. We first show that it
cannot be zero.

If
\[
\zeta_{4}^{(i)} \frac{r_{i}+s_{i}\sqrt{\core \left( N_{\alpha} \right)}}
{r_{i}-s_{i}\sqrt{\core \left( N_{\alpha} \right)}}
=\zeta_{4}^{(j)} \frac{r_{j}+s_{j}\sqrt{\core \left( N_{\alpha} \right)}}
{r_{j}-s_{j}\sqrt{\core \left( N_{\alpha} \right)}},
\]
then from our expressions for $f_{i}\sqrt{y_{i}}$ and $f_{j}\sqrt{y_{j}}$ from
Proposition~\ref{prop:quad-rep} we have
\[
\frac{\left( r_{i}+s_{i}\sqrt{\core \left( N_{\alpha} \right)} \right)^{4}}{f_{i}^{2}y_{i}}
= \pm \frac{\left( r_{j}+s_{j}\sqrt{\core \left( N_{\alpha} \right)} \right)^{4}}{f_{j}^{2}y_{j}}.
\]

From \eqref{eq:3.2}, it follows that
\[
\left( x_{i} + N_{\varepsilon^{i}}\sqrt{N_{\alpha}} \right) y_{j}
= \pm \left( x_{j}  + N_{\varepsilon^{j}}\sqrt{N_{\alpha}} \right) y_{i}.
\]

Comparing the imaginary parts of both sides of this equation, we find that
$y_{i}=y_{j}$, but this contradicts our assumption that $y_{j}>y_{i}$. Hence
the left-hand side of \eqref{eq:3.7} cannot be $0$.

Let $x+y\sqrt{\core \left( N_{\alpha} \right)}= \left( r_{i}-s_{i}\sqrt{\core \left( N_{\alpha} \right)} \right)
\left( r_{j}+s_{j}\sqrt{\core \left( N_{\alpha} \right)} \right)$.
We can write
\begin{align}
\label{eq:gap-lb}
&\zeta_{4}^{(i)} \frac{r_{i}+s_{i}\sqrt{\core \left( N_{\alpha} \right)}}{r_{i}-s_{i}\sqrt{\core \left( N_{\alpha} \right)}}
- \zeta_{4}^{(j)} \frac{r_{j}+s_{j}\sqrt{\core \left( N_{\alpha} \right)}}{r_{j}-s_{j}\sqrt{\core \left( N_{\alpha} \right)}} \\
= & \frac{2\zeta_{4}^{(i)}x-\left( \zeta_{4}^{(i)}+\zeta_{4}^{(j)} \right) \left( x+y\sqrt{\core \left( N_{\alpha} \right)} \right)}
{\left( r_{i}-s_{i}\sqrt{\core \left( N_{\alpha} \right)} \right) \left( r_{j}-s_{j}\sqrt{\core \left( N_{\alpha} \right)} \right)}. \nonumber
\end{align}

The numerator on the right-hand side is
\begin{equation}
\label{eq:num-exp}
\begin{array}{rl}
-2\zeta_{4}^{(i)}y\sqrt{\core \left( N_{\alpha} \right)} & \text{if $\zeta_{4}^{(j)}=\zeta_{4}^{(i)}$,} \\
2\zeta_{4}^{(i)}x                  & \text{if $\zeta_{4}^{(j)}=-\zeta_{4}^{(i)}$,} \\
\zeta_{4}^{(i)} \left( 1-\sqrt{-1} \right) \left( x-y\sqrt{-\core \left( N_{\alpha} \right)} \right) & \text{if $\zeta_{4}^{(j)}=\sqrt{-1}\zeta_{4}^{(i)}$ and} \\
\zeta_{4}^{(i)} \left( 1+\sqrt{-1} \right) \left( x+y\sqrt{-\core \left( N_{\alpha} \right)} \right) & \text{if $\zeta_{4}^{(j)}=-\sqrt{-1}\zeta_{4}^{(i)}$.}
\end{array}
\end{equation}
Here we use $\sqrt{-1}$ to denote $\exp \left( 2\pi i/4 \right)$.

Using \eqref{eq:fy-rel}, we find that
\begin{equation}
\label{eq:abs}
\left| r_{i}-s_{i}\sqrt{\core \left( N_{\alpha} \right)} \right|
\left| r_{j}-s_{j}\sqrt{\core \left( N_{\alpha} \right)} \right|
=\frac{\sqrt{f_{i}f_{j}}\left( y_{i}y_{j} \right)^{1/4}}{b}.
\end{equation}

At this point, our proofs of the two parts of the lemma separate.

\vspace*{1.0mm}

(a) Suppose that $-N_{\alpha}$ is a square. From \eqref{eq:num-exp}, we see that
at least one of $1 \pm \sqrt{-1}$ always divides the numerator of the right-hand side
of \eqref{eq:gap-lb} and so
\begin{align*}
& \left| \zeta_{4}^{(i)} \frac{r_{i}+s_{i}\sqrt{\core \left( N_{\alpha} \right)}}{r_{i}-s_{i}\sqrt{\core \left( N_{\alpha} \right)}}
       - \zeta_{4}^{(j)} \frac{r_{j}+s_{j}\sqrt{\core \left( N_{\alpha} \right)}}{r_{j}-s_{j}\sqrt{\core \left( N_{\alpha} \right)}} \right| \\
\geq & \frac{|1+\sqrt{-1}|}
          {\left| \left( r_{i}-s_{i}\sqrt{\core \left( N_{\alpha} \right)} \right)\left( r_{j}-s_{j}\sqrt{\core \left( N_{\alpha} \right)} \right) \right|}
= \frac{\sqrt{2} \, b}{\sqrt{f_{i}f_{j}} \left( y_{i}y_{j} \right)^{1/4}},
\end{align*}
using \eqref{eq:abs} above.

Furthermore, from Proposition~\ref{prop:quad-rep}(b), we know that $f_{i},f_{j}|b^{2}$,
since $-N_{\alpha}$ is a square. So
\[
\left| \zeta_{4}^{(i)} \frac{r_{i}+s_{i}\sqrt{\core \left( N_{\alpha} \right)}}{r_{i}-s_{i}\sqrt{\core \left( N_{\alpha} \right)}}
       - \zeta_{4}^{(j)} \frac{r_{j}+s_{j}\sqrt{\core \left( N_{\alpha} \right)}}{r_{j}-s_{j}\sqrt{\core \left( N_{\alpha} \right)}} \right|
\geq \frac{\sqrt{2}}{b\left( y_{i}y_{j} \right)^{1/4}}.
\]

Combining this with \eqref{eq:3.7}, we have
\begin{equation}
\label{eq:gapLB-a}
\frac{\sqrt{2}}{b \left( y_{i}y_{j} \right)^{1/4}}
< 0.5051 b^{2} \sqrt{\frac{\left| N_{\alpha} \right|}{a^{2}+\left| N_{\alpha} \right|}}
\left( \frac{1}{y_{i}} + \frac{1}{y_{j}} \right).
\end{equation}

From $y_{j}>y_{i}$, it follows that $1/y_{i}+1/y_{j}<2/y_{i}$ and so
\[
y_{j}> \frac{4}{1.0102^{4} b^{12}}
\left( \frac{a^{2}+\left| N_{\alpha} \right|}{\left| N_{\alpha} \right|} \right)^{2} y_{i}^{3}.
\]

Since $N_{\alpha}=a^{2}-b^{4}d<0$, we have
\begin{equation}
\label{eq:gap-a1}
y_{j}
> \frac{3.84d^{2}}{b^{4}\left| N_{\alpha} \right|^{2}} y_{i}^{3}.
\end{equation}

We can use this gap principle to improve its constant term. Combining
\eqref{eq:gap-a1} with $y_{i} \geq b^{2}\left| N_{\alpha} \right|/d$,
we obtain
\[
y_{j}
>\frac{3.84d^{2}}{b^{4}\left| N_{\alpha} \right|^{2}}
\left( \frac{b^{2}\left| N_{\alpha} \right|}{d} \right)^{2} y_{i}
=3.84y_{i}.
\]

Applying this to \eqref{eq:gapLB-a} yields
\[
\frac{\sqrt{2}}{b\left( y_{i}y_{j} \right)^{1/4}}
< 0.5051 b^{2} \sqrt{\frac{\left| N_{\alpha} \right|}{a^{2}+\left| N_{\alpha} \right|}}
\frac{1+1/3.84}{y_{i}}.
\]

This implies that
\[
y_{j}> \frac{24.34}{b^{12}} \left( \frac{a^{2}+\left| N_{\alpha} \right|}{\left| N_{\alpha} \right|} \right)^{2} y_{i}^{3}
=\frac{24.34d^{2}}{b^{4}\left| N_{\alpha} \right|^{2}} y_{i}^{3}.
\]

Repeating the process from \eqref{eq:gap-a1} onwards with $3.84$ replaced by
$24.34$, we obtain
\[
y_{j}> \frac{52.31d^{2}}{b^{4}\left| N_{\alpha} \right|^{2}} y_{i}^{3}.
\]

Repeating it again with $52.31$ instead of $3.84$, we improve the constant to $56.97$.
Repeating it one final time with $56.97$, we obtain the inequality in part~(a).

\vspace*{2.0mm}

(b) We again use \eqref{eq:gap-lb} and \eqref{eq:num-exp}. Notice that
$x+y\sqrt{-\core \left( N_{\alpha} \right)}$ on the left-hand side of
\eqref{eq:num-exp} for $\zeta_{4}^{(j)}=\pm \sqrt{-1}\zeta_{4}^{(i)}$ can be as
small as $1/ \left( 2\sqrt{-\core \left( N_{\alpha} \right)}y \right)$ when
$-N_{\alpha}$ is not a square, so the proof of part~(a) breaks down here.

Among the three values $\zeta_{4}^{(k_{1})}$, $\zeta_{4}^{(k_{2})}$ and
$\zeta_{4}^{(k_{3})}$, there must be at least one pair such that one member of
the pair is $\pm 1$ times the other. Choose any such pair and let $i$ and $j$
be the associated indices, ordered so that $y_{i}<y_{j}$.

Thus,
\[
\left| \zeta_{4}^{(i)} \frac{r_{i}+s_{i}\sqrt{\core \left( N_{\alpha} \right)}}{r_{i}-s_{i}\sqrt{\core \left( N_{\alpha} \right)}}
       - \zeta_{4}^{(j)} \frac{r_{j}+s_{j}\sqrt{\core \left( N_{\alpha} \right)}}{r_{j}-s_{j}\sqrt{\core \left( N_{\alpha} \right)}} \right|
\geq \frac{2b}{\sqrt{f_{i}f_{j}} \left( y_{i}y_{j} \right)^{1/4}}.
\]

The argument is the same as in the proof of part~(a) except we have an extra
factor of $\sqrt{2}$ on the right-hand side of this inequality. Thus
\begin{equation}
\label{eq:gapLB-b}
\frac{2b}{\sqrt{f_{i}f_{j}} \left( y_{i}y_{j} \right)^{1/4}}
< 0.5051b^{2} \sqrt{\frac{\left| N_{\alpha} \right|}{a^{2}+\left| N_{\alpha} \right|}}
\left( \frac{1}{y_{i}} + \frac{1}{y_{j}} \right).
\end{equation}

As in the proof of part~(a), this gives
\[
y_{j}> 15.36 \left( \frac{b^{2}d}{f_{i}f_{j} \left| N_{\alpha} \right|} \right)^{2} y_{i}^{3},
\]
which establishes the first lower bound for $y_{j}$ in part~(b).

If $y_{i}>f_{i}f_{j} \left| N_{\alpha} \right|\sqrt{1.09/15.36}/ \left( b^{2}d \right)$,
then this lower bound for $y_{j}$ yields $y_{j}>1.09y_{i}$. Applying this to
\eqref{eq:gapLB-b}, we obtain
\[
y_{j}>18.25 \left( \frac{b^{2}d}{f_{i}f_{j} \left| N_{\alpha} \right|} \right)^{2} y_{i}^{3}.
\]

Repeating this process eight more times yields
\[
y_{j}>182 \left( \frac{b^{2}d}{f_{i}f_{j} \left| N_{\alpha} \right|} \right)^{2} y_{i}^{3}.
\]

Since $f_{i}f_{j} \leq \left( 4b^{2} \right)^{2} \left| N_{\alpha} \right|$,
if $y_{i} \geq 4.27b^{2}\left| N_{\alpha} \right|^{2}/d$, then
$y_{i}>f_{i}f_{j} \left| N_{\alpha} \right|\sqrt{1.09/15.36}/ \left( b^{2}d \right)$,
as required.
\end{proof}

\subsection{Miscellaneous Lemmas}

Here we collect some results that we will need for bounding quantities that
arise in the proof of the main result in the following section (Section~\ref{sect:prop-11}).

\begin{lemma}
\label{lem:zeta}
Let the $y_{k}$'s be defined as in $\eqref{eq:yk-defn}$ with the notation and assumptions
there. Suppose that $N_{\alpha}<0$.

{\rm (a)} Let $y_{k} \geq 4\sqrt{\left| N_{\alpha} \right|/d}$ be a square and put
\[
\omega_{k}=\left( x_{k} + N_{\varepsilon^{k}}\sqrt{N_{\alpha}} \right)
/\left( x_{k} - N_{\varepsilon^{k}}\sqrt{N_{\alpha}} \right)
=e^{i\varphi_{k}}
\]
with $-\pi<\varphi_{k} \leq \pi$.
Then
\[
\left| \varphi_{k} \right| <\frac{2.29\sqrt{\left| N_{\alpha} \right|}}{\left| x_{k} \right|}
<0.6.
\]

{\rm (b)} Let $y_{k}$ and $y_{\ell}$ be two squares with $k,\ell \neq 0$,
$y_{\ell}>y_{k} \geq 4\sqrt{\left| N_{\alpha} \right|/d}$,
$\omega_{k}$ be as in part~{\rm(a)} and put
\[
x+y\sqrt{\core \left( N_{\alpha} \right)}
=\left(  r_{k}-s_{k}\sqrt{\core \left( N_{\alpha} \right)} \right)
\left( r_{\ell}+s_{\ell}\sqrt{\core \left( N_{\alpha} \right)} \right)
\]
with
$r_{k}$, $r_{\ell}$, $s_{k}$ and $s_{\ell}$ as in Proposition~$\ref{prop:quad-rep}$.
Furthermore, suppose that the quantities $\zeta_{4}^{(k)}$ and $\zeta_{4}^{(\ell)}$
defined in \eqref{eq:zetaj-defn} in the proof of Lemma~$\ref{lem:gap}$ satisfy
$\zeta_{4}^{(k)}= \pm \zeta_{4}^{(\ell)}$. Then
\[
\min_{0 \leq j \leq 3} \left| \omega_{k}^{1/4} - \zeta_{4}^{j} \frac{x+y\sqrt{\core \left( N_{\alpha} \right)}}{x-y\sqrt{\core \left( N_{\alpha} \right)}} \right|
\]
occurs for either $j=0$ or $j=2$, where $\zeta_{4}$ is a primitive $4$-th root
of unity.
\end{lemma}

\begin{proof}
(a) We can write $\omega_{k} = \left( \left| x_{k} \right| \pm \sqrt{N_{\alpha}} \right)^{2}/ \left( x_{k}^{2}+ \left| N_{\alpha} \right| \right)$,
so with $\omega_{k}=e^{i\varphi_{k}}$, we have
$\left| \tan \left( \varphi_{k} \right) \right| = \left| 2x_{k}\sqrt{\left| N_{\alpha} \right|}/\left( x_{k}^{2}+N_{\alpha} \right) \right|$.
From \eqref{eq:yk-defn} and $N_{\alpha}>-b^{4}d$, we have
\[
x_{k}^{2}+N_{\alpha} = dy_{k}^{2}+2N_{\alpha}
=dy_{k}^{2} \left( 1 + \frac{2N_{\alpha}}{dy_{k}^{2}} \right)
\geq 0.875dy_{k}^{2} > 0.875x_{k}^{2},
\]
since $y_{k} \geq 4\sqrt{\left| N_{\alpha} \right|/d}$ and
$x_{k}^{2}-dy_{k}^{2}=N_{\alpha}<0$.

From $\left| \varphi_{k} \right| \leq \left| \tan \left( \varphi_{k} \right) \right|$, we have
\[
\left| \varphi_{k} \right|
\leq \left| \frac{2\sqrt{\left| N_{\alpha} \right|}x_{k}}{x_{k}^{2}+N_{\alpha}} \right|
< \left| \frac{2\sqrt{\left| N_{\alpha} \right|}x_{k}}{0.875x_{k}^{2}} \right|
<\frac{2.29\sqrt{\left| N_{\alpha} \right|}}{\left| x_{k} \right|}.
\]

Since $y_{k} \geq 4\sqrt{\left| N_{\alpha} \right|/d}$, we have
$dy_{k}^{2} \geq 16\left| N_{\alpha} \right|$
and so
\[
x_{k}^{2}=dy_{k}^{2}+N_{\alpha} \geq 15\left| N_{\alpha} \right|.
\]

Combining this with the inequality above it yields
$\left| \varphi_{k} \right|<2.29\sqrt{\left| N_{\alpha} \right|}/\left| x_{k} \right|
<2.29/\sqrt{15}<0.6$.

(b) Recall that from \eqref{eq:3.3} in the proof of Lemma~\ref{lem:gap}, we can write
\begin{align*}
& \left( a+\sqrt{N_{\alpha}} \right) \left( r_{i} + s_{i} \sqrt{\core \left( N_{\alpha} \right)} \right)^{4}
- \left( a-\sqrt{N_{\alpha}} \right) \left( r_{i} - s_{i} \sqrt{\core \left( N_{\alpha} \right)} \right)^{4} \\
& = \pm 2 f_{i}^{2}\sqrt{N_{\alpha}}
\end{align*}
for $i=k,\ell$.
There are other
cases according to the signs in the result in Proposition~\ref{prop:quad-rep},
but the argument for them is identical to the argument that follows below.

Using this relationship for $i=\ell$ and using our expressions in Proposition~\ref{prop:quad-rep}
for $x_{k} + N_{\varepsilon^{k}}\sqrt{N_{\alpha}}$ and $f_{k}\sqrt{y_{k}}$, we have
\begin{align*}
&
\left( x_{k} + N_{\varepsilon^{k}}\sqrt{N_{\alpha}} \right) \left( r_{k} - s_{k} \sqrt{\core\left( N_{\alpha} \right)} \right)^{4}
\left( r_{\ell} + s_{\ell} \sqrt{\core\left( N_{\alpha} \right)} \right)^{4} \\
& - \left( x_{k} - N_{\varepsilon^{k}}\sqrt{N_{\alpha}} \right) \left( r_{k} + s_{k} \sqrt{\core\left( N_{\alpha} \right)} \right)^{4}
\left( r_{\ell} - s_{\ell} \sqrt{\core\left( N_{\alpha} \right)} \right)^{4} \\
= &
\frac{\left( r_{k}^{2} - \core\left( N_{\alpha} \right)s_{k}^{2} \right)^{4}}
{f_{k}^{2}}
\left[ \left( a+ \sqrt{N_{\alpha}} \right) \left( r_{\ell} + s_{\ell}  \sqrt{\core\left( N_{\alpha} \right)} \right)^{4} \right. \\
& \hspace*{33.0mm} \left. - \left( a- \sqrt{N_{\alpha}} \right) \left( r_{\ell} - s_{\ell} \sqrt{\core\left( N_{\alpha} \right)} \right)^{4} \right] \\
= &
\frac{f_{k}^{2}y_{k}^{2}}{b^{4}}
\left[ \left( a+ \sqrt{N_{\alpha}} \right) \left( r_{\ell} + s_{\ell}  \sqrt{\core\left( N_{\alpha} \right)} \right)^{4} \right. \\
& \hspace*{10.0mm} \left. - \left( a- \sqrt{N_{\alpha}} \right) \left( r_{\ell} - s_{\ell} \sqrt{\core\left( N_{\alpha} \right)} \right)^{4} \right].
\end{align*}

Applying equation~\eqref{eq:3.3} with $j=\ell$ to the last expression and with
$x+y \sqrt{\core \left( N_{\alpha} \right)}$
as defined in the statement of this lemma, we have
\begin{align}
\label{eq:27}
    |f(x,y)|
= & \left| \left( x_{k} + N_{\varepsilon^{k}}\sqrt{N_{\alpha}} \right) \left( x+y \sqrt{\core\left( N_{\alpha} \right)} \right)^{4} \right. \\
  & \left. - \left( x_{k} - N_{\varepsilon^{k}}\sqrt{N_{\alpha}} \right) \left( x-y \sqrt{\core\left( N_{\alpha} \right)} \right)^{4} \right| \nonumber \\
= & \frac{2f_{k}^{2}f_{\ell}^{2}y_{k}^{2}}{b^{4}} \sqrt{\left| N_{\alpha} \right|}. \nonumber
\end{align}

Let $\zeta_{4}$ be the $4$-th root of unity satisfying
\[
\left| \omega_{k}^{1/4} - \zeta_{4} \frac{x-y\sqrt{\core\left( N_{\alpha} \right)}}{x+y\sqrt{\core\left( N_{\alpha} \right)}} \right|
= \min_{0 \leq j \leq 3} \left| \omega_{k}^{1/4} - e^{2j\pi i/4} \frac{x-y\sqrt{\core\left( N_{\alpha} \right)}}{x+y\sqrt{\core\left( N_{\alpha} \right)}} \right|.
\]

From \eqref{eq:27}, our expression in the statement of this lemma for
$x+y\sqrt{\core\left( N_{\alpha} \right)}$, the expressions for $y_{k}$ and
$y_{\ell}$ in Proposition~\ref{prop:quad-rep}, and \eqref{eq:yk-defn}
(which implies that $\left| x_{k}-N_{\varepsilon^{k}}\sqrt{N_{\alpha}} \right|^{2}
=x_{k}^{2}-N_{\alpha}=dy_{k}^{2}$),
we have
\begin{align}
\label{eq:omega-UB}
& \left| \omega_{k} - \left( \frac{x-y\sqrt{\core\left( N_{\alpha} \right)}}{x+y\sqrt{\core\left( N_{\alpha} \right)}} \right)^{4} \right| \\
= & \frac{2\sqrt{\left| N_{\alpha} \right|}f_{k}^{2}f_{\ell}^{2}y_{k}^{2}}
{b^{4}\left| x_{k}-N_{\varepsilon^{k}}\sqrt{N_{\alpha}} \right| \left| r_{k} \mp s_{k}\sqrt{\core\left( N_{\alpha} \right)} \right|^{4}
\left| r_{\ell} \mp s_{\ell}\sqrt{\core\left( N_{\alpha} \right)} \right|^{4}} \nonumber \\
= & \frac{2\sqrt{\left| N_{\alpha} \right|}}{\sqrt{d} \, y_{\ell}}
\leq \frac{1}{2}. \nonumber
\end{align}
since $y_{\ell} \geq 4\sqrt{\left| N_{\alpha} \right|/d}$.

By Lemma~\ref{lem:omega-bnd}(b) with $c_{0}=1/2+0.0001$,
\begin{equation}
\label{eq:omega14-UB}
\left| \omega_{k}^{1/4} - \zeta_{4} \frac{x-y\sqrt{\core \left( N_{\alpha} \right)}}
{x+y\sqrt{\core\left( N_{\alpha} \right)}} \right|
<0.1263.
\end{equation}

From \eqref{eq:3.7}, along with $y_{k}, y_{\ell} \geq 4\sqrt{\left| N_{\alpha} \right|/d}$, we have
\begin{align*}
\left| \frac{\zeta_{4}^{(k)}}{\zeta_{4}^{(\ell)}}
\frac{x-y\sqrt{\core \left( N_{\alpha} \right)}}{x+y\sqrt{\core \left( N_{\alpha} \right)}} - 1 \right|
&=\left| \zeta_{4}^{(k)} \frac{r_{k}+s_{k}\sqrt{\core \left( N_{\alpha} \right)}}{r_{k}-s_{k}\sqrt{\core \left( N_{\alpha} \right)}}
- \zeta_{4}^{(\ell)} \frac{r_{\ell}+s_{\ell}\sqrt{\core \left( N_{\alpha} \right)}}{r_{\ell}-s_{\ell}\sqrt{\core \left( N_{\alpha} \right)}} \right| \\
&< 0.5051b^{2} \frac{\sqrt{\left| N_{\alpha} \right|}}{\sqrt{a^{2}+\left| N_{\alpha} \right|}}
\left( \frac{1}{y_{k}} + \frac{1}{y_{\ell}} \right)
<0.253.
\end{align*}

From part~(a), we have $\left| \varphi_{k} \right|<0.6$, so $\left| \omega_{k}^{1/4}-1 \right|<0.15$ and
\[
\left| \omega_{k}^{1/4} - \frac{\zeta_{4}^{(k)}}{\zeta_{4}^{(\ell)}}
\frac{x-y\sqrt{\core \left( N_{\alpha} \right)}}{x+y\sqrt{\core \left( N_{\alpha} \right)}} \right|
\leq \left| \omega_{k}^{1/4}-1 \right|
+ \left| \frac{\zeta_{4}^{(k)}}{\zeta_{4}^{(\ell)}}
\frac{x-y\sqrt{\core \left( N_{\alpha} \right)}}{x+y\sqrt{\core \left( N_{\alpha} \right)}} - 1 \right|
<0.403.
\]

Recalling \eqref{eq:omega14-UB}, it follows that
\begin{align*}
& \left| \omega_{k}^{1/4} - \zeta_{4}'\frac{\left( x-y\sqrt{\core \left( N_{\alpha} \right)} \right)}
{\left( x+y\sqrt{\core \left( N_{\alpha} \right)} \right)} \right| \\
= & \left| \omega_{k}^{1/4}
          -\zeta_{4}\frac{x-y\sqrt{\core \left( N_{\alpha} \right)}}{x+y\sqrt{\core \left( N_{\alpha} \right)}}
          +\left( \zeta_{4}-\zeta_{4}' \right)\frac{x-y\sqrt{\core \left( N_{\alpha} \right)}}{x+y\sqrt{\core \left( N_{\alpha} \right)}}
  \right| \\
& \geq \left| \zeta_{4}-\zeta_{4}' \right| - \left| \omega_{k}^{1/4}
   -\zeta_{4}\frac{x-y\sqrt{\core \left( N_{\alpha} \right)}}{x+y\sqrt{\core \left( N_{\alpha} \right)}} \right|
> \sqrt{2}-0.127,
\end{align*}
for any $4$-th root of unity, $\zeta_{4}'$, with $\zeta_{4}' \neq \zeta_{4}$.
Since this exceeds $0.403$, it follows that
$\zeta_{4}=\zeta_{4}^{(k)}/\zeta_{4}^{(\ell)}=\pm 1$, the last equality holding
by our assumption in the statement of this lemma.
\end{proof}

In Lemma~\ref{lem:gcd}, we establish a gcd result for elements of a generalisation
of our sequences. This will help us prove Lemmas~\ref{lem:g-gcd} and \ref{lem:gNd4}.
Lemma~\ref{lem:gNd4} is used in the proof of Proposition~\ref{prop:4.1}, as
well as in the proof of Theorem~\ref{thm:1.3-seq-new}. Lemma~\ref{lem:gNd4}
will be particularly important for showing that the possible exceptions to
Conjecture~\ref{conj:3-seq} in the statement of Theorem~\ref{thm:1.3-seq-new}
can only possibly occur for $u=1$ or $u=2$.

\begin{lemma}
\label{lem:gcd}
Let $\left( x_{k} \right)_{k=0}^{\infty}$ and $\left( y_{k} \right)_{k=0}^{\infty}$
be sequences defined by $x_{k}+y_{k}\sqrt{d} = \left( a+b\sqrt{d} \right) \varepsilon^{k}$,
where $a,b,d$ are positive integers, $d$ is not a square and $\varepsilon
=\left( t+u\sqrt{d} \right)/2 \neq \pm 1$, with $t$ and $u$ integers, is a unit
in $\cO_{\bbQ \left( \sqrt{d} \right)}$.

If $x_{k}$ and $y_{k}$ are both integers and $\gcd(a,b)$ is odd, then
$\gcd \left( x_{k}, y_{k} \right)/\gcd(a,b)=1$ or $2$.
\end{lemma}

The additional hypothesis that $\gcd(a,b)$ is odd was omitted from the published
version of this lemma. If $\gcd(a,b)$ is even, then
$\gcd \left( x_{k}, y_{k} \right)/\gcd(a,b)=1/2$ also occurs.

\begin{proof}
We may assume without loss of generality that $\gcd(a,b)=1$. Otherwise, consider
below the sequences $\left( x_{k}/\gcd(a,b) \right)_{k=0}^{k=\infty}$ and
$\left( y_{k}/\gcd(a,b) \right)_{k=0}^{k=\infty}$ instead of $\left( x_{k} \right)_{k=0}^{k=\infty}$
and $\left( y_{k} \right)_{k=0}^{k=\infty}$. Since $x_{0}=a$, $x_{1}= \left( at+bdu \right)/2$,
$x_{-1}= \pm \left( at-dbu \right)/2$,
$y_{0}=b$, $y_{1}=\left( au+bt \right)/2$, $y_{-1}=\pm \left( au-bt \right)/2$,
the recurrence relation satisfied by both sequences has integer coefficients and
$\gcd(a,b)$ is odd, we see that $x_{k}/\gcd(a,b)$ and $y_{k}/\gcd(a,b)$ are
integers if and only if $x_{k}$ and $y_{k}$ are integers.

We start with some relationships that we will require for the proof itself.

Throughout the proof, we will use the fact that both the $x_{k}$'s and $y_{k}$'s
satisfy the recurrence relation
\begin{equation}
\label{eq:gcd-recur}
u_{k} = \Tr \left( \varepsilon \right) u_{k-1}-N_{\varepsilon} u_{k-2}
= tu_{k-1} \pm u_{k-2}.
\end{equation}

Next, we show that
\begin{equation}
\label{eq:gcd-rel}
ux_{k}-ty_{k}=-2N_{\varepsilon}y_{k-1}
\quad \text{ and } \quad
tx_{k}-duy_{k}=2N_{\varepsilon}x_{k-1}
\end{equation}
for $k \geq 1$.

For $k=1$, we have $x_{1}=(at+bdu)/2$ and $y_{1}=(au+bt)/2$, so we find
$ux_{1}-ty_{1}=-2bN_{\varepsilon}=-2y_{0}N_{\varepsilon}$. Similarly,
$tx_{1}-duy_{1}=2N_{\varepsilon}x_{0}$ holds.

Similarly, we can prove the relationship holds for $k=2$.

For $k \geq 3$, we use induction and the recurrence relations in \eqref{eq:gcd-recur}:
\begin{align*}
ux_{k}-ty_{k}
&= u \Tr \left( \varepsilon \right) x_{k-1}-uN_{\varepsilon} x_{k-2}
- t\Tr \left( \varepsilon \right) y_{k-1}+tN_{\varepsilon} y_{k-2} \\
&= \Tr \left( \varepsilon \right) (-2)N_{\varepsilon} y_{k-2}
- N_{\varepsilon} (-2)N_{\varepsilon} y_{k-3} \\
&= (-2) N_{\varepsilon} y_{k-1}.
\end{align*}

The proof of $tx_{k}-duy_{k}=2N_{\varepsilon}x_{k-1}$ is identical, so we omit
it here.

We will also need
\begin{align}
\label{eq:gcd-rel2}
\frac{tu}{N_{\varepsilon}} x_{k}
-\frac{t^{2}+du^{2}}{2N_{\varepsilon}} y_{k}
&=-2N_{\varepsilon}y_{k-2}
\quad \text{ and } \\
\frac{t^{2}+du^{2}}{2N_{\varepsilon}} x_{k}
-\frac{dtu}{N_{\varepsilon}} y_{k}
&=2N_{\varepsilon}x_{k-2}, \nonumber
\end{align}
for $k \geq 2$ and
\begin{align}
\label{eq:gcd-rel3}
\left( \frac{u\left( 3t^{2}+du^{2} \right)}{4N_{\varepsilon}^{2}} \right) x_{k}
-\left( \frac{t\left( t^{2}+3du^{2} \right)}{4N_{\varepsilon}^{2}} \right) y_{k}
&=-2N_{\varepsilon}y_{k-3}
\quad \text{ and } \\
\left( \frac{t\left( t^{2}+3du^{2} \right)}{4N_{\varepsilon}^{2}} \right) x_{k}
-\left( \frac{du\left( 3t^{2}+du^{2} \right)}{4N_{\varepsilon}^{2}} \right) y_{k}
&=2N_{\varepsilon}x_{k-3}, \nonumber
\end{align}
for $k \geq 3$, which both follow from \eqref{eq:gcd-rel}.

Writing $3t^{2}+du^{2}=-4N_{\varepsilon}+4t^{2}$ and
$t^{2}+3du^{2}=4N_{\varepsilon}+4du^{2}$, we see that the coefficients of $x_{k}$
and $y_{k}$ in the relationships in \eqref{eq:gcd-rel3}
are integers.

We divide the remainder of the proof into cases according to the parity of $t$ and $u$.

\vspace*{3.0mm}

(i) Suppose that $t$ and $u$ are both even. Then $x_{k}$ and $y_{k}$ are both integers.
Furthermore in this case, the coefficients of $x_{k}$ and $y_{k}$
in \eqref{eq:gcd-rel} can all be divided by $2$ to eliminate the factors of $2$
on the right-hand sides. So in this case, we find that
\[
\gcd \left( x_{k}, y_{k} \right) | \gcd \left( x_{k-1}, y_{k-1} \right) |
\cdots | \gcd \left( x_{0}, y_{0} \right)=1.
\]

\vspace*{3.0mm}

(ii) Next suppose that $t$ and $u$ are both odd. In this case, we must have $d \equiv 5 \pmod{8}$,
since $\varepsilon$ is a unit and hence $t^{2}-du^{2} \equiv 4 \pmod{8}$.

(ii-a) Suppose that $a$ and $b$ have opposite
parities. In this case, $y_{1}=(au+bt)/2$ cannot be an integer. Similarly,
$2x_{1}$ is odd.

Since $\left( \left( t + u\sqrt{d} \right)/2 \right)^{2}
=\left( t^{2}+u^{2}d \right)/4+(tu/2)\sqrt{d}= \left( t_{2}+u_{2}\sqrt{d} \right)/2$,
where $t^{2}+u^{2}d=\pm 4 +2du^{2} \equiv 2 \pmod{4}$, it follows that $t_{2}$
and $u_{2}$ are both odd. Hence $x_{2}$ and $y_{2}$ are not integers. However,
\[
\left( \frac{t + u\sqrt{d}}{2} \right)^{3}
= \frac{t^{3}+3tu^{2}d}{8} + \frac{3t^{2}u+u^{3}d}{8} \sqrt{d}
= \frac{t_{3}+u_{3}\sqrt{d}}{2},
\]
where $t^{3}+3tu^{2}d=t \left( t^{2}+3u^{2}d \right)
=t \left( \pm 4 +4u^{2}d \right) \equiv 0 \pmod{8}$ and
$3t^{2}u+u^{3}d=-u \left( -4t^{2} \pm 4 \right) \equiv 0 \pmod{8}$.
So both $t_{3}$ and $u_{3}$ are even.
Hence $x_{3}$ and $y_{3}$ are both integers.

By an inductive argument using the recurrence relations for the $x_{k}$'s and
$y_{k}$'s in \eqref{eq:gcd-recur}, one finds that $x_{k}$ and $y_{k}$ are both integers
if and only if $3|k$, when $t$ and $u$ are both odd and $a$ and $b$ have opposite
parities. So we consider only $\gcd \left( x_{3k}, y_{3k} \right)$.

Since $t$ and $du$ are both odd in this case, we have $3t^{2}+du^{2}=-4N_{\varepsilon}+4t^{2} \equiv 0 \pmod{8}$
and $t^{2}+3du^{2}=4N_{\varepsilon}+4du^{2} \equiv 0 \pmod{8}$, so the coefficients
of $x_{k}$ and $y_{k}$ in \eqref{eq:gcd-rel3}
can all be divided by $2$ to eliminate the factors of $2$ on the right-hand
sides there. Thus,
$\gcd \left( x_{k}, y_{k} \right)|\gcd \left( x_{k-3}, y_{k-3} \right)|
\cdots | \gcd \left( x_{0}, y_{0} \right)=1$.

\vspace*{1.0mm}

(ii-b) If $a$ and $b$ have the same parity, then $y_{1}=(au+bt)/2 \in \bbZ$.
Similarly, $x_{1} \in \bbZ$. Using the recurrence relations in \eqref{eq:gcd-recur},
the $x_{k}$'s and $y_{k}$'s are integers.

As in case~(ii-a), we can remove the factor of $2$ on the right-hand side of
\eqref{eq:gcd-rel3}. So we find that if $k \equiv k_{1} \pmod{3}$ with
$0 \leq k_{1} \leq 2$, then
$\gcd \left( x_{k}, y_{k} \right) | \gcd \left( x_{k_{1}}, y_{k_{1}} \right)$.

So if $3|k$, then 
$\gcd \left( x_{k}, y_{k} \right) | \gcd \left( x_{0}, y_{0} \right)$.

If $k \equiv 1 \pmod{3}$, then
$\gcd \left( x_{k}, y_{k} \right) | \gcd \left( x_{1}, y_{1} \right)$.
From \eqref{eq:gcd-rel},
\[
\gcd \left( x_{1}, y_{1} \right)| 2\gcd \left( x_{0}, y_{0} \right).
\]

If $k \equiv 2 \pmod{3}$, then
$\gcd \left( x_{k}, y_{k} \right) | \gcd \left( x_{2}, y_{2} \right)$.
From \eqref{eq:gcd-rel},
\[
\gcd \left( x_{2}, y_{2} \right) | \left( 2\gcd \left( x_{1}, y_{1} \right) \right)
| \left( 4\gcd \left( x_{0}, y_{0} \right) \right).
\]

We now consider the parity of $x_{k}$ and $y_{k}$. Since $\gcd(a,b)=1$, at least
one of $x_{0}$ or $y_{0}$ is odd.

Suppose first that $y_{0}$ is odd.
Since $t$ is odd, we find that $y_{2}=ty_{1} \pm y_{0}$ is even
if $y_{1}$ is odd and $y_{2}$ is odd otherwise.

In the first case (when $y_{1}$ is odd), we have
$\gcd \left( x_{1}, y_{1} \right)| \gcd \left( x_{0}, y_{0} \right)$
and hence
\[
\gcd \left( x_{2}, y_{2} \right) | \left( 2\gcd \left( x_{1}, y_{1} \right) \right)
| \left( 2 \gcd \left( x_{0}, y_{0} \right) \right).
\]

In the second case (when $y_{1}$ is even), we have
$\gcd \left( x_{2}, y_{2} \right)| \gcd \left( x_{1}, y_{1} \right)$
and hence
\[
\gcd \left( x_{2}, y_{2} \right) | \gcd \left( x_{1}, y_{1} \right)
| \left( 2\gcd \left( x_{0}, y_{0} \right) \right).
\]

Now we suppose that $x_{0}$ is odd.
Since $t$ is odd, we find that $x_{2}=tx_{1}+x_{0}$ is even
if $x_{1}$ is odd and $x_{2}$ is odd otherwise. As above, we find that
\[
\gcd \left( x_{2}, y_{2} \right)
| \left( 2 \gcd \left( x_{0}, y_{0} \right) \right).
\]

\vspace*{3.0mm}

(iii) Now suppose that $t$ is even and $u$ is odd. In this case, $4|d$.

Since $y_{1}=(au+bt)/2$, we see that $y_{1} \in \bbZ$ if and only if $a$ is even.
Hence all the $y_{k}$'s are integers if $a$ is even.
Since $t$ is even, using the recurrence relation for the $y_{k}$'s in
\eqref{eq:gcd-recur}, we see that for $a$ odd,
$y_{k} \in \bbZ$ when $k$ is even and $y_{k} \not\in \bbZ$ otherwise.

Since $x_{1}=(at+bdu)/2$ and $t$ and $d$ are both even, we see that $x_{1} \in \bbZ$.
Hence, by the recurrence relation for the $x_{k}$'s in \eqref{eq:gcd-recur},
$x_{k} \in \bbZ$ always holds.

Here we use \eqref{eq:gcd-rel2}, observing that the coefficients of $x_{k}$ and
$y_{k}$ on the left-hand sides are all even. So we have relationships with integer
coefficients and without the factor of $2$ on the right-hand sides. Hence
\[
\gcd \left( x_{k}, y_{k} \right)
\text{ divides } \left\{
\begin{array}{ll}
\gcd \left( x_{0}, y_{0} \right)=1 & \text{if $k$ is even}, \\
\gcd \left( x_{1}, y_{1} \right) | \left( 2 \gcd \left( x_{0}, y_{0} \right) \right)=2 & \text{if $k$ is odd}.
\end{array}
\right.
\]

\vspace*{3.0mm}

(iv) We cannot have $t$ odd and $u$ even because $t^{2}-du^{2} \equiv 0 \pmod{4}$.
\end{proof}

We will use Lemma~\ref{lem:g-gcd} in the proof of Lemma~\ref{lem:gNd4}.

\begin{lemma}
\label{lem:g-gcd}
Let the sequences $\left( x_{k} \right)_{k=-\infty}^{\infty}$ and
$\left( y_{k} \right)_{k=-\infty}^{\infty}$ be as defined in \eqref{eq:yk-defn},
with the notation and assumptions there. Suppose that $k \neq 0$, $\gcd(a,b)=1$,
$N_{\alpha}<0$ and that $x_{k}$ and $y_{k}$ are both integers.
Using the notation of Subsection~$\ref{subsect:const}$ with $t'=\core \left( N_{\alpha} \right)$,
$u_{1}=2x_{k}$ and $u_{2}=\pm 2\sqrt{N_{\alpha}/\core\left( N_{\alpha} \right)}$,
then
\[
g^{2}/\gcd \left( a^{2}, db^{4} \right)=2^{m}, \text{ where $m \geq 0$.}
\]
\end{lemma}

\begin{remark-nonum}
With more work, one can show that $m=0,1,2,4$.
\end{remark-nonum}

\begin{proof}
We break the proof into several cases.

(1) Suppose that $p>2$ is prime.

Put $v_{p} \left( g^{2} \right)=v_{p} \left( g_{1}^{2}g_{2}/g_{3} \right)=\ell$.
Since $g_{3}$ is a power of $2$, it follows that $v_{p} \left( g_{3} \right)=0$
and hence $\ell \geq 0$. We will prove that
$v_{p} \left( \gcd \left( a^{2}, d \right) \right)=\ell$.
We also put
$v_{p} \left( g_{1}^{2} \right)=\ell_{1}$ and $v_{p} \left( g_{2} \right)=\ell_{2}$.
Note that $\ell_{2}=0$ or $1$, since $g_{2}|\core \left( N_{\alpha} \right)$
(from the definition of $g_{2}$, we have $g_{2}|t'$ and here $t'=\core \left( N_{\alpha} \right)$).
Since $v_{p} \left( g_{3} \right)=0$, we have $\ell=\ell_{1}+\ell_{2}$.

Since $g_{1}|u_{2}$, by the definition of $g_{1}$, and from the expression for
$u_{2}$ here, we have $p^{\ell_{1}}| \left( N_{\alpha}/\core \left( N_{\alpha} \right) \right)$.
Similarly, from $g_{2}|\core \left( N_{\alpha} \right)$, it follows that
$p^{\ell_{2}}|\core \left( N_{\alpha} \right)$. So $p^{\ell}|N_{\alpha}$.

In the same way, we find that $p^{\ell_{1}+2\ell_{2}}|x_{k}^{2}$. So
$p^{\ell}|x_{k}^{2}$ too.

By Lemma~\ref{lem:gcd} and $\gcd(a,b)=1$, since $p$ is odd,
$p \nmid \gcd \left( x_{k}, y_{k} \right)$. Hence, from
$x_{k}^{2}-dy_{k}^{2}=N_{\alpha}$, we have $p^{\ell}|d$.
Combining this
with $p^{\ell}|N_{\alpha}$ and $N_{\alpha}=a^{2}-db^{4}$, we find that
$p^{\ell}|a^{2}$. So $p^{\ell}| \gcd \left( a^{2}, d \right)$.

We now show that $p^{\ell+1} \nmid \gcd \left( a^{2}, d \right)$, so that
the denominator of $g^{2}/\gcd \left( a^{2}, d \right)$ has no odd prime
factors.

Suppose that $p^{\ell+1}| \gcd \left( a^{2}, d \right)$. Then
$p^{\ell+1}|N_{\alpha}$. Also
$p^{\ell+1}|x_{k}^{2}$ follows from $x_{k}^{2}-dy_{k}^{2}=N_{\alpha}$.
By the definition of $\ell_{1}$ and $g_{1}$, this means that
$p^{\ell_{1}} \parallel N_{\alpha}/\core \left( N_{\alpha} \right)$.

If $p \nmid \core \left( N_{\alpha} \right)$, then $p^{\ell_{1}} \parallel N_{\alpha}$,
but $\ell_{1}<\ell+1=\ell_{1}+\ell_{2}+1$, contradicting 
$p^{\ell+1} | N_{\alpha}$.

If $p \mid \core \left( N_{\alpha} \right)$, then $p^{\ell_{1}+1} \parallel N_{\alpha}$.
We must have $\ell_{2}=0$, since otherwise $\ell_{1}+1<\ell_{1}$ and we get a
contradiction to $p^{\ell+1} | N_{\alpha}$. In this case,
$p \nmid \left( x_{k}^{2}/g_{1}^{2} \right)$. So $p^{\ell_{1}} \parallel x_{k}^{2}$.
But this contradicts $p^{\ell+1}|x_{k}^{2}$.
Hence $p^{\ell+1} \nmid \gcd \left( a^{2}, d \right)$.
Therefore, for all primes $p>2$, $v_{p} \left( g^{2} / \gcd \left( a^{2}, d \right) \right)=0$.

\vspace*{1.0mm}

(2) We now consider $p=2$.

Put
$v_{2} \left( \gcd \left( a^{2}, db^{4} \right) \right)=\ell$.
Since $\gcd(a,b)=1$, we also have $v_{2} \left( \gcd \left( a^{2}, d \right) \right)=\ell$.
So $2^{\ell}| N_{\alpha}$ too.
From $x_{k}^{2}-dy_{k}^{2}=N_{\alpha}$, we find
that $2^{\ell}|x_{k}^{2}$. Hence $2^{\ell}| \gcd \left( x_{k}^{2}, N_{\alpha} \right)$.

Furthermore, if $v_{2} \left( a^{2} \right)=v_{2} \left( db^{4} \right)=\ell$,
then $2^{\ell+1}| N_{\alpha}$.

\vspace*{1.0mm}

(2.1) Suppose that $2 \nmid \core \left( N_{\alpha} \right)$.

We have
$v_{2} \left( g_{1}^{2} \right) = v_{2} \left( 4 \gcd \left( x_{k}^{2}, N_{\alpha} \right) \right)
\geq \ell+2$. Also $v_{2} \left( g_{2} \right)=1$. Hence $v_{2} \left( g^{2} \right) = v_{2} \left( g_{1}^{2}/g_{3} \right)
\geq \ell$, recalling that $g_{3}=1,2$ or $4$.

\vspace*{1.0mm}

(2.2) Suppose that $2 | \core \left( N_{\alpha} \right)$.

Here
$v_{2} \left( g_{1}^{2} \right) = v_{2} \left( \gcd \left( 4x_{k}^{2}, 4N_{\alpha}/ \core \left( N_{\alpha} \right) \right) \right)
= v_{2} \left( 4\gcd \left( x_{k}^{2}, N_{\alpha}/2 \right) \right)
\geq \ell+1$. We break this case into subcases.

(2.2.1) If $v_{2} \left( a^{2} \right)=v_{2} \left( db^{4} \right)=\ell$,
then $v_{2} \left( N_{\alpha} \right) \geq \ell+1$.
So $v_{2} \left( g_{1}^{2} \right) = v_{2} \left( 4 \gcd \left( x_{k}^{2}, N_{\alpha}/2 \right) \right)
\geq \ell+2$. Here $v_{2} \left( g_{2} \right) \geq 1$, so as in case~(2.1), we have
$v_{2} \left( g^{2} \right) \geq \ell$, as required.

(2.2.2) If $v_{2} \left( a^{2} \right) \neq v_{2} \left( db^{4} \right)$, then
$v_{2} \left( N_{\alpha} \right)=\ell$. So $v_{2} \left( g_{1}^{2} \right)
=v_{2} \left( 4N_{\alpha}/\core \left( N_{\alpha} \right) \right)= \ell+1$,
since $2|\core \left( N_{\alpha} \right)$.
From $v_{2} \left( g_{1}^{2} \right)=\ell+1$, it also follows that
$\ell$ must be odd.

From $x_{k}^{2}-dy_{k}^{2}=N_{\alpha}$ and $v_{2} \left( N_{\alpha} \right)=\ell$,
we have either $v_{2} \left( x_{k}^{2} \right)=\ell$ and
$v_{2} \left( dy_{k}^{2} \right) \geq \ell+1$; or else
$v_{2} \left( x_{k}^{2} \right) \geq \ell+1$ and $v_{2} \left( dy_{k}^{2} \right)=\ell$.
The first case is not possible because $\ell$ is odd.
In the latter case, we must have $v_{2}(d)=\ell$ and $2 \nmid y_{k}$. Also,
\[
v_{2} \left( u_{1}/g_{1} \right)=(1/2) v_{2} \left( 4x_{k}^{2}/g_{1}^{2} \right)
\geq (1/2) \left( \ell+3-(\ell+1) \right)=1,
\]
so $v_{2} \left( g_{2} \right)=1$.
Hence
$v_{2} \left( g^{2} \right)=v_{2} \left( g_{1}^{2}g_{2}/g_{3} \right)
=v_{2} \left( g_{1}^{2}g_{2}/4 \right)
=\ell$.
\end{proof}

\begin{lemma}
\label{lem:gNd4}
Suppose that $a$, $b$, $d$, $g$, $t'$, $u_{1}$, $u_{2}$, $x_{k}$ and
$y_{k}$ are as in Lemma~$\ref{lem:g-gcd}$, $d'=u_{2}^{2}t'/g^{2}$, as defined in
$\eqref{eq:d-defn}$, $\cN_{d',4}$ is as defined in 
$\eqref{eq:ndn-defn}$
and put $b'=N_{\alpha}/\gcd \left( a^{2}, db^{4} \right)$. Then
\begin{equation}
\label{eq:gn-values}
|g|\cN_{d',4}
=\sqrt{\gcd \left( a^{2}, db^{4} \right)} \, 2^{1+\min \left( 2, v_{2}(b')/2 \right)}
\geq 2^{1+\min \left( 2, v_{2}\left( N_{\alpha} \right)/2 \right)}.
\end{equation}
\end{lemma}

\begin{proof}
Since $\gcd \left( a^{2}, db^{4} \right) \geq 2^{v_{2} \left( \gcd \left( a^{2}, db^{4} \right) \right)}$,
the inequality in equation~\eqref{eq:gn-values} holds.

To establish the equality in equation~\eqref{eq:gn-values}, we use Lemma~\ref{lem:g-gcd}.
From that lemma, we can write $g^{2}=2^{m}\gcd \left( a^{2}, db^{4} \right)$, where
$m \geq 0$. From \eqref{eq:ndn-defn} with $n=4$, we find that
$\cN_{d',4}^{2}=2^{\min \left( v_{2}(d'), 6 \right)}$.
So, squaring both sides of the equality in \eqref{eq:gn-values},
it is equivalent to
\[
2^{\min \left( 6+m, v_{2}(d')+m \right)}
=2^{\min \left( 6, v_{2}(4b') \right)}.
\]

From the definition of $d'$ and Lemma~\ref{lem:g-gcd}, we have
$d'=4N_{\alpha}/g^{2}=4N_{\alpha}/ \left( 2^{m} \gcd \left( a^{2}, db^{4} \right) \right)$
Applying this, along with the definition of $b'$, it follows that the equality
in \eqref{eq:gn-values} is equivalent to
\[
\min \left( 6+m, v_{2} \left( 4N_{\alpha}/\gcd \left( a^{2}, db^{4} \right) \right) \right)
=\min \left( 6, v_{2} \left( 4N_{\alpha}/\gcd \left( a^{2}, db^{4} \right) \right) \right).
\]

To show that this equality holds, we prove that if
$v_{2} \left( 4N_{\alpha}/\gcd \left( a^{2}, db^{4} \right) \right)>6$, then $m=0$
(i.e., $g^{2}=\gcd \left( a^{2}, db^{4} \right)$).
Putting $\ell=v_{2} \left( \gcd \left( a^{2}, db^{4} \right) \right)$, we must
have $v_{2} \left( a^{2} \right)=v_{2} \left( db^{4} \right)=v_{2} \left( d \right)=\ell$.
The second-last equality holds because $\gcd(a,b)=1$.
From Lemma~\ref{lem:gcd}, we have $\gcd \left( x_{k}, y_{k} \right)=1,2$. We
will show that $y_{k}$ is odd and hence $\gcd \left( x_{k}, y_{k} \right)=2$ is
not possible here.

We assume that $y_{k}$ is even and derive a contradiction. From
$x_{k}^{2}-dy_{k}^{2}=N_{\alpha}$, we have
$\left( x_{k}/2^{\ell/2+1} \right)^{2} - \left( d/2^{\ell} \right) \left( y_{k}/2 \right)^{2}
=N_{\alpha}/2^{\ell+2}$. Since $v_{2} \left( 4N_{\alpha}/\gcd \left( a^{2}, db^{4} \right) \right)>6$,
we have $v_{2} \left( N_{\alpha} \right) \geq \ell+5$.
It follows that $x_{k}/2^{\ell/2+1}$ and $y_{k}/2$ are both odd. Arguing mod $8$,
we find that $d \equiv 1 \pmod{8}$ and that $t$ and $u$ must both be even.

From \eqref{eq:yPM1}, we can write
$4y_{k}=\left( b^{2}\left( 2du^{2} \pm 4 \right) \pm 2atu \right)$ for some positive
integers $t$ and $u$ with $t^{2}-du^{2}=\pm 4$. Since $t$ and $u$ are both even,
it follows that $v_{2} \left( 4y_{k} \right)=v_{2} \left( 2du^{2} \pm 4 \right)=2$.
This contradicts the assumption that $y_{k}$ is even.

From $y_{k}$ is odd and
$\left( x_{k}/2^{\ell/2} \right)^{2} - \left( d/2^{\ell} \right) y_{k}^{2}
=N_{\alpha}/2^{\ell}$, we have
\[
v_{2} \left( g_{1}^{2} \right) = \min \left( v_{2} \left( 4x_{k}^{2} \right), v_{2} \left( 4N_{\alpha}/\core\left( N_{\alpha} \right) \right) \right)
= v_{2} \left( 4x_{k}^{2} \right) = \ell+2.
\]

This also implies that $v_{2} \left( 4x_{k}^{2}/g_{1}^{2} \right)=0$, so $v_{2} \left( g_{2} \right)=0$.

Furthermore, since $v_{2} \left( 4x_{k}^{2} \right)< v_{2} \left( 4N_{\alpha}/\core\left( N_{\alpha} \right) \right)$,
we also have $g_{3}=4$. Hence $g^{2}=\gcd \left( a^{2}, db^{4} \right)$, as desired
and the lemma follows.
\end{proof}

The quantities in Lemma~\ref{lem:eq-LB} are related to the quantities
$E$ and $Q$ that we will use in the proof of Proposition~\ref{prop:4.1}.

\begin{lemma}
\label{lem:eq-LB}
Suppose that $b=1$, $x_{k}$ and $y_{k}$ are both integers with $y_{k}>1$
and that $g$ and $\cN_{d',4}$ are as above.

If $d \geq 105$, then
\begin{equation}
\label{eq:E-LB}
\frac{0.1832 |g|\cN_{d',4} \sqrt{d} \, y_{k}}{\left| N_{\alpha} \right|} > 1.13
\end{equation}
and
\begin{equation}
\label{eq:Q-LB}
\frac{21.12\sqrt{d} \, y_{k}}{|g|\cN_{d',4}} > 217.
\end{equation}
\end{lemma}

\begin{remark-nonum}
Our choice of the lower bound for $d$ comes from an example of $E<1$ for $d=104$:
$a=9$, $d=104$, $N_{\alpha}=-23$, $k=-1$, $x_{k}=-61$ and $y_{k}=11$ where
$E=0.973\ldots$.
\end{remark-nonum}

\begin{proof}
We first obtain an analytic lower bound for $d$ such that \eqref{eq:E-LB}
and \eqref{eq:Q-LB} hold for larger $d$.

We have $|g|\cN_{d',4} \geq 2^{1+\min \left( v_{2}(N_{\alpha})/2, 2 \right)} \geq 2$
from \eqref{eq:gn-values}. Since $b=1$, from Lemma~\ref{lem:Y-LB}(c),
we find that $y_{k} \geq \left| N_{\alpha} \right|/4$ holds, so
\[
\frac{0.1832 |g|\cN_{d',4} \sqrt{d} \, y_{k}}{\left| N_{\alpha} \right|}
>0.091\sqrt{d}.
\]

For $d \geq 155$, we find that $0.091\sqrt{d}>1.13$ holds.

Similarly,
\[
\frac{21.12\sqrt{d} \, y_{k}}{|g|\cN_{d',4}}
\geq \frac{21.12\sqrt{d} \, \left( \left| N_{\alpha} \right|/4 \right)}{8\sqrt{\left| N_{\alpha} \right|}}
\geq 0.66\sqrt{d},
\]
where we use $|g|\cN_{d',4} \leq 8\sqrt{\left| N_{\alpha} \right|}$ from \eqref{eq:gn-values}
and Lemma~\ref{lem:Y-LB}(c) with $u \geq 1$ to establish the first inequality.
For $d \geq 109,000$, we find that $0.66\sqrt{d}>217$ holds.

We complete the proof computationally, checking all the remaining pairs, $(a,d)$.
Since $d$ is bounded and $N_{\alpha}<0$ (so $a^{2}<d$), there are only finitely
many such pairs.

For each pair, we check \eqref{eq:E-LB} and \eqref{eq:Q-LB} by computing their
left-hand sides for all $k \neq 0$ such that
$0.1832 \cdot 2^{1+\min \left( v_{2}(b)/2, 2 \right)}\sqrt{d}y_{k}/b<2$ and
$21.12y_{k}/2^{1+\min \left( v_{2}(b)/2, 2 \right)}<300$. No counterexamples to
\eqref{eq:E-LB} and \eqref{eq:Q-LB} were found.
A PARI/GP program took 112 seconds to run on a Windows 10 laptop
with an Intel i7-9750H processor and 16~GB of RAM.

The smallest value of the left-hand side of \eqref{eq:E-LB} found for $d \geq 105$
was $1.139\ldots$,
for $a=11$, $d=140$, $N_{\alpha}=-19$, $k=-1$, $x_{k}=-59$ and $y_{k}=5$.
The smallest value of the left-hand side of \eqref{eq:Q-LB} found for $d \geq 105$
was $217.3\ldots$,
for $a=10$, $d=140$, $N_{\alpha}=-40$, $k=-1$, $x_{k}=-130$ and $y_{k}=11$.

\end{proof}

\section{Proposition~\ref{prop:4.1} and its proof}
\label{sect:prop-11}

Our results in Subsection~\ref{subsect:results} follow from
the next proposition.

\begin{proposition}
\label{prop:4.1}
Let $\left( y_{k} \right)_{k=-\infty}^{\infty}$
be defined by \eqref{eq:yk-defn}.

If $b=1$ and $-N_{\alpha}$ is a square, then there is at most one integer
square among all distinct elements of $\left( y_{k} \right)_{k=-\infty}^{\infty}$
which satisfies
\[
y_{k}>\max \left( 1, \frac{76\left| N_{\alpha} \right|^{3/2}}{\sqrt{d} \left( |g|\cN_{d',4} \right)^{2}} \right).
\]
\end{proposition}

Conjecture~\ref{conj:3-seq} would immediately follow when $-N_{\alpha}$ is a
square, if we could replace the lower bound for $y_{k}$ in
Proposition~\ref{prop:4.1}
with $\max \left( 1, \left| N_{\alpha} \right|/4 \right)$.
This is due to Lemma~\ref{lem:Y-LB}.

\subsection{Prerequisites}
\label{subsect:preq}

In this subsection, we collect some inequalities what will be required in the
subsections that follow.

We will suppose that there are two distinct squares,
$y_{k}$ and $y_{\ell}$, in the sequence with $y_{\ell}>y_{k}>1$. So $y_{\ell}>y_{k} \geq 4
\geq \max \left( 4\sqrt{\left| N_{\alpha} \right|/d}, \left| N_{\alpha} \right|/d \right)$
(since $b=1$ and $N_{\alpha}=a^{2}-bd^{4}<0$, so $\left| N_{\alpha} \right|<d$),
as required in our lemmas above.

We shall initially assume that
\begin{equation}
\label{eq:d-LB}
d \geq 105.
\end{equation}

This is the condition in Lemma~\ref{lem:eq-LB}, which will
allow us to bound $E$ and $Q$ from below. This assumption shall be removed at
the end of the proof in Subsection~\ref{subsect:step-v}.

As in Lemma~\ref{lem:zeta}, we put
$\omega_{k}=\left( x_{k}+N_{\varepsilon^{k}}\sqrt{N_{\alpha}} \right)
/\left( x_{k}-N_{\varepsilon^{k}}\sqrt{N_{\alpha}} \right)$
and let $\zeta_{4}$ be the $4$-th root of unity such that
\[
\left| \omega_{k}^{1/4} - \zeta_{4} \frac{x-y\sqrt{\core \left( N_{\alpha} \right)}}{x+y\sqrt{\core \left( N_{\alpha} \right)}} \right|
= \min_{0 \leq j \leq 3} \left| \omega_{k}^{1/4} - e^{2j\pi i/4} \frac{x-y\sqrt{\core \left( N_{\alpha} \right)}}{x+y\sqrt{\core \left( N_{\alpha} \right)}} \right|,
\]
where $x+y\sqrt{\core \left( N_{\alpha} \right)}=\left(  r_{k}-s_{k}\sqrt{\core \left( N_{\alpha} \right)} \right)
\left( r_{\ell}+s_{\ell}\sqrt{\core \left( N_{\alpha} \right)} \right)$ with
$\left( r_{k}, s_{k} \right)$ and $\left( r_{\ell}, s_{\ell} \right)$ as in
Proposition~$\ref{prop:quad-rep}$, which are associated with $\left( x_{k}, y_{k} \right)$
and $\left( x_{\ell}, y_{\ell} \right)$, respectively.
From \eqref{eq:omega14-UB} in the proof of Lemma~\ref{lem:zeta}(b), we have
\[
\left| \omega_{k}^{1/4} - \zeta_{4} \frac{x-y\sqrt{\core\left( N_{\alpha} \right)}}{x+y\sqrt{\core\left( N_{\alpha} \right)}} \right|
<0.127.
\]

Thus we can apply Lemma~\ref{lem:omega-bnd}(a) with $c_{1}=0.127$ to find that
\begin{equation}
\label{eq:29}
\frac{2\sqrt{\left| N_{\alpha} \right|}}{\sqrt{d} \, y_{\ell}}
= \left| \omega_{k} - \left( \frac{x-y\sqrt{\core \left( N_{\alpha} \right)}}{x+y\sqrt{\core \left( N_{\alpha} \right)}} \right)^{4} \right|
> 3.959 \left| \omega_{k}^{1/4} - \zeta_{4} \frac{x-y\sqrt{\core \left( N_{\alpha} \right)}}{x+y\sqrt{\core \left( N_{\alpha} \right)}} \right|.
\end{equation}

The equality on the left-hand side is from the equalities in \eqref{eq:omega-UB}.

By our choice of $k$ and $\ell$, and by Lemma~\ref{lem:zeta}(b) when $-N_{\alpha}$
is not a square, $\zeta_{4}=\pm 1 \in \bbQ \left( \sqrt{\core\left( N_{\alpha} \right)} \right)$. This
is important for us here as $\zeta_{4} \left( x-y\sqrt{\core \left( N_{\alpha} \right)} \right)
/ \left( x+y\sqrt{\core\left( N_{\alpha} \right)} \right)$
must be in an imaginary quadratic field in order to apply Lemma~\ref{lem:2.1}
to obtain a lower bound for the rightmost quantity in \eqref{eq:29}.

We need to derive a lower bound for the far-right quantity in \eqref{eq:29}.
To do so, we shall use the lower bounds in Lemma~\ref{lem:2.1} with
a sequence of good approximations $p_{r}/q_{r}$ obtained from the hypergeometric
functions. So we collect here the required quantities.

Since $y_{k} \geq 4$ and $N_{\alpha}>-db^{4}=-d$ with $b=1$ here, we obtain
\begin{equation}
\label{eq:25}
x_{k}^{2}= dy_{k}^{2}+N_{\alpha}
> d \left( y_{k}^{2}-1 \right)
\geq 0.9375dy_{k}^{2}.
\end{equation}

So
\begin{equation}
\label{eq:26}
\sqrt{x_{k}^{2}-N_{\alpha}}= \sqrt{dy_{k}^{2}}
< 1.04 \left| x_{k} \right|.
\end{equation}

Using the notation of Subsection~\ref{subsect:const}, let $t'=\core \left( N_{\alpha} \right)$,
$u_{1}=2x_{k}$, $u_{2}=2\sqrt{N_{\alpha}/\core \left( N_{\alpha} \right)}$ and
$d'$ is as defined in \eqref{eq:d-defn}.

Substituting these quantities along with $\cD_{4}=e^{1.68}$ from Lemma~\ref{lem:denom-est}(a)
into the definition of $E$ in \eqref{eq:e-defn} and applying \eqref{eq:25}, we have
\begin{align}
\label{eq:E-LB1}
E &= \frac{|g|\cN_{d',4} \left|  \left| u_{1} \right|  + \sqrt{u_{1}^{2}-t'u_{2}^{2}} \right|}{\cD_{4}u_{2}^{2}|t'|}
= \frac{|g|\cN_{d',4}\left| \left| 2x_{k} \right| + 2\sqrt{x_{k}^{2}-N_{\alpha}} \right|}{4e^{1.68} \left| N_{\alpha} \right|} \\
&> \frac{|g|\cN_{d',4}\left| \left| x_{k} \right| + \sqrt{x_{k}^{2}-N_{\alpha}} \right|}{10.74 \left| N_{\alpha} \right|}
> \frac{|g|\cN_{d',4} \left( 1+\sqrt{0.9375} \right) \sqrt{d} \, y_{k}}{10.74 \left| N_{\alpha} \right|} \nonumber \\
&> \frac{0.1832|g|\cN_{d',4}\sqrt{d} \, y_{k}}{\left| N_{\alpha} \right|}. \nonumber
\end{align}

From \eqref{eq:Q-LB} in Lemma~\ref{lem:eq-LB}, we have $E>1.13>1$, as required
for its use with Lemma~\ref{lem:2.1}.

Similarly, using \eqref{eq:25} and \eqref{eq:Q-LB} in Lemma~\ref{lem:eq-LB},
we have
\begin{equation}
\label{eq:Q-LB1}
Q > \frac{2e^{1.68}\left( 1+\sqrt{0.9375} \right) \sqrt{d}y_{k}}{|g|\cN_{d',4}}
> 217,
\end{equation}
so the condition $Q>1$ in Lemma~\ref{lem:2.1} is satisfied.

From $N_{\alpha}<0$ and the equality in \eqref{eq:26}, we have
$x_{k}<\sqrt{x_{k}^{2}-N_{\alpha}}=\sqrt{d} \, y_{k}$,
so from Lemma~\ref{lem:denom-est}(a) we have
\begin{equation}
\label{eq:Q-UB2}
Q = \frac{e^{1.68} \left| \left| 2x_{k} \right| + 2\sqrt{x_{k}^{2}-N_{\alpha}}\right|}{|g|\cN_{d',4}}
<4\frac{e^{1.68} \sqrt{d} \, y_{k}}{|g|\cN_{d',4}}
<\frac{21.47\sqrt{d} \, y_{k}}{|g|\cN_{d',4}}.
\end{equation}

Recall from \eqref{eq:k-UB} that we take $k_{0}=0.89$.

Writing $\omega_{k}=e^{i\varphi_{k}}$, with $-\pi < \varphi_{k} \leq \pi$, from
Lemmas~\ref{lem:denom-est}(a) and \ref{lem:zeta}(a), we can take
\begin{equation}
\label{eq:ell-UB}
\ell_{0}=\cC_{4,2} \left| \varphi_{k} \right|<0.458\sqrt{\left| N_{\alpha} \right|}/ \left| x_{k} \right|.
\end{equation}

Also from Lemma~\ref{lem:zeta}(a), we have $\left| \varphi_{k} \right|<0.6$,
so the condition $\left| \omega_{k}-1 \right|<1$ in Lemma~\ref{lem:hypg} is satisfied.

Let $q=x+y\sqrt{\core \left( N_{\alpha} \right)}=\left( r_{k}-s_{k}\sqrt{\core \left( N_{\alpha} \right)} \right) \left( r_{\ell}+s_{\ell}\sqrt{\core \left( N_{\alpha} \right)} \right)$
and $p=x-y\sqrt{\core \left( N_{\alpha} \right)}$. Recall from \eqref{eq:abs}
with $b=1$ that
\begin{equation}
\label{eq:qAbs}
|q| = \sqrt{f_{k}f_{\ell}} \left( y_{k}y_{\ell} \right)^{1/4}.
\end{equation}

We are now ready to deduce the required contradiction from the assumption that
there are two sufficiently large squares when $-N_{\alpha}$ is a square. We will
break the proof into five parts.

With $r_{0}$ as in Lemma~\ref{lem:2.1}, we separate the case of $\zeta_{4}p/q \neq p_{r_{0}}/q_{r_{0}}$
for all $4$-th roots of unity, $\zeta_{4}$, from the case of
$\zeta_{4}p/q=p_{r_{0}}/q_{r_{0}}$ for some $4$-th
root of unity, $\zeta_{4}$. In the first case, Lemma~\ref{lem:2.1} provides a
suitable lower bound for the approximation. But in the second case, Lemma~\ref{lem:2.1}
is not strong enough, so we work directly with the approximations themselves.

\subsection{$r_{0}=1$ and $\zeta_{4}p/q \neq p_{1}/q_{1}$ for all $4$-th roots of unity, $\zeta_{4}$}
\label{subsect:step-i}

We start by determining an upper bound for $y_{\ell}$ for all $r_{0} \geq 1$ when
$\zeta_{4}p/q \neq p_{r_{0}}/q_{r_{0}}$, since we will also need such a result
in Subsection~\ref{subsect:step-iii}.

From \eqref{eq:29}, along with Lemma~\ref{lem:2.1}(b) and \eqref{eq:qAbs},
we have
\begin{equation}
\label{eq:32}
\frac{2\sqrt{\left| N_{\alpha} \right|}}{\sqrt{d} \, y_{\ell}}
> 3.959 \left| \omega_{k}^{1/4} - \zeta_{4} \frac{x-y\sqrt{\core \left( N_{\alpha} \right)}}{x+y\sqrt{\core \left( N_{\alpha} \right)}} \right|
> \frac{3.959(1-c)}{k_{0}Q^{r_{0}}\sqrt{f_{k}f_{\ell}} \left( y_{k}y_{\ell} \right)^{1/4}}.
\end{equation}

Applying \eqref{eq:k-UB} and \eqref{eq:Q-UB2} to \eqref{eq:32}, we obtain
\[
\frac{2\sqrt{\left| N_{\alpha} \right|}}{\sqrt{d} \, y_{\ell}}
> \frac{3.959(1-c)}
  {0.89 \left( 21.47\sqrt{d} \, y_{k} / \left( |g|\cN_{d',4} \right) \right)^{r_{0}}\sqrt{f_{k}f_{\ell}} \left( y_{k}y_{\ell} \right)^{1/4}}.
\]

After taking the fourth power of both sides and rearranging, we find that
\begin{equation}
\label{eq:y2UB-step1}
\left( \left| N_{\alpha} \right| f_{k}f_{\ell} \right)^{2}\left( \frac{0.45}{1-c} \right)^{4}
\left( \frac{21.47}{|g|\cN_{d',4}} \right)^{4r_{0}}
d^{2r_{0}-2} y_{k}^{4r_{0}+1}
>y_{\ell}^{3}.
\end{equation}

Specialising to the case when $r_{0}=1$ and using $|g|\cN_{d',4} \geq 2$ from
\eqref{eq:gn-values}, we have
\begin{equation}
\label{eq:y3UB}
y_{\ell}^{3} < 545(1-c)^{-4} \left( \left| N_{\alpha} \right|f_{k}f_{\ell} \right)^{2}y_{k}^{5}.
\end{equation}

We will now combine \eqref{eq:y3UB} with the gap principle in Lemma~\ref{lem:gap}(a)
to show that this case cannot occur.

Since $-N_{\alpha}$ is a square, we have $f_{k}=f_{\ell}=1$ from
Proposition~\ref{prop:quad-rep}(b), so combining the upper bound for $y_{\ell}^{3}$
in \eqref{eq:y3UB} with Lemma~\ref{lem:gap}(a) and cancelling the common factor
of $y_{k}^{5}$ on both sides, we find that
\[
57.32^{3} \left( \frac{d}{\left| N_{\alpha} \right|} \right)^{6}y_{k}^{4}
< 545(1-c)^{-4}N_{\alpha}^{2},
\]
which we can rewrite as
\[
y_{k} < \frac{0.24}{1-c}\frac{\left| N_{\alpha} \right|^{2}}{d^{3/2}}
< \frac{\left| 0.24N_{\alpha} \right|}{(1-c)d^{1/2}}.
\]
The last inequality holds because $\left| N_{\alpha} \right|=d-a^{2}<d$.

By \eqref{eq:d-LB}, we have $0.24/ \left( (1-c) d^{1/2} \right)<0.25$ if $c<0.9$.
We will see in Subsection~\ref{subsect:step-iii} that $c=0.75$ is a good choice.
So $y_{k}<\left| N_{\alpha} \right|/4$, which contradicts parts~(b) and (c) of
Lemma~\ref{lem:Y-LB}.
Hence this case cannot hold.

\subsection{$r_{0}=1$ and $\zeta_{4}p/q = p_{1}/q_{1}$ for some $4$-th root of unity, $\zeta_{4}$}
\label{subsect:step-ii}

As in Subsection~\ref{subsect:step-i}, we start by proving an upper bound for
$y_{\ell}$ that holds for all $r_{0} \geq 1$ with $\zeta_{4}p/q=p_{r_{0}}/q_{r_{0}}$
for some $4$-th root of unity, $\zeta_{4}$.

From the definitions of $p_{r_{0}}$, $q_{r_{0}}$ and $R_{r_{0}}$ in \eqref{eq:7},
along with parts~(a) and (e) of Lemma~\ref{lem:hypg}, we have
\begin{align}
\label{eq:omega14-LB}
  & \left| \omega_{k}^{1/4} - \zeta_{4} \frac{p}{q} \right|
    = \frac{1}{q_{r_{0}}} \left| q_{r_{0}}\omega_{k}^{1/4} - p_{r_{0}} \right| \nonumber \\
= & \left| \frac{N_{d',4,r_{0}}}{D_{4,r_{0}} Y_{1,4,r_{0}}(\omega_{k}) \left( \frac{u_{1}-u_{2}\sqrt{t'}}{2g} \right)^{r_{0}}} \right|
   \left| \frac{D_{4,r_{0}}}{N_{d',4,r_{0}}} R_{1,4,r_{0}}(\omega_{k}) \left( \frac{u_{1}-u_{2}\sqrt{t'}}{2g} \right)^{r_{0}} \right| \nonumber \\
= & \left| \frac{(\omega_{k}-1)^{2r_{0}+1}}{Y_{1,4,r_{0}}(\omega_{k})} \frac{(1/4) \cdots (r_{0}+1/4)}{(r_{0}+1) \cdots (2r_{0}+1)} 
		   {} _{2}F_{1} \left( r_{0}+3/4, r_{0}+1; 2r_{0}+2; 1-\omega_{k} \right) \right| \nonumber \\ 
\geq & \left| \frac{(\omega_{k}-1)^{2r_{0}+1}}{Y_{1,4,r_{0}}(\omega_{k})} \frac{(1/4) \cdots (r_{0}+1/4)}{(r_{0}+1) \cdots (2r_{0}+1)} \right|.
\end{align}

From Lemma~\ref{lem:hypg}(d), we have
\begin{equation}
\label{eq:step2-yUB}
\left| Y_{1,4,r_{0}} \left( \omega_{k} \right) \right|
< 1.072\frac{r_{0}!\Gamma(3/4)}{\Gamma(r_{0}+3/4)} \left| 1 + \sqrt{\omega_{k}} \right|^{2r_{0}}.
\end{equation}

We can write $\omega_{k}-1$ as $\left( 2N_{\alpha}+2x_{k}N_{\varepsilon^{k}}\sqrt{N_{\alpha}} \right)
/ \left( x_{k}^{2}-N_{\alpha} \right)$,
so
\[
\left| \omega_{k}-1 \right|=\sqrt{\frac{4\left| N_{\alpha} \right|}{x_{k}^{2}-N_{\alpha}}}
=2\sqrt{\frac{\left| N_{\alpha} \right|}{d}} \frac{1}{y_{k}}.
\]

Since $\left| \sqrt{\omega_{k}}+1 \right|<2$, it follows that
\[
\left| \sqrt{\omega_{k}}-1 \right|
=\frac{\left| \omega_{k}-1 \right|}{\left| \sqrt{\omega_{k}}+1 \right|}
>\sqrt{\frac{\left| N_{\alpha} \right|}{d}} \frac{1}{y_{k}}.
\]

From these two inequalities, we also obtain
\[
\left| \frac{(\omega_{k}-1)^{2r_{0}+1}}{\left( 1+\sqrt{\omega_{k}} \right)^{2r_{0}}} \right|
=\left| \left( \omega_{k}-1 \right) \left( \sqrt{\omega_{k}}-1 \right)^{2r_{0}} \right|
> 2 \left( \frac{\left| N_{\alpha} \right|}{d} \right)^{r_{0}+1/2}
\left( \frac{1}{y_{k}} \right)^{2r_{0}+1}.
\]

Applying this inequality together with the upper bound for $\left| Y_{1,4,r_{0}} \left( \omega_{k} \right) \right|$
in \eqref{eq:step2-yUB} and the inequalities \eqref{eq:gamma-bnds}
in Lemma~\ref{lem:denom-est}(b) to \eqref{eq:omega14-LB}, it follows that
\begin{align*}
\left| \omega_{k}^{1/4} - \zeta_{4} \frac{p}{q} \right|
&> \frac{0.2915}{4^{r_{0}} \cdot r_{0}^{1/2}}
\left( \frac{\left| N_{\alpha} \right|}{d} \right)^{r_{0}+1/2}
\left( \frac{1}{y_{k}} \right)^{2r_{0}+1}.
\end{align*}

Applying \eqref{eq:29} with this inequality, we obtain
\[
\frac{2\sqrt{\left| N_{\alpha} \right|}}{\sqrt{d} \, y_{\ell}}
> \frac{1.154}{4^{r_{0}} \cdot r_{0}^{1/2}}
\left( \frac{\left| N_{\alpha} \right|}{d} \right)^{r_{0}+1/2}
\left( \frac{1}{y_{k}} \right)^{2r_{0}+1},
\]
so
\begin{equation}
\label{eq:y2UB-step2}
1.734r_{0}^{1/2} \left(4\frac{d}{\left| N_{\alpha} \right|} \right)^{r_{0}}y_{k}^{2r_{0}+1}
>y_{\ell}.
\end{equation}

We now specialise to the case of $r_{0}=1$ and apply the gap principle in
Lemma~\ref{lem:gap}.

Our gap principle in Lemma~\ref{lem:gap}(a) with $b=1$, along with \eqref{eq:y2UB-step2},
implies that
\[
6.96 \frac{d}{\left| N_{\alpha} \right|} y_{k}^{3}>y_{\ell}
>57.32 \frac{d^{2}}{\left| N_{\alpha} \right|^{2}} y_{k}^{3}
>57.32 \frac{d}{\left| N_{\alpha} \right|} y_{k}^{3},
\]
since $-N_{\alpha}$ is a square and $d<\left| N_{\alpha} \right|$. But this
inequality is impossible. Hence this case cannot hold.

\subsection{$r_{0}>1$, $\zeta_{4}p/q \neq p_{r_{0}}/q_{r_{0}}$ for all $4$-th
roots of unity, $\zeta_{4}$}
\label{subsect:step-iii}

Here we establish a stronger gap principle for $y_{k}$ and $y_{\ell}$ than the
one in Lemma~\ref{lem:gap}. We then use this with the upper bound for $y_{\ell}$
in \eqref{eq:y2UB-step1} to obtain a contradiction.

We start by deriving a lower bound for $y_{\ell}$ that holds in both this step
and in the next step.

From the definition of $r_{0}$ in Lemma~\ref{lem:2.1},
along with \eqref{eq:Q-LB1} and $E>1$, we have
\begin{equation}
\label{eq:qLB-step3}
|q| \geq \frac{c(Q-1)}{\ell_{0}(Q-1/E)}E^{r_{0}-1}
> 0.995cE^{r_{0}-1}/\ell_{0}.
\end{equation}

Recall that $|q|=\sqrt{f_{k}f_{\ell}} \left( y_{k}y_{\ell} \right)^{1/4}$ by \eqref{eq:qAbs}.
Thus
\[
\left( y_{k}y_{\ell} \right)^{1/4}
> \frac{0.995cE^{r_{0}-1}}{\ell_{0}\sqrt{f_{k}f_{\ell}}}.
\]

Applying \eqref{eq:E-LB1} and \eqref{eq:ell-UB}, and then \eqref{eq:25}, to this
inequality, we obtain
\begin{align*}
\left( y_{k}y_{\ell} \right)^{1/4}
&> \frac{0.995c \left| x_{k} \right|}{0.458\sqrt{\left| N_{\alpha} \right| f_{k}f_{\ell}}} \left( \frac{0.1832|g|\cN_{d',4}\sqrt{d} \, y_{k}}{\left| N_{\alpha} \right|} \right)^{r_{0}-1} \\
&> \frac{2.103c\sqrt{d} \, y_{k}}{\sqrt{\left| N_{\alpha} \right|f_{k}f_{\ell}}}
  \left( \frac{0.1832|g|\cN_{d',4}\sqrt{d} \, y_{k}}{\left| N_{\alpha} \right|} \right)^{r_{0}-1}.
\end{align*}

Taking the fourth power of both sides and rearranging, we find that this
inequality implies
\begin{equation}
\label{eq:y2LB-step3}
y_{\ell} > \left( \frac{11.47}{|g|\cN_{d',4}} \sqrt{\frac{\left| N_{\alpha} \right|}{f_{k}f_{\ell}}} \right)^{4} c^{4}
           \left( \frac{0.1832 |g|\cN_{d',4}}{\left| N_{\alpha} \right|} \right)^{4r_{0}} d^{2r_{0}}
y_{k}^{4r_{0}-1}.
\end{equation}

With this lower bound for $y_{\ell}$, we now focus for the rest of this subsection
on when $\zeta_{4}p/q \neq p_{r_{0}}/q_{r_{0}}$ for all $4$-th
roots of unity, $\zeta_{4}$.

We now take the third power of both sides of this inequality and combine it
with the upper bound for $y_{\ell}^{3}$ in \eqref{eq:y2UB-step1}, obtaining
\begin{align}
\label{eq:y1UB-step3-a}
&
\left( \left| N_{\alpha} \right|f_{k}f_{\ell} \right)^{2}\left( \frac{0.45}{1-c} \right)^{4}
\left( \frac{21.47}{|g|\cN_{d',4}} \right)^{4r_{0}}
d^{2r_{0}-2} y_{k}^{4r_{0}+1} \\
& >
\left( \frac{11.47}{|g|\cN_{d',4}} \sqrt{\frac{\left| N_{\alpha} \right|}{f_{k}f_{\ell}}} \right)^{12} c^{12}
           \left( \frac{0.1832 |g|\cN_{d',4}}{\left| N_{\alpha} \right|} \right)^{12r_{0}}
		   d^{6r_{0}} y_{k}^{12r_{0}-3}. \nonumber
\end{align}

Using elementary calculus, $c^{12}(1-c)^{4}$ is monotonically
increasing for $0<c \leq 0.75$. So we put $c=0.75$
and find that for such $c$, $c^{12}(1-c)^{4}>0.000124$. Applying this to
\eqref{eq:y1UB-step3-a} and simplifying, we have
\begin{equation}
\label{eq:y1UB-step3-b}
\left( f_{k}f_{\ell} \right)^{8}
      \frac{0.00078\left| N_{\alpha} \right|^{2} \left( |g|\cN_{d',4} \right)^{4}}{d^{4}}
      \left( \frac{1.22 \cdot 10^{7}}
                  {\left( |g|\cN_{d',4} \right)^{8}} \right)^{2r_{0}-1}
> \left( \frac{y_{k}^{4}d^{2}}{\left| N_{\alpha} \right|^{6}} \right)^{2r_{0}-1}.
\end{equation}

Since $-N_{\alpha}$ is a square, we have $f_{k}=f_{\ell}=1$
from Proposition~\ref{prop:quad-rep}(b), so
\[
\frac{0.00078\left| N_{\alpha} \right|^{2} \left( |g|\cN_{d',4} \right)^{4}}{d^{4}}
\left( 1.22 \cdot 10^{7} \right)^{2r_{0}-1}
> \left( \frac{y_{k}^{4}d^{2}\left( |g|\cN_{d',4} \right)^{8}}{\left| N_{\alpha} \right|^{6}} \right)^{2r_{0}-1}
\]
and then
\[
\frac{0.00078\left| N_{\alpha} \right|^{2} \left( |g|\cN_{d',4} \right)^{4}}{d^{4}}
> \left( \frac{y_{k}^{4}d^{2}\left( |g|\cN_{d',4} \right)^{8}}
{59.2^{4} \left| N_{\alpha} \right|^{6}} \right)^{2r_{0}-1}.
\]

Using the equality in \eqref{eq:gn-values} of Lemma~\ref{lem:gNd4}, we can compute
that $\left| N_{\alpha} \right|^{2} \left( |g|\cN_{d',4} \right)^{4}/d^{4} \leq 2^{20}/17^{4}$
(the max value occurs when $d=17a^{2}$ so that $\left| N_{\alpha} \right|=16a^{2}$)
to obtain
\begin{equation}
\label{eq:ykUB-step3a}
0.01
> \left( \frac{y_{k}^{4}d^{2}\left( |g|\cN_{d',4} \right)^{8}}
{59.2^{4} \left| N_{\alpha} \right|^{6}} \right)^{2r_{0}-1}.
\end{equation}

But if
\begin{equation}
\label{eq:ykUB-step3b}
y_{k}>59.2\left| N_{\alpha} \right|^{3/2}/ \left( \sqrt{d} \left( |g|\cN_{d',4} \right)^{2} \right),
\end{equation}
then the right-hand side is greater than $1$, so this is not possible.

\subsection{$r_{0}>1$ and $\zeta_{4}p/q = p_{r_{0}}/q_{r_{0}}$ for some $4$-th root of unity, $\zeta_{4}$}
\label{subsect:step-iv}

We now combine our upper bound for $y_{\ell}$ in \eqref{eq:y2UB-step2} with our lower
bound for $y_{\ell}$ in \eqref{eq:y2LB-step3}. Thus
\[
1.734r_{0}^{1/2} \left(4\frac{d}{\left| N_{\alpha} \right|} \right)^{r_{0}}y_{k}^{2r_{0}+1}
> \left( \frac{11.47}{|g|\cN_{d',4}} \sqrt{\frac{\left| N_{\alpha} \right|}{f_{k}f_{\ell}}} \right)^{4} c^{4}
  \left( \frac{0.1832 |g|\cN_{d',4}}{\left| N_{\alpha} \right|} \right)^{4r_{0}} d^{2r_{0}}
  y_{k}^{4r_{0}-1}
\]
and so
\[
1.734r_{0}^{1/2}
> \left( \frac{11.47}{|g|\cN_{d',4}} \sqrt{\frac{\left| N_{\alpha} \right|}{f_{k}f_{\ell}}} \right)^{4} c^{4}
  \left( \frac{0.1832^{4} |g|^{4}\cN_{d',4}^{4}}{4\left| N_{\alpha} \right|^{3}} \right)^{r_{0}} d^{r_{0}}
  y_{k}^{2r_{0}-2}.
\]

We can show that $0.1832^{4r_{0}}/r_{0}^{1/2}>0.175^{4r_{0}}$, with the minimum
being attained at $r_{0}=3$. Applying this, along with $c=0.75$ and collecting
the terms taken to the power $r_{0}-1$, yields
\[
1
> \frac{0.7405d}{f_{k}^{2}f_{\ell}^{2} \left| N_{\alpha} \right|}
  \left( \frac{0.0002344 |g|^{4}\cN_{d',4}^{4}}{\left| N_{\alpha} \right|^{3}} dy_{k}^{2} \right)^{r_{0}-1}.
\]

Combined with $d>N_{\alpha}$, this implies that
\begin{equation}
\label{eq:step-iv}
1
> \frac{0.7405}{f_{k}^{2}f_{\ell}^{2}}
  \left( \frac{0.0002344 |g|^{4}\cN_{d',4}^{4}}{\left| N_{\alpha} \right|^{3}} dy_{k}^{2} \right)^{r_{0}-1}.
\end{equation}

We now proceed similarly to the way we did in Subsection~\ref{subsect:step-iii}.

Since $-N_{\alpha}$ is a square, we have $f_{k}=f_{\ell}=1$ from Proposition~\ref{prop:quad-rep}(b).
Also, as $r_{0}-1 \geq 1$, it follows that $0.7405 \cdot 0.0002344^{r_{0}-1}
> 0.0001735^{r_{0}-1}$. So the above inequality implies
\[
1>
\left( \frac{0.0001735 |g|^{4}\cN_{d',4}^{4}}{\left| N_{\alpha} \right|^{3}} dy_{k}^{2} \right)^{r_{0}-1}
> \left( y_{k}^{2}\frac{d \left( |g|\cN_{d',4} \right)^{4}}
    {5764\left| N_{\alpha} \right|^{3}} \right)^{r_{0}-1}.
\]

But if
\begin{equation}
\label{eq:ykUB-step4}
y_{k} \geq \frac{76\left| N_{\alpha} \right|^{3/2}}{\sqrt{d} \left( |g|\cN_{d',4} \right)^{2}},
\end{equation}
then the right-hand side is greater than $1$, so this is not possible.

\subsection{Small $d$}
\label{subsect:step-v}

To complete the proof of Proposition~\ref{prop:4.1}, we need to remove the
assumption on $d$ in \eqref{eq:d-LB}.

We check directly all pairs $(a,d)$ of positive integers with $2 \leq d \leq 104$
not a square and $-N_{\alpha}=d-a^{2}$ a square.
Note that there are only finitely many such pairs, since $N_{\alpha}=a^{2}-d<0$.

If $y_{k}=y^{2}$ is a square, then from \eqref{eq:yk-defn}, we know that
$x_{k}^{2}-dy^{4}=N_{\alpha}$. By Theorem~1.1 of \cite{Akh}, our theorem holds
for $N_{\alpha}=-1$, so we may also assume that $N_{\alpha} \leq -4$. There are
59 such pairs $(a,d)$.

For each such pair,
we solved $x^{2}-dy^{4}=N_{\alpha}$ using Magma (version V2.28-2) \cite{Magma}
and its \verb+IntegralQuarticPoints()+ function. For the 59 equations, this
calculation took 17.66 seconds using MAGMA's online calculator. 21 of the equations
had at least two solutions in positive integers. Only two of these had three
solutions in positive integers:\\
$x^{2}-17y^{4}=-16$ has the solutions $(x,y)=(1, 1),(16, 2)$, $(103, 5)$,\\
$x^{2}-68y^{4}=-64$ has the solutions $(x,y)=(2, 1),(32, 2)$, $(206, 5)$.\\
None of the equations had more solutions.

For $x^{2}-17y^{4}=-16$, since $\left( 103+5^{2}\sqrt{17} \right)/\left( 1+\sqrt{17} \right)
=\left( 161+39\sqrt{17} \right)/8$, the solutions $(103,5)$ and $(1,1)$
arise from different sequences. Hence Proposition~\ref{prop:4.1} holds for
$(a,d)=(1,17)$.
Similarly, Proposition~\ref{prop:4.1} holds for $(a,d)=(2,68)$ too.
This completes the proof of Proposition~\ref{prop:4.1}.

\section{Proof of Theorem~\ref{thm:1.3-seq-new}}
\label{sect:proofs}

If $N_{\alpha}$ is even, then $|g|\cN_{d',4} \geq 4$ by Lemma~\ref{lem:gNd4}.
So the right-hand side of the inequality in Proposition~\ref{prop:4.1} is
$\max \left( 1, 4.75 \left| N_{\alpha} \right|^{3/2}/\sqrt{d} \right)$.
But for $u \geq 5$ and $k \neq 0$, we have $y_{k} \geq 25\left| N_{\alpha} \right|/4$
by Lemma~\ref{lem:Y-LB}(c). So Theorem~\ref{thm:1.3-seq-new} holds for
$u \geq 5$ and we need only consider $1 \leq u \leq 4$.

Similarly, if $N_{\alpha}$ is odd and $u \geq 9$, then we have $y_{k} \geq 81\left| N_{\alpha} \right|/4
>\max \left( 1, 19\left| N_{\alpha} \right|^{3/2}/\sqrt{d} \right)$ for $k \neq 0$
and we need only consider $1 \leq u \leq 8$.

Next we treat the case when $K<-1$, where, as in Lemma~\ref{lem:Y-LB}, $K$
is the largest negative integer such that $y_{K}>1$.

\begin{lemma}
\label{lem:K-2-proof}
If $K<-1$, then there are at most two distinct integer squares among the $y_{k}$'s.
\end{lemma}

\begin{proof}
From Lemma~\ref{lem:K-value}, we know that $K<-1$ can only happen in the following
two cases:\\
{\rm (i)} $a \geq 1$, $d=a^{2}+4$, $t=a$ and $u=1$, where $\alpha=2\varepsilon$ and $N_{\alpha}=-4$, \\
{\rm (ii)} $a \geq 1$, $d=a^{2}+1$, $t=2a$ and $u=2$, where $\alpha=\varepsilon$ and $N_{\alpha}=-1$.

In case~(i), we have $x_{k}^{2}-dy_{k}^{2}=-4$.

If $d=a^{2}+4$ is even, then it is divisible by $4$ and $x_{k}$ is also even. So
the equation becomes $\left( x_{k}/2 \right)^{2}-(d/4)y_{k}^{2}=-1$.
From Theorem~1.1 of \cite{Akh}, there are at most two distinct squares among the
$y_{k}$'s here, so the lemma holds in this case.

If $d$ is odd, the conditions in Theorem~1 of \cite{L4} hold with $A=d=a^{2}+4$,
$B=1$ both odd and the minimal solution of equation~(1) there in odd positive
integers being $(1,a)$. Hence there are at most two distinct squares among the
$y_{k}$'s in this case too.

In case~(ii), we can also apply Theorem~1.1 of \cite{Akh}.
\end{proof}

For the remainder of this section, we may assume that $K=-1$.

\begin{lemma}
\label{lem:kpm1}
If $y_{\ell}>y_{k}>1$ are two squares, then $k= \pm 1$.
\end{lemma}

\begin{proof}
From Theorem~1.1 of \cite{Akh}, we cannot have three distinct squares among the
$y_{k}$'s for $N_{\alpha}=-1$, so we may assume $\left| N_{\alpha} \right| \geq 4$
for $N_{\alpha}$ even and $\left| N_{\alpha} \right| \geq 9$ for $N_{\alpha}$ odd.

We will suppose that $y_{k}$ is a square with $|k|>1$, but show that this
is not possible.

From Lemma~\ref{lem:Y-LB}(c) and $\left| N_{\alpha} \right|=d-a^{2}<d$,
we obtain $y_{k} \geq \left| N_{\alpha} \right|^{2}u^{4}/10$. So, for $u \geq 2$
and $\left| N_{\alpha} \right| \geq 4$,
we have $y_{k} \geq (8/5)\left| N_{\alpha} \right|^{2} \geq (32/5)\left| N_{\alpha} \right|
>(32/5)\left| N_{\alpha} \right|^{3/2}/\sqrt{d}$.

However, we saw at the start of this section that for $N_{\alpha}$ even,
there can be no squares $y_{k}$ and $y_{\ell}$ satisfying $y_{\ell}>y_{k}
>\max \left( 1, 4.75\left| N_{\alpha} \right|^{3/2}/\sqrt{d} \right)$.
So we need only consider
$u=1$ when $N_{\alpha}$ is even.

Similarly, if $-N_{\alpha}$ is an odd square with $N_{\alpha} \leq -9$ and
$u \geq 3$, then we have
$y_{k}>(81/10)\left| N_{\alpha} \right|^{2} \geq (729/10)\left| N_{\alpha} \right|
>(729/10)\left| N_{\alpha} \right|^{3/2}/\sqrt{d}$.
But, we saw at the start of this section that for $N_{\alpha}$ odd,
there can be no squares $y_{k}$ and $y_{\ell}$ satisfying $y_{\ell}>y_{k}
>\max \left( 1, 19\left| N_{\alpha} \right|^{3/2}/\sqrt{d} \right)$.
So we need only consider $u \leq 2$ when $N_{\alpha}$ is odd.

For $u=2$, we have $y_{k}>(8/5)\left| N_{\alpha} \right|^{2}$ for $|k| \geq 2$.
So if $\left| N_{\alpha} \right| \geq 25$ is odd, then
we have $y_{k} \geq (200/5)\left| N_{\alpha} \right|>40\left| N_{\alpha} \right|^{3/2}/\sqrt{d}$.
So, for odd $N_{\alpha}$, we need only consider $N_{\alpha}=-9$.
In this case, $N_{\varepsilon}=\left( t^{2}-4d \right)/4=\pm 1$, so
$t^{2}-4a^{2}=-4N_{\alpha} \pm 4=32,40$. This means that
$\left( t,a,d \right)=(6, 1,10)$ (recall from Subsection~\ref{subsect:notation}
that we only consider positive values of $t$ and $u$).
Here we have $y_{\pm 2} \geq y_{-2}=493 > 19.75\left| N_{\alpha} \right|^{3/2}/\sqrt{d}$.

For $u=1$, we proceed similarly. Here we have $y_{k}>\left| N_{\alpha} \right|^{2}/10$
for $|k| \geq 2$, so we need to consider $4 \leq \left| N_{\alpha} \right| \leq 36$
when $-N_{\alpha}$ is an even square and $9 \leq \left| N_{\alpha} \right| \leq 169$
when $-N_{\alpha}$ is an odd square.

For each value of $N_{\alpha}$, the norm of $\varepsilon$ leads us to an
equation of the form $t^{2}-a^{2}=-N_{\alpha} \pm 4$. We solve each of these
equations over the integers using PARI/GP and calculate $y_{k}$ and $y_{-k}$
for $k \geq 2$ until the lower bound from Proposition~\ref{prop:4.1} is exceeded
(we never had to go beyond $|k|=5$).
No squares were found.
\end{proof}

\begin{lemma}
\label{lem:yPM1}
If both $y_{\pm 1}$ are squares, then $y_{k}$ is not a square for any $|k|>1$.
\end{lemma}

\begin{proof}
From Lemmas~\ref{lem:gap}(a) and ~\ref{lem:Y-LB}(c), we have
\[
y_{1} > \frac{57.32d^{2}}{\left| N_{\alpha} \right|^{2}} y_{-1}^{3}
\geq \frac{57.32d^{2}}{\left| N_{\alpha} \right|^{2}} \left( \frac{\left| N_{\alpha} \right| u^{2}}{4} \right)^{3}
= \frac{57.32d^{2}\left| N_{\alpha} \right| u^{6}}{64}.
\]

By Proposition~\ref{prop:4.1} and Lemma~\ref{lem:gNd4}, if the quantity on the
right-hand side exceeds
$\max \left( 1, 19\left| N_{\alpha} \right|^{3/2}/\sqrt{d} \right)$, then the lemma
follows. Since $d>\left| N_{\alpha} \right|$, this holds if $d^{2}u^{6}>19 \cdot 64/57.32$,
which holds unless $d^{2}<21.3$ and $u=1$. This leaves $d=2$ or $3$, but for neither
of these does $u=1$ yield a unit.
\end{proof}

\begin{lemma}
\label{lem:u-cong}
Suppose that $-N_{\alpha}$ is an odd square and $y_{\pm 1}$ is a square integer.

\noindent
{\rm (a)} If $t^{2}-du^{2}=4$, then $t \equiv 2 \pmod{4}$ and $u \equiv 0 \pmod{4}$.

\noindent
{\rm (b)} If $t^{2}-du^{2}=-4$, then $u \equiv 2 \pmod{4}$.
\end{lemma}

\begin{proof}
For $t^{2}-du^{2}= \pm 4$, if $t$ is odd, then
$d \equiv u^{2} \equiv 1 \pmod{4}$. Since $N_{\alpha}$ is odd, it follows that
$a$ is even. Therefore $4y_{\pm 1}=2du^{2} \pm 2atu \pm 4 \equiv 2 \pmod{4}$,
which is not possible since $y_{\pm 1} \in \bbZ$.
Therefore, if $N_{\alpha}$ is odd, $t$ must be even.

(a) Next we show that if $t^{2}-du^{2}=4$, then $4\nmid t$. If $4 \mid t$, then $du^{2} \equiv 12 \pmod{16}$.
Since $-N_{\alpha}$ is an odd
square, we have $-N_{\alpha} \equiv 1,9 \pmod{16}$. If $a$ is odd, then $a^{2} \equiv 1 \pmod{8}$.
So $d \equiv 2 \pmod{8}$ and hence $u^{2} \equiv 6 \pmod{8}$, which is not possible.
If $a$ is even, then
$a^{2} \equiv 0,4 \pmod{16}$. Since $-N_{\alpha}=d-a^{2} \equiv 1,9 \pmod{16}$,
it follows that $d \equiv 1 \pmod{8}$. So $u^{2} \equiv 12 \pmod{16}$,
which is also impossible. So $4 \nmid t$.

Therefore, we must have $t \equiv 2 \pmod{4}$, if $t^{2}-du^{2}=4$ and $N_{\alpha}$
is an odd square.

From $t \equiv 2 \pmod{4}$ and $t^{2}-du^{2}=4$, we have $du^{2} \equiv 0 \pmod{16}$.

Since $N_{\alpha}$ is odd, we have $-N_{\alpha} \equiv 1,9 \pmod{16}$.
Also $a^{2} \equiv 0,1,4,9 \pmod{16}$. Hence $d \equiv 1,2,5,9,10,13 \pmod{16}$ and 
so $8|u^{2}$. Thus $4|u$.

\vspace*{1.0mm}

(b) We proceed similarly. Suppose that $4 \mid t$, then $du^{2} \equiv 4 \pmod{16}$.
As in the proof of part~(a), if $a$ is odd, then
$d \equiv 2 \pmod{8}$. Hence $u^{2} \equiv 2 \pmod{8}$, which is not possible.

If $a$ is even, as in the proof of part~(a),
then $d \equiv 1 \pmod{8}$. So $u^{2} \equiv 4 \pmod{16}$.
I.e., $u \equiv 2 \pmod{4}$.

If $t \equiv 2 \pmod{4}$, then $du^{2} \equiv 8 \pmod{16}$.
Again, if $a$ is odd, then $d \equiv 2 \pmod{8}$ and so $u^{2} \equiv 4 \pmod{8}$.

If $a$ is even, then $d \equiv 1 \pmod{8}$ and $u \equiv 8 \pmod{16}$,
which is not possible.
\end{proof}

\begin{lemma}
\label{lem:nrmA-odd}
If $-N_{\alpha}$ is an odd square, $1 \leq u \leq 8$ and $y_{\pm 1} \in \bbZ$ a
square with
$y_{\pm 1} < \max \left( 1, \frac{76\left| N_{\alpha} \right|^{3/2}}{\sqrt{d} \left( |g|\cN_{d',4} \right)^{2}} \right)$,
then we must have $u=2$, $t^{2}-du^{2}=-4$ and $\gcd \left( a^{2}, d \right)=1$.
\end{lemma}

\begin{proof}
We proceed in steps.

(1) We start by showing that $t^{2}-du^{2}=4$ is not possible.

From Lemma~\ref{lem:u-cong}(a), if $t^{2}-du^{2}=4$ and $-N_{\alpha}$ is an odd
square, then $4|u$, so we are left with $u=4$ or $8$ here.

(1-i) First, we eliminate $u=4$.

We can write
$N_{\alpha}=-(2n+1)^{2}$, so
$t^{2}-\left( a^{2}+(2n+1)^{2} \right)u^{2}=4$ becomes
$t^{2}-64n^{2}-64n=16a^{2}+20$.
But $t^{2}\equiv 20 \pmod{64}$, has no solution, so $a$ must be odd.
Hence $t^{2}-64n^{2}-64n=16a^{2}+20$ implies $t^{2} \equiv 36 \pmod{64}$.
So $t \equiv \pm 6 \pmod{16}$.

Expanding the expressions for $4y_{\pm 1}$, with $a=2a_{1}+1$, $N_{\alpha}=-(2n+1)^{2}$
and $t=16t_{1} \pm 6$, we find that $4y_{\pm 1} \equiv 20 \pmod{32}$. But this
congruence has no solution with $y_{\pm 1}$ a square.
Hence $u=4$ is not possible.

\vspace*{1.0mm}

(1-ii) Now we eliminate $u=8$.

Arguing modulo $9$, from $t^{2}-64d=4$, we see that
$d \equiv 0,3,5,6 \pmod{9}$.  From this, $-N_{\alpha}=d-a^{2}$ being a square and
the squares modulo $9$ being $0,1,4,7$, we must have $\left( a^{2} \bmod 9, d \bmod 9 \right)=(0,0)$, $(1,5)$, $(4,5)$ or $(7,5)$.

We will show that $d \equiv 5 \pmod{9}$ is not possible.
In this case, from $t^{2}-64d=4$, we find that $t^{2} \equiv 0 \pmod{9}$.
So $t \equiv 0 \pmod{3}$.
Since $d \equiv 5 \pmod{9}$, we have $y_{\pm 1}=1+32d \pm 4at \equiv 2 \pm 4at \pmod{3}$.
Since $3|t$, it follows that $y_{\pm 1} \equiv 2 \pmod{3}$, which can never be
square. Hence $9|\gcd \left( a^{2}, d \right)$.

(Since we will use this same argument in several places in the proof of this
lemma and the next, we wrote programs (in both Maple and PARI/GP) to automate these steps
and prove that $\gcd \left( a^{2}, d \right)$ has certain factors (typically,
powers of $2$ and $3$). In fact, using this same code, we could have eliminated
case~(1-i) above too.)

Therefore by Lemma~\ref{lem:gNd4}, we have $\left( |g|\cN_{d',4} \right)^{2} \geq 36$.
Hence
$\frac{76\left| N_{\alpha} \right|^{3/2}}{\sqrt{d} \left( |g|\cN_{d',4} \right)^{2}}
<(76/36)\left| N_{\alpha} \right|$. But by
Lemma~\ref{lem:Y-LB}(c), we know that $y_{\pm 1} \geq 16\left| N_{\alpha} \right|$ here.
So the case where $u=8$, $-N_{\alpha}$ is an odd square and $t^{2}-du^{2}=4$ is
excluded.

Thus for $N_{\alpha}$ odd, we do not need to consider $t^{2}-du^{2}=4$.

\vspace*{1.0mm}

(2) Now we consider $t^{2}-du^{2}=-4$.

(2-i) From Lemma~\ref{lem:u-cong}(b), we have $u \equiv 2 \pmod{4}$.
For $u=6$, $t^{2}-du^{2}=-4$ is not possible modulo $9$.
So we must have $u=2$ if $t^{2}-du^{2}=-4$.

\vspace*{1.0mm}

(2-ii) We show that $\gcd \left( a^{2}, d \right)=1$ when $u=2$.

We suppose otherwise.
Since $N_{\alpha}$ is odd, $\gcd \left( a^{2}, d \right)$ must be odd
and since $d= \left( t^{2}+4 \right)/4$,
any odd prime factor, $p$, of the gcd must also be a factor of $t^{2}+4$. This
means that $t^{2} \equiv -4 \pmod{p}$. Hence $p \equiv 1 \pmod{4}$. Since
$p^{2}|a^{2}$ and $-N_{\alpha}=d-a^{2}$ is a perfect square, this also implies that
$p^{2}|d$. Therefore by Lemma~\ref{lem:gNd4}, we have $\left( |g|\cN_{d',4} \right)^{2} \geq 100$.
Hence the theorem holds if $y_{\pm 1}>(76/100)\left| N_{\alpha} \right|$, by
Proposition~\ref{prop:4.1}. But by
Lemma~\ref{lem:Y-LB}(c), we know that $y_{\pm 1} \geq \left| N_{\alpha} \right|$ here.
So the theorem holds if $u=2$, $-N_{\alpha}$ is an odd square and
$\gcd \left( a^{2}, d \right)>1$.
\end{proof}

We now prove an analogous lemma for $N_{\alpha}$ even.

\begin{lemma}
\label{lem:nrmA-even}
If $-N_{\alpha}$ is an even square, $1 \leq u \leq 4$ and
$y_{\pm 1} < \max \left( 1, \frac{76\left| N_{\alpha} \right|^{3/2}}{\sqrt{d} \left( |g|\cN_{d',4} \right)^{2}} \right)$,
then we must have $u=1$, $t^{2}-du^{2}=-4$, $N_{\alpha} \equiv 12 \pmod{16}$ and
$\gcd \left( a^{2}, d \right)=1,4$.
\end{lemma}

\begin{proof}
As in the proof of Lemma~\ref{lem:nrmA-odd}, we proceed in steps.

(1) We start by showing that $t^{2}-du^{2}=4$ is not possible.

(1-i) We show that $u=4$, $-N_{\alpha}$ an even square and $t^{2}-du^{2}=4$ is
not possible.

Here we can argue as for $u=8$ in the proof of
Lemma~\ref{lem:nrmA-odd} and use the Maple program mentioned there modulo $9$ to
show that $9 | \gcd \left( a^{2}, d \right)$.
Therefore by Lemma~\ref{lem:gNd4},
we have $\left( |g|\cN_{d',4} \right)^{2} \geq 144$.
Hence
$\frac{76\left| N_{\alpha} \right|^{3/2}}{\sqrt{d} \left( |g|\cN_{d',4} \right)^{2}}
\leq (76/144)\left| N_{\alpha} \right|$.
But by
Lemma~\ref{lem:Y-LB}(c), we know that $y_{\pm 1} \geq 4\left| N_{\alpha} \right|$ here.
So we can exclude $u=4$ and $t^{2}-du^{2}=4$ from consideration.

\vspace*{1.0mm}

(1-ii) We show that $u=3$, $-N_{\alpha}$ an even square and $t^{2}-du^{2}=4$ is
not possible.

We will again argue as for $u=8$ in the proof of
Lemma~\ref{lem:nrmA-odd} and use the Maple program mentioned there modulo $32$ to
show that $64 | \gcd \left( a^{2}, d \right)$.

Therefore by Lemma~\ref{lem:gNd4}, we have $\left( |g|\cN_{d',4} \right)^{2} \geq 256$.
Hence the theorem holds if $y_{\pm 1}>(76/256)\left| N_{\alpha} \right|$, by
Proposition~\ref{prop:4.1}. But by
Lemma~\ref{lem:Y-LB}(c), we know that $y_{\pm 1} \geq (9/4)\left| N_{\alpha} \right|$ here.

\vspace*{1.0mm}

(1-iii) We show that $u=2$, $-N_{\alpha}$ an even square and $t^{2}-du^{2}=4$ is
not possible.

Using the argument for $u=8$ in the proof of
Lemma~\ref{lem:nrmA-odd} and the Maple program mentioned there modulo $16$,
we obtain $16 | \gcd \left( a^{2}, d \right)$.

Doing the same modulo $9$, we obtain $9 | \gcd \left( a^{2}, d \right)$.
Combining this with $16|\gcd \left( a^{2},d \right)$ and applying Lemma~\ref{lem:gNd4},
we have $|g|\cN_{d',4} \geq 24$.
Hence
$\frac{76\left| N_{\alpha} \right|^{3/2}}{\sqrt{d} \left( |g|\cN_{d',4} \right)^{2}}
\leq \left( 76/24^{2} \right)\left| N_{\alpha} \right|$.
But by
Lemma~\ref{lem:Y-LB}(c), we know that $y_{\pm 1} \geq \left| N_{\alpha} \right| u^{2}/4
= \left| N_{\alpha} \right|$ here.
So we can exclude $u=2$ and $t^{2}-du^{2}=4$ from consideration.

\vspace*{1.0mm}

(1-iv) We show that $u=1$, $-N_{\alpha}$ an even square and $t^{2}-du^{2}=4$ is
not possible.

Here too we use the argument and Maple program from the proof of Lemma~\ref{lem:nrmA-odd},
first modulo $9$ and then modulo $16$. From the latter, we obtain
$64|\gcd \left( a^{2}, d \right)$. So by Lemma~\ref{lem:gNd4},
$|g|\cN_{d',4} \geq 48$. Hence
$\frac{76\left| N_{\alpha} \right|^{3/2}}{\sqrt{d} \left( |g|\cN_{d',4} \right)^{2}}
\leq \left( 76/48^{2} \right)\left| N_{\alpha} \right|$.
But by
Lemma~\ref{lem:Y-LB}(c), we know that $y_{\pm 1} \geq \left| N_{\alpha} \right|/4$
here.
So we can exclude $u=1$ and $t^{2}-du^{2}=4$ from consideration.

\vspace*{1.0mm}

(2) We now consider $t^{2}-du^{2}=-4$.

Independent of the parity of $N_{\alpha}$, for $u=3$ and $4$, $t^{2}-du^{2}=-4$
is not possible modulo $9$ and $16$, respectively.
So we can only have $u=1,2$ with $t^{2}-du^{2}=-4$

(2-i) We show that if $u=1$, $-N_{\alpha}$ is an even square and $t^{2}-du^{2}=-4$,
then $N_{\alpha} \equiv 12 \pmod{16}$.

We find that $y_{\pm 1} \equiv 0,1,4,9 \pmod{16}$ is not possible if
$-N_{\alpha}=b_{1}^{2}$ where $b_{1} \equiv 4 \pmod{8}$.

If $8|b_{1}$, then
$a \equiv 2 \pmod{4}$ and $|g|\cN_{d',4} \geq 16$.
So
$\frac{76 \left| N_{\alpha} \right|^{3/2}}{\left( \sqrt{d} \left( |g| \cN_{d',4} \right)^{2} \right)}
\leq \left( 76/16^{2} \right) \left| N_{\alpha} \right|^{3/2}/\sqrt{d}$.
However $76/16^{2}$ is bigger than $1/4$, so we must work a bit harder to eliminate
$8|b_{1}$.

If $a^{2}\geq 0.292d$, then
\[
\frac{76 \left| N_{\alpha} \right|^{3/2}}{\left( \sqrt{d} \left( |g| \cN_{d',4} \right)^{2} \right)}
<\left( 76 \sqrt{1-0.292}/16^{2} \right) \left| N_{\alpha} \right|<0.2498 \left| N_{\alpha} \right|.
\]

By Lemma~\ref{lem:Y-LB}(c), we know that $y_{\pm 1} \geq \left| N_{\alpha} \right|/4$.
So $a^{2} \geq 0.292d$ is excluded.

Suppose that $d \geq 80$. Then $t^{2}=d-4$, so $t^{2} \geq (76/80)d$.
If $a^{2}<0.292d$, then $t-a>\sqrt{76d/80}-\sqrt{0.292d}$ and so
$(t-a)^{2}/4>0.047d>0.047 \left| N_{\alpha} \right|$.
We can write $4y_{\pm 1} = \left( t \pm au \right)^{2}-N_{\alpha}u^{2}$.
So here with $u=1$, we have
$y_{\pm 1} \geq \left( t-a \right)^{2}/4+ \left| N_{\alpha} \right|/4>0.297\left| N_{\alpha} \right|$.
But $76 \left| N_{\alpha} \right|/16^{2}=0.296875\left| N_{\alpha} \right|$.
So we can exclude $8|b_{1}$, provided $d \geq 80$.

For $d<80$ with $u=1$ and
$-N_{a}$ an even square divisble by $64$, there is just one possibility:
$d=68$, $t=8$, $u=1$, $a=2$, so $N_{\alpha}=-64$ and $y_{-1}=25$.
But here, 
$\frac{76 \left| N_{\alpha} \right|^{3/2}}{\left( \sqrt{d} \left( |g| \cN_{d',4} \right)^{2} \right)}
<\left( 76/16^{2} \right)\left| N_{\alpha} \right|^{3/2}/\sqrt{d}
<19$, so this case is excluded too and hence we can exclude $8|b_{1}$ altogether.

Hence $b_{1} \equiv 2 \pmod{4}$ and so $-N_{\alpha} \equiv 4 \pmod{16}$.

\vspace*{1.0mm}

(2-ii) We show that if $u=1$, $-N_{\alpha}$ is an even square, $t^{2}-du^{2}=-4$,
and $N_{\alpha} \equiv 12 \pmod{16}$, then $\gcd \left( a^{2}, d \right)=1$ or $4$.

We have $v_{2} \left( \gcd \left( a^{2}, d \right) \right)=0$ or $2$.
Since $-N_{\alpha} \equiv 4 \pmod{16}$, we have $v_{2} \left( b' \right)=2$ in the
first case and $0$ in the second case.

We consider first $\gcd \left( a^{2}, d \right)>4$ odd.
So by Lemma~\ref{lem:gNd4} and the argument in the previous paragraph,
$|g|\cN_{d',4}=4\sqrt{\gcd \left( a^{2}, d \right)}$. Hence
$\frac{76\left| N_{\alpha} \right|^{3/2}}{\sqrt{d} \left( |g|\cN_{d',4} \right)^{2}}
=\frac{76\left| N_{\alpha} \right|^{3/2}}{16\sqrt{d}\gcd \left( a^{2}, d \right)}$.
If $\gcd \left( a^{2}, d \right) \geq 25$, then the right-hand side is less
than $0.19\left| N_{\alpha} \right|<\left| N_{\alpha} \right|/4$, the right-hand
side being the lower bound for $y_{\pm 1}$ from Lemma~\ref{lem:Y-LB}(c).

Since $d=t^{2}+4$, arguing modulo $3$, we see that $3 \nmid d$ and hence
$\gcd \left( a^{2}, d \right)=9$ is not possible.

Now consider $\gcd \left( a^{2}, d \right)>4$ even. Then
$v_{2} \left( \gcd \left( a^{2}, d \right) \right)=2$ (so $\gcd \left( a^{2}, d \right)
\equiv 4 \pmod{16}$).
So by Lemma~\ref{lem:gNd4} and the argument above,
$|g|\cN_{d',4}=2\sqrt{\gcd \left( a^{2}, d \right)}$. Hence
$\frac{76\left| N_{\alpha} \right|^{3/2}}{\sqrt{d} \left( |g|\cN_{d',4} \right)^{2}}
=\frac{76\left| N_{\alpha} \right|^{3/2}}{4\sqrt{d}\gcd \left( a^{2}, d \right)}$.

We saw above that $3 \nmid d$, so we can ignore $\gcd \left( a^{2}, d \right)=6^{2}$.
If $\gcd \left( a^{2}, d \right) \geq 100$, then
$\frac{76\left| N_{\alpha} \right|^{3/2}}{4\sqrt{d}\gcd \left( a^{2}, d \right)}
\leq 0.19\left| N_{\alpha} \right|
<\left| N_{\alpha} \right|/4$, the lower bound for $y_{\pm 1}$ from
Lemma~\ref{lem:Y-LB}(c).

\vspace*{1.0mm}

(2-iii) We show that $u=2$, $-N_{\alpha}$ an even square and $t^{2}-du^{2}=-4$
is not possible.

Arguing modulo $64$, we have $d \equiv a^{2} \equiv 1 \pmod{16}$ and $8|t$.
So $16|N_{\alpha}$ and, in the notation of Lemma~\ref{lem:gNd4}, $16|b'$.
Thus $|g|\cN_{d',4}=8\sqrt{\gcd \left( a^{2}, d \right)}$, where
$\gcd \left( a^{2}, d \right)$ is odd.

It follows that
$\frac{76 \left| N_{\alpha} \right|^{3/2}}{\left( \sqrt{d} \left( |g| \cN_{d',4} \right)^{2} \right)}
\leq \frac{76 \left| N_{\alpha} \right|^{3/2}}{64\left( \sqrt{d} \gcd \left( a^{2}, d \right) \right)}$.
However $76/64$ is bigger than $1$, so we proceed as in case~(2-i) above.

If $a^{2}\geq 0.292d$, then
\[
\frac{76 \left| N_{\alpha} \right|^{3/2}}{\left( 64 \sqrt{d} \gcd \left( a^{2}, d \right) \right)}
<\frac{76 \sqrt{1-0.292}}{64} \left| N_{\alpha} \right|<0.9992 \left| N_{\alpha} \right|.
\]

By Lemma~\ref{lem:Y-LB}(c), we know that $y_{\pm 1} \geq \left| N_{\alpha} \right|$.
So $a^{2} \geq 0.292d$ is excluded.

Suppose that $d \geq 2$. Then $t^{2}=4d-4$, so $t^{2} \geq 2d$.
If $a^{2}<0.292d$, then $t-a>\sqrt{2d}-\sqrt{0.292d}$ and so
$(t-a)^{2}/4>0.19d>0.19 \left| N_{\alpha} \right|$.
We can write $4y_{\pm 1} = \left( t \pm au \right)^{2}-N_{\alpha}u^{2}$.
So here with $u=2$, we have
$y_{\pm 1} \geq \left( t-a \right)^{2}/4+ \left| N_{\alpha} \right|>1.19\left| N_{\alpha} \right|$.
But $76 \left| N_{\alpha} \right|/64=1.1875\left| N_{\alpha} \right|$.
So we can exclude the case where $u=2$, $-N_{\alpha}$ an even square and $t^{2}-du^{2}=-4$.
\end{proof}

Theorem~\ref{thm:1.3-seq-new}(c) now follows from the bounds on $u$ that we obtained
from Proposition~\ref{prop:4.1} at the start of Section~\ref{sect:proofs},
along with Lemmas~\ref{lem:K-2-proof}, \ref{lem:kpm1}, \ref{lem:nrmA-odd}
and \ref{lem:nrmA-even}. Furthermore, from Lemma~\ref{lem:yPM1}, if $y_{-1}$ and
$y_{1}$ are both squares, then there are no further squares. So we may assume
that precisely one of $y_{\pm 1}$ is a square. From Lemma~\ref{lem:Y-LB}(c), we
obtain $y_{k}>19 \left| N_{\alpha} \right|^{3/2}/\sqrt{d}$ for $|k| \geq 2$,
provided that $d>190$ when $u=1$ and that $d>11$ when $u=2$. This leaves $26$ values
of $d$. For each of these, we compute $y_{-2}$ directly for each possibility
of $N_{\alpha}$ and compare the value to the bound in Proposition~\ref{prop:4.1}.
Where the bound was exceeded, Proposition~\ref{prop:4.1} tells us that there are
no further squares, completing the proof for them.
There were only four cases where $y_{-2}$ did not exceed that bound:
$(a,b,d,t,u)=(1,1,2,2,2)$, $(1,1,5,1,1)$, $(2,1,8,2,1)$, $(3,1,13,3,1)$.
These were treated in Subsection~\ref{subsect:step-v}.

\section{Examples}
\label{sect:egs}

In this section, we give examples showing that our conjectures and results are
best possible.

\subsection{Examples for Conjecture~\ref{conj:1-seq}}
\label{subsect:conj1-egs}

\begin{center}
\begin{table}[h]
\begin{tabular}{rrll}\hline
$a$        & $b$      & indices, $k$      & $\sqrt{y_{k}}$\\ \hline
         1 &      $3$ & [0,  -1,  3,  -5] & [3, 5, 31, 167]\\
      1019 &     $27$ & [0,   1, -3,  -7] & [27, 65, 29, 983]\\
       167 &     $13$ & [0,   1, -3,   4] & [13,29,71, 407]\\
       157 &     $29$ & [0,  -1,  3,  -4] & [29, 47, 307, 649]\\
         1 &     $41$ & [0,  -1, -9,  11] & [41, 71, 80753, 470861]\\
      1633 &     $65$ & [0,  -1, -4,   7] & [65, 97, 1331, 24791]\\
     48479 &    $211$ & [0,  -3,  4,  -7] & [211, 1007, 6743, 34205]\\
     45649 &    $677$ & [0,  -1, -4,   4] & [677, 1133, 15679, 16825]\\
   1940147 &   $1217$ & [0,  -3,  4, -11] & [1217, 3289, 40573, 3794239]\\
    600589 &   $2213$ & [0,  -1, -4,   4] & [2213, 3673, 50801, 55415]\\
  20509501 &   $8689$ & [0,  -1,  3,  -4] & [8689, 13619, 94393, 187603]\\
 255488029 &  $13457$ & [0,  -1,  3,  -4] & [13457, 5683, 189241, 15821]\\
 409660129 &  $17023$ & [0,  -1, -4,  -8] & [17023, 7073, 7949, 269495]\\
3032771269 &  $46313$ & [0,  -1, -4,  -8] & [46313, 19213, 15269, 516625]\\ \hline
\end{tabular}
\caption{Examples for $d=2$ with $t=u=2$}
\label{table:d2}
\end{table}
\end{center}

In addition to the examples for $d=2$ in Table~\ref{table:d2}, we also found examples
with four squares for $(d,t,u,a,b)=(3,4,2,672,91)$,
$(6,10,4,78,7)$, $(6,10,4,34986,149)$,\\
$(6,10,4,3663828,2257)$, $(30,22,4,826320, 1111)$ and $(37,12,2,138,5)$.

At least for $d=2$ and $d=6$, it appears there may be infinitely many such examples.

\subsection{Examples for Conjecture~\ref{conj:3-seq}}
\label{subsect:conj3-egs}

Let $n \geq 5$ be an odd integer. Put $a=\left( n^{2}-9 \right)/4$ and
$d=\left( n^{4}-2n^{2}+17 \right)/16$, so that $N_{\alpha}=4-n^{2}$. With $t=\left( n^{2}-1 \right)/2$ and $u=2$,
$\varepsilon=\left( t+u\sqrt{d} \right)/2$ is a unit in
$\cO_{\bbQ \left( \sqrt{d} \right)}$. We have $y_{1}=\left( n^{2}-3 \right)^{2}/4$
and $y_{-1}=n^{2}$.
This example shows that Conjecture~\ref{conj:3-seq} is best possible, if true.

Let $n>5$ satisfy $n \equiv 1 \pmod{4}$ and not divisible by $5$. Put
$a=\left( n-5 \right)/4$, $d=\left( n^{2}+6n+25 \right)/16$. Here
$N_{\alpha}=-n$.
With $t=2a+4$ and $u=2$, $\varepsilon=\left( t+u\sqrt{d} \right)/2$
is a unit in $\cO_{\bbQ \left( \sqrt{d} \right)}$ and
we have $y_{1}=(n+1)^{2}/4$.
If $n$ is a square, then this example shows that Conjecture~\ref{conj:3-seq} is
best possible when $\left| N_{\alpha} \right|$ is a perfect square.

\subsection{Examples for Conjecture~\ref{conj:2-seq}}
\label{subsect:conj2-egs}

The examples in Subsection~\ref{subsect:conj1-egs} also apply in this more
general case to show there can be at least four distinct squares.

Here are two examples showing that there are at least three distinct squares when
$\left| N_{\alpha} \right|$ is a square:\\
$(d,t,u,a,b)=(5,1,1,43,3)$: $N_{\alpha}=38^{2}$, $y_{-3}=5^{2}$, $y_{0}=3^{2}$, $y_{11}=53^{2}$;\\
$(d,t,u,a,b)=(10,6,2,1,1)$: $N_{\alpha}=-3^{2}$, $y_{0}=1^{2}$, $y_{1}=2^{2}$, $y_{2}=5^{2}$.

\vspace*{1.0mm}

Here are two examples showing that the same is true when
$\left| N_{\alpha} \right|$ is a prime power:\\
$(d,t,u,a,b)=(5,1,1,153,4)$: $N_{\alpha}=22129$, $y_{-3}=11^{2}$, $y_{0}=4^{2}$, $y_{30}=8862^{2}$;\\
$(d,t,u,a,b)=(51,100,14,2,1)$: $N_{\alpha}=-47$, $y_{-1}=6^{2}$, $y_{0}=1^{2}$, $y_{1}=8^{2}$.

\vspace*{1.0mm}

Here is an example showing that if $\left| N_{\alpha} \right|$ is not
a prime power, but $\core \left( \left| N_{\alpha} \right| \right)$ is a prime,
then there can be distinct four squares (unlike in Conjecture~\ref{conj:1-seq}):\\
$(d,t,u,a,b)=(5,1,1,7,1)$: $N_{\alpha}=2^{2}\cdot 11$, $y_{-9}=9^{2}$, $y_{0}=1$, $y_{1}=2^{2}$, $y_{3}=3^{2}$.

And finally, an example showing that if $\core \left( \left| N_{\alpha} \right| \right)$
is twice an odd prime, then there can be distinct four squares (again unlike in Conjecture~\ref{conj:1-seq}):\\
$(d,t,u,a,b)=(6,10,4,2,3)$: $N_{\alpha}=-2\cdot 241$, $y_{-3}=63^{2}$, $y_{0}=3^{2}$, $y_{1}=7^{2}$, $y_{3}=69^{2}$.

\subsection{Examples for Theorem~\ref{thm:1.3-seq-new}(a)}
\label{subsect:thm14a-egs}

From \eqref{eq:rel-sqrs} with $b=u=1$, we find that $(2t+a)^{2}-8y_{1}=a^{2}-8$.
Since $-N_{\alpha}=d-a^{2}=t^{2}+4-a^{2}$ is a square, we put $a=2$, so if $y_{1}$
is a square, then $2t+2$ is an element of the sequence $\left( t_{n} \right)_{n \geq 0}$
with $t_{0}=6$, $t_{1}=238$ and $t_{n+1}=34u_{n}-u_{n-1}$. The corresponding values
of $\sqrt{y_{1}}$ are the elements of the sequence $\left( u_{n} \right)_{n \geq 0}$
with $u_{0}=5$, $u_{1}=169$ and $u_{n+1}=34u_{n}-u_{n-1}$. Then $y_{1}$ is approximately
$t^{2}/2 \approx \left| N_{\alpha} \right|/2 \approx \left| N_{\alpha} \right|^{3/2}/ \left( 2\sqrt{d} \right)$.

In fact, such $y_{1}$'s are the ones that are squares satisfying the conditions in Theorem~\ref{thm:1.3-seq-new}(a)
with the smallest ratio compared to $\left| N_{\alpha} \right|^{3/2}/\sqrt{d}$,
the quantity in Proposition~\ref{prop:4.1}.

\subsection{Examples for Lemma~\ref{lem:Y-LB}}
\label{subsect:YLB-egs}

We can show that the lower bound in Lemma~\ref{lem:Y-LB}(c) is actually best
possible. Let $n$ be a positive integer such that $a=2n^{2}-3$ is not divisible
by $5$ and put $d=4n^{4}-8n^{2}+8$. We have $N_{\alpha}=1-4n^{2}$. Then $(t,u)=(a+1,1)$,
$y_{-1}=n^{2}$ and so
$y_{-1}/ \left| N_{\alpha} \right| \rightarrow 1/4$ from above as $n \rightarrow +\infty$.
There exist families with the same limit for $y_{-1}/\left| N_{\alpha} \right|$
when $-N_{\alpha}$ is a square too.

The sequences given in Lemma~\ref{lem:K-value} show that for $k=-1$, the lower
bound in Lemma~\ref{lem:Y-LB}(c) can be attained too. At least, when the lower
bound is $1$.

\section{Corrections from published version}

The most significant one is the change from 20 August 2025 for Lemma~3.10.

\vspace*{1.0mm}

\noindent
2 Aug 2024:\\
removed the PV footnote that should not be present.

\vspace*{1.0mm}

\noindent
25 Aug 2024:\\
in Lemma~3.9, added comment at the end that $\zeta_{4}$ is any primitive $4$-th root of unity.

\vspace*{1.0mm}

\noindent
1 Oct 2024:\\
in the statement of Lemma~3.9(b) added the condition that $k,\ell \neq 0$.\\
This is needed for use of Prop~3.1.

\vspace*{1.0mm}

\noindent
8 Dec 2024:\\
used $0.1263$ and $3.96$ consistently throughout paper\\
had some $0.127$ and $3.959$ and some $0.1263$ and $3.96$ before\\
(see 3 May 2025 change below too).

\vspace*{1.0mm}

\noindent
20 Jan 2025:\\
towards bottom of page 299 of published file, had $d'=u_{2}^{2}t'$, but should be
$d'=u_{2}^{2}t'/g^{2}$.\\
arxiv location: page~7, line~-8.

\vspace*{1.0mm}

\noindent
3 Feb 2025:\\
displayed formula on line~2 of the proof of Lemma~5.3:\\
$y_{-1}$ on the right-hand side of first inequality should be $y_{-1}^{3}$.
(but $y_{-1}^{3}$ correctly used afterwards)

\vspace*{1.0mm}

\noindent
27 Feb 2025:\\
Near beginning of Subsection, ``Construction of Approximations''\\
change ``. We let $u=\ldots$'' to ``such that $u=\ldots$''.

\vspace*{1.0mm}

\noindent
23 Apr 2025:\\
added $N_{\alpha}<0$ to proof of Lemma~3.5(c) when $k=1$.

\vspace*{1.0mm}

\noindent
23 Apr 2025:\\
--a few lines after that, $u_{1}=u$ should be $u_{1}=u/2$\\
--also removed redundant ()'s in expressions for $y_{\pm 1}$ in equation~(1.2)\\
--journal name for reference [17] corrected.

\vspace*{1.0mm}

\noindent
3 May 2025:\\
third line of Section~4.2:\\
``From the equality in (4.2)...'' should read ``From (4.2)...''

Also reverted the numerical changes from 8 Dec 2024 and only changed $0.127$ to $0.1263$
in equation~(3.31) so that the argument there to use (3.22) (where we need $0.1263$)
would be correct. No other changes are needed.

\vspace*{1.0mm}

\noindent
5 May 2025:\\
where appropriate (i.e., where they arise from the Representation Proposition),
all $\omega$'s, $\varepsilon$'s and $\varphi$'s changed to
$\omega_{k}$'s, $\varepsilon^{k}$'s and $\varphi_{k}$'s, respectively.

\newpage

\noindent
22 May 2025:\\
page 297, line 3:\\
in the definition of $Y_{m,n,r}$ in Section~2, I added the argument, $z$, that was missing.\\
I.e., changed ``$Y_{m,n,r}$'' to ``$Y_{m,n,r}(z)$''.\\
arxiv location: page~5, line~14.

\vspace*{1.0mm}

\noindent
22 June 2025:\\
page~306, after line~4. State what $g_{1}'^{2}$ and $g_{1}''^{2}$ are.\\
I did not previously state this explicitly before.\\
arxiv location: page~12, line~-6.

\vspace*{1.0mm}

\noindent
3 July 2025:\\
for max, use $\left( \right)$ rather than $\left\{ \right\}$\\
E.g., in gap principle (Lemma~3.8)\\
referee suggestion for the $-N_{\alpha}$ square for any $b$ paper (the ``...$y_0=b^2$ (I)'' paper).

\vspace*{1.0mm}

\noindent
31 July 2025:\\
page~330, line~-7: remove period before ``to show that''.\\
arxiv location: page~33, sentence after equation~(4.12).

\vspace*{1.0mm}

\noindent
10 Aug 2025:
page~300, last line in Section~2: give this an equation number~(2.10).\\
arxiv location: page~8, last line in Section~2.

\vspace*{1.0mm}

\noindent
20 Aug 2025:\\
Lemma~3.10 requires the condition that $\gcd(a,b)$ is odd.\\
The existing proof is fine with that condition added.


\begin{thebibliography}{10}
\bibitem{Akh}
S. Akhtari,
\emph{The Diophantine equation $aX^{4}-bY^{2} = 1$},
Journal f\"{u}r die reine und angewandte Mathematik {\bf 630} (2009), 33--57.

\bibitem{Baker1}
A. Baker,
\emph{Rational approximations to certain algebraic numbers},
Proc. London. Math. Soc. (3) {\bf 14} (1964), 385--398.

\bibitem{Baker2}
A. Baker,
\emph{Rational approximations to $\sqrt[3]{2}$ and other algebraic numbers},
Quart. J. Math. Oxford {\bf 15} (1964), 375--383. 

\bibitem{Magma}
W.~Bosma, J.~Cannon, C.~Playoust,
\emph{The Magma algebra system. I. The user language},
J.~Symbolic Comput. {\bf 24} (1997), 235--265. 

\bibitem{BMS1}
Y. Bugeaud, M. Mignotte, S. Siksek,
\emph{Classical and modular approaches to exponential Diophantine equations I. Fibonacci and Lucas perfect powers},
Ann. Math. {\bf 163} (2006), 969--1018.

\bibitem{Chen1}
Chen Jian Hua, P. Voutier,
\emph{Complete solution of the diophantine equation $X^{2}+1=dY^{4}$ and a related
family of quartic Thue equations},
J. Number Theory {\bf 62} (1997), 71--99.

\bibitem{Evert1}
J.-H. Evertse,
\emph{On the Representation of Integers by Binary Cubic Forms of Positive Discriminant},
Invent. Math. {\bf 73} (1983), 117--138.

\bibitem{Liou}
J. Liouville,
\emph{Sur des classes tr\`{e}s-\'{e}tendues de quantit\'{e}s dont la valeur
n'est ni alg\'{e}brique, ni m\^{e}me r\'{e}ductible \`{a} des irrationnelles
alg\'{e}briques},
C. R. Acad. Sci. Paris, S\'{e}r. {\bf A 18} (1844) 883--885.

\bibitem{L4}
W. Ljunggren,
\emph{On the Diophantine equation $Ax^{4}-By^{2}=C$ ($C=1,4$)},
Math. Scand. {\bf 21} (1967) 149--158.

\bibitem{Pari}
The PARI~Group, PARI/GP version {\tt 2.12.0}, Univ. Bordeaux, 2019,
\url{http://pari.math.u-bordeaux.fr/}.

\bibitem{Sch}
W. M. Schmidt,
\emph{The zero multiplicity of linear recurrence sequences},
Acta Math. {\bf 182} (1999), 243--282.

\bibitem{Siegel1}
C. L. Siegel,
\emph{\"{U}ber einige Anwendungen diophantischer Approximationen},
Abh, Preuss. Akad. Wiss. {\bf 1} (1929), 41--69.

\bibitem{Siegel2}
C. L. Siegel,
\emph{Die Gleichung $ax^{n}-by^{n}=c$},
Math. Ann. {\bf 114} (1937), 57--68.

\bibitem{St}
C. L. Stewart,
\emph{On divisors of Lucas and Lehmer numbers},
Acta Math. {\bf 211} (2013), 291--314.

\bibitem{V2}
P. M. Voutier,
\emph{Thue's Fundamentaltheorem, I: The General Case},
Acta Arith. {\bf 143} (2010), 101--144.

\bibitem{V3}
P. M. Voutier,
\emph{Thue's Fundamentaltheorem, II: Further Refinements and Examples},
J. Number Theory {\bf 160} (2016), 215--236.

\bibitem{V4}
P. M. Voutier,
\emph{Improved constants for effective irrationality measures from hypergeometric functions},
Combinatorics and Number Theory {\bf 11} (2022), 161--180
\url{https://doi.org/10.2140/moscow.2022.11.161}.

\bibitem{V6}
P. M. Voutier,
\emph{Sharp bounds on the number of squares in recurrence sequences and solutions of $X^{2}-\left( a^{2}+b \right) Y^{4}=-b$},
Research Number Theory (accepted) \url{https://arxiv.org/abs/1807.04116}.
\end{thebibliography}
\end{document}